# Modularity of the Rankin-Selberg $L$-series, and multiplicity one for $\mathrm{SL}(2)$

By Dinakar Ramakrishnan

*In memory of my father Sundaram Ramakrishnan* (SRK)

## Contents

## 1. Introduction

Let $f, g$ be primitive cusp forms, holomorphic or otherwise, on the upper half-plane $\mathcal{H}$ of levels $N, M$ respectively, with (unitarily normalized) $L$-functions

$$L(s, f) = \sum_{n \geq 1} \frac{a_n}{n^s} = \prod_{p\,\mathrm{prime}} [(1 - \alpha_p p^{-s})(1 - \beta_p p^{-s})]^{-1}$$

and

$$L(s, g) = \sum_{n \geq 1} \frac{b_n}{n^s} = \prod_{p\,\mathrm{prime}} [(1 - \alpha'_p p^{-s})(1 - \beta'_p p^{-s})]^{-1}.$$

When $p$ does not divide $N$ (resp. $M$), the inverse roots $\alpha_p, \beta_p$ (resp. $\alpha'_p, \beta'_p$) are nonzero with sum $a_p$ (resp. $b_p$). For every $p$ prime to $NM$, set

$$L_p(s, f \times g) = [(1 - \alpha_p \alpha'_p p^{-s})(1 - \alpha_p \beta'_p p^{-s})(1 - \beta_p \alpha'_p p^{-s})(1 - \beta_p \beta'_p p^{-s})]^{-1}.$$

Let $L^*(s, f \times g)$ denote the (incomplete Euler) product of $L_p(s, f \times g)$ over all $p$ not dividing $NM$. This is closely related to the convolution $L$-series $\sum_{n \geq 1} a_n b_n n^{-s}$, whose miraculous properties were first studied by Rankin and Selberg.

A fundamental question, first raised by Langlands, is to know whether this Rankin-Selberg product is modular, i.e., if there exists an automorphic form $f \boxtimes g$ on $\mathrm{GL}(4)/\mathbb{Q}$ whose standard $L$-function equals $L^*(s, f \times g)$ after removing the ramified and archimedean factors. The first main result of this paper is to answer it in the affirmative, in fact with the base field $\mathbb{Q}$ replaced by any number field $F$ (see Theorem M, §3). Our proof uses a mixture of converse theorems, base change and descent, and it also appeals to the local regularity properties of Eisenstein series and the scalar products of their truncations. Briefly, experts have long suspected that this result should follow from the converse theorem for $\mathrm{GL}(4)$ requiring only twists by $\mathrm{GL}(2)$ and $\mathrm{GL}(1)$, which has recently been published by Cogdell and Piatetski-Shapiro ([CoPS]). While this is morally true, *three difficulties* crop up when one tries to implement this principle, and new ideas are required to surmount the difficulties. The *first difficulty* arises because, given cuspidal automorphic representations $\pi, \pi'$ on $\mathrm{GL}(2)/F$, one needs to have a definition of a candidate for $\pi \boxtimes \pi'$ as an *admissible* representation of $\mathrm{GL}(4, \mathbb{A}_F)$ (which is needed before we can test its analytic properties to show modularity). If one admits the local Langlands correspondence, then the local candidates $\pi_v \boxtimes \pi'_v$ can be defined as the admissible representations of $\mathrm{GL}(4, F_v)$ corresponding to $\sigma_v \otimes \sigma'_v$, where $\sigma_v$ (resp. $\sigma'_v$) is the 2-dimensional representation of the Weil-Deligne group $W'_{F_v}$ associated to $\pi_v$ (resp. $\pi'_v$) by Kutzko ([Ku]). (Since the work on this paper was completed, two preprints, one by M. Harris and R. L. Taylor, and the other by G. Henniart, have appeared, establishing the local correspondence for $\mathrm{GL}(n)$. But we feel that it will be satisfying not to have to appeal to it here; global arguments should always be able to circumvent fine local difficulties.) We get around this problem by appealing to the *base change* results of Arthur and Clozel ([AC]). To be precise, we first make use of the fact that every supercuspidal representation becomes a principal series representation after a finite normal solvable base change, so that a candidate for $\pi \boxtimes \pi'$ can be defined over suitable (infinite families of) global, solvable extensions $K/F$. Then, after proving modularity upstairs, we perform a simultaneous descent via an inductive argument in cyclic layers. At the end, we get as a byproduct the definition of $\pi_v \boxtimes \pi'_v$ at *any* $v$ (and any $F$). Once one has the admissible $\Pi$ on $\mathrm{GL}(4)/K$, one needs, for modularity via the converse theorem, good analytic information on $L(s, \Pi \times \eta)$, for any cuspidal $\eta$ on $\mathrm{GL}(m)/K$, for $m \leq 2$. The $m = 1$ case is easy by the Rankin-Selberg theory as extended to general $K$ by Jacquet ([J]). But the $m = 2$ case is subtle and leads to problems. One of them, which is the *second difficulty*, is caused by not knowing the equality, at *every* place, with each other and with $L(s, \Pi \times \eta)$, of the *three candidates* for the *triple product L-function* $L(s, \pi \times \pi' \times \eta)$.

The first is defined formally as an Euler product of $\{L(s, \sigma_v \otimes \sigma'_v \otimes \tau_v)\}$, where $\tau_v$ is the 2-dimensional of $W'_{F_v}$ associated to $\eta_v$; the second candidate is defined by the integral representation of Garrett ([G]), as generalized by Piatetski-Shapiro and Rallis ([PS-R2]); and the third is given by the machinery of Langlands-Shahidi ([La2], [Sh1]). Though the unramified local factors are known to be the same in all cases, something close to an equality is essential as no candidate has all the desired properties, and they have complementary strengths. (One simply cannot avoid dealing with the bad places!)

By a careful analysis and synthesis of known results due to Ikeda and others, and by using in addition some global arguments involving the works of Langlands and Tunnell on the Artin conjecture, we manage to prove at the end that the *first* candidate has *all* the desired properties under some local restrictions which can be achieved by solvable base change. So base change is used for yet another reason! (At the archimedean places, we make use of the work of Ikeda ([Ik2]) on the computation of the triple product $L$-factors for unramified representations.) We prove later the equality, at *each* place $v$, of $L(s, \pi_v \times \pi'_v \times \eta_v)$ (resp. $\varepsilon(s, \pi_v \times \pi'_v \times \eta_v)$) with the $\mathrm{GL}(4) \times \mathrm{GL}(2)$-factor $L(s, (\pi_v \boxtimes \pi'_v) \times \tau_v)$ (resp. $\varepsilon(s, (\pi_v \boxtimes \pi'_v) \times \eta_v)$).

The *third* and final *difficulty* is the question of boundedness in vertical strips of the triple product $L$-function. This is absolutely crucial for applying the converse theorem, and this does not seem to be elementary as, for example, in the case of the standard $L$-function of $\mathrm{GL}(2)$. Our method is to use Arthur's truncation ([A1]) of the noncuspidal Eisenstein series $E(g, s)$ on $\mathrm{GSp}(6)/F$ which occurs in the integral representation of $L(s, \pi \times \pi' \times \eta)$, and reduce the problem, via local regularity of eigenfunctions of the Laplacian, to the norms of the truncated Eisenstein series, making use of the fact that the intervening constant terms are in our case expressible in terms of abelian $L$-functions. The reader is referred to Section 3.1 for a fuller discussion of the strategy of proof of the various parts. The details of proof occupy the Sections 3.2–3.7.

As a consequence of Theorem M, we settle a conjecture of Labesse and Langlands ([LL], [Lab1]) asserting that the space of cusp forms on $\mathrm{SL}(2)$ has multiplicity one (see Theorem 4.1.1). To see what this means concretely, consider $f, g$ as above with trivial characters such that

$$a_p^2 = b_p^2$$

for almost all $p$. Then multiplicity one (in this context) implies that there exists a quadratic Dirichlet character $\chi$ such that (for almost all $p$)

$$a_p = b_p \chi(p).$$

(In fact, if $N, M$ are in addition square-free, $\chi$ must be trivial.)

When we started on this project in 1994, we were able to settle quickly the case when $f$ and $g$ are holomorphic, but then learned that this case had been known to various people including D. Blasius and J-P. Serre; here the idea is to make use of the associated $\ell$-adic Galois representations (see [K, pp. 90–91], and [Ra2], where there is also a mod $\ell$ analog and a density result). But this method *does not* work for Maass forms, and the starting point for this paper was our realization in fall 1994 that both cases could be tackled simultaneously if one knew of the existence of $f \boxtimes g$. We first managed to prove multiplicity one for $\mathrm{SL}(2)$ over $\mathbb{Q}$ in fall 96 by some special tricks and a weaker form of Theorem M. The proof given here works over arbitrary $F$, but is shorter, partly because some of the earlier arguments over $\mathbb{Q}$ have been transplanted to the proof of Theorem M. If there is any creativity in this paper, it is perhaps foremost in the application of Theorem M to this problem, though it is not the most technically difficult part. It should also be remarked that it has been expected for some time now that the Labesse-Langlands conjecture should follow from the adjoint square lifting from $\mathrm{SL}(2)$ to $\mathrm{PGL}(3)$, more precisely from a careful comparison of the stable trace formula for $\mathrm{SL}(2)$ and the twisted trace formula for $\mathrm{PGL}(3)$ relative to $g \to {}^t g^{-1}$. This approach has been expounded by Y. Flicker in a series of papers culminating in [F]. But we ignore the question of whether this program has been completed, as our approach is *totally* different and (hopefully) has independent interest. We would also like to note in passing that, as shown by D. Blasius ([B]), multiplicity one fails for $\mathrm{SL}(n)$, for every $n > 2$.

A corollary of Theorem M is the deduction of the standard analytic properties of 4-fold convolutions $L(s, f_1 \times f_2 \times f_3 \times f_4)$ for quadruples $(f_1, f_2, f_3, f_4)$ of primitive cusp forms. This is done by appealing, in addition, to the higher Rankin-Selberg theory for $\mathrm{GL}(4) \times \mathrm{GL}(4)$ due to Jacquet, Piatetski-Shapiro and Shalika ([JPSS1]), and Shahidi ([Sh 1,3,5]). We show, moreover, that $L(s, \mathrm{sym}^2(f_1) \times f_2 \times f_3)$ and $L(s, \mathrm{sym}^3(f_1) \times f_2)$ have meromorphic continuations and functional equations (see Theorem 4.2.1) and are nonzero in $\{\Re(s) \geq 1\}$ without any pole there except possibly at $s = 1$.

By using the existence of $\boxtimes$ on $\mathrm{GL}(2) \times \mathrm{GL}(2)$ in conjunction with a global argument we prove an identity, even at the ramified places $v$, equating the Rankin-Selberg $L$ and $\varepsilon$-factors of $\pi_{1,v} \boxtimes \pi_{2,v}, \pi_{3,v} \boxtimes \pi_{4,v}$ with those of $\sigma_{1,v} \otimes \sigma_{2,v} \otimes \sigma_{3,v} \otimes \sigma_{4,v}$, where $\sigma_{j,v}$ denotes, for each $j \leq 4$, the 2-dimensional representation of $W'_{F_v}$ associated to $\pi_{j,v}$ (cf. [Ku]). This gives a little bit of information (see Prop. 4.3.3) on the local Langlands correspondence for $\mathrm{GL}(4)$. Moreover, we use the identities and the method of [PR] to prove that when $\pi_j$ has trivial central character for every $j$, the *global root number* $W((\pi_1 \boxtimes \pi_2) \times (\pi_3 \boxtimes \pi_4))$ is 1 (see Theorem 4.3.4).

We then turn the tables and prove, in Section 4.4, that as a *consequence* of Theorem M, the three candidates for the triple product $L$-functions on $\mathrm{GL}(2)$ all agree at *all* the places. (Recall that we needed to know that they were very nearly the same to get Theorem M in the first place.) This may be of independent interest.

The final application of our main result is the proof of the Tate conjecture for 4-fold products $V$ of modular curves, asserting in particular that the order of pole at $s = 2$ of the $L$-function over any solvable (normal) number field $K$ of the Galois module $W_\ell := H^4_{\text{et}}(V_{\overline{\mathbb{Q}}}, \mathbb{Q}_\ell)$ equals the rank of the group of $K$-rational codimension 2 Tate cycles on $V_{\overline{\mathbb{Q}}}$ (see Theorem 4.5.1). Moreover we show, in line with the works of Ribet ([Ri1]) and V. K. Murty ([Mu]) on the *Jacobians* of modular curves, that the latter number can also be computed with the *Tate cycles* replaced by the *algebraic cycles* modulo homological equivalence if the level of at least one of the curves is square-free. We refer to Chapter 5 for a precise statement.

We would like to express our gratitude to Ilya Piatetski-Shapiro for his continued interest in this project, and for kindly writing down, with J. Cogdell, the form of the converse theorem for GL(4) which we need ([CoPS]). Thanks are also due to T. Ikeda for writing down his calculations of the archimedean factors of the triple product $L$-functions ([Ik2]), to S. Rallis for useful remarks on these $L$-functions, to F. Shahidi for explaining his approach to the same via Langlands's theory of Eisenstein series and for commenting on an earlier version, to my colleague T. Wolff for helpful conversations on an analytic lemma we use in Section 3.4, and to many others, including H. Jacquet, R. P. Langlands, J. Rogawski and P. Sarnak, who have shown encouragement and interest. Special thanks must go to J. Cogdell for reading the earlier and the revised versions thoroughly and making crucial remarks. Part of the technical typing of this paper was done by Cherie Galvez, whom we thank. Finally, we would like to express our appreciation to the following: the National Science Foundation for support through the grants DMS-9501151 and DMS-9801328, Université Paris-sud, Orsay, where we spent a fruitful month during September 1996, the DePrima Mathematics House in Sea Ranch, CA, for inviting us to visit and work there during August 1996 and 1998, MATSCIENCE, India, for hospitality in February 98, and — last, but not the least — the MSRI, Berkeley, for (twice) providing the right climate to work in; this project was started (in 1994) and essentially ended there.

## 2. Notation and preliminaries

2.1. Let $\overline{\mathbb{Q}}$ denote the algebraic closure of $\mathbb{Q}$ in $\mathbb{C}$. For any subfield $K$ of $\overline{\mathbb{Q}}$, let $\text{Gal}(\overline{\mathbb{Q}}/K)$ denote the Galois group of $\overline{\mathbb{Q}}$ over $K$, together with the profinite topology. For any number field $F$ with ring of integers $\mathfrak{O}_F$, let $\Sigma(F)$ (resp. $\Sigma_\infty(F)$, resp. $\Sigma_0(F)$) denote the places (resp. archimedean places, resp. finite places) of $F$. At each $v$, let $F_v$ denote the completion of $F$ relative to $v$, and if moreover $v$ is finite, let $\mathfrak{O}_v$, $\mathfrak{P}_v$, $\mathbb{F}_v$, $|.|_v$ and $Nv = q_v$ respectively denote the ring of integers of $F_v$, the maximal ideal, the residue field, the normalized absolute value and the norm of $v$.

When $F = \mathbb{Q}$, there is a unique archimedean place $\infty$ given by the canonical imbedding of $\mathbb{Q}$ in $\mathbb{R}$, and every finite place $v$ corresponds to a rational prime number $p$, in which case we will write $\mathbb{Z}_p$ instead of $\mathcal{O}_v$. Denote by $\mathbb{A}_F$ the ring of adeles of $F$, which is the restricted direct product of $\{F_v\}$ relative to $\{\mathcal{O}_v\}$, equipped with its usual locally compact topology. Let $I_F$ denote the group of ideles of $F$, and $C_F$ the idele class group $I_F/F^*$.

We will fix a nontrivial unitary character $\psi$ of $\mathbb{A}_F$ which is trivial on $F$, and let $\psi_v$ be the $v$-component of $\psi$. Various quantities, such as the $\varepsilon$-factors, will depend on this choice, which we will suppress in our notation.

In this paper we will systematically use the powerful language of automorphic representations, though at relevant places we will indicate briefly the essence of what we do in the classical language. Given any primitive cusp form $f$ on the upper half-plane of level $N$ and character $\omega$, there exists a (unique) cuspidal automorphic representation $\pi = \pi(f)$ of $\text{GL}(2, \mathbb{A}_\mathbb{Q})$ of conductor $N$ and central character $\omega$ such that, at every prime $p$, the $p$-Euler factor of $L(s, f)$ agrees with the $L$-factor $L(s, \pi_p)$ of the $p$-component $\pi_p$. We refer the reader to the expository monograph [Ge] to understand how to go back and forth between the two approaches.

2.2. For every algebraic group $G$ over $F$, let $G(\mathbb{A}_F)$ denote the restricted direct product $\prod'_v G(F_v)$, endowed with the usual locally compact topology. By a *cuspidal* representation of $G(\mathbb{A}_F) = G(F_\infty) \times G(\mathbb{A}_{F,0})$, we will always mean an irreducible, unitary, cuspidal automorphic representation. Such a representation is in particular a restricted tensor product $\pi = \otimes'_v \pi_v = \pi_\infty \otimes \pi_0$, where each $\pi_v$ is an (irreducible) admissible representation of $G(F_v)$ for $v$ finite, and an admissible (Lie $G_v, K_v$)-module for $v$ archimedean, with $K_v$ denoting a maximal compact modulo center subgroup of $G(F_v)$; $\pi_0$ (resp. $\pi_\infty$) is the restricted tensor product of $\pi_v$ over all finite (resp. archimedean) places $v$. By definition, $\pi_v$ must be unramified at almost all $v$.

For any irreducible, automorphic representation $\pi$ of $\text{GL}(n, \mathbb{A}_F)$, let $L(s, \pi) = L(s, \pi_\infty) L(s, \pi_0)$ denote the associated *standard* $L$-function ([J]) of $\pi$; it has an Euler product expansion

(2.2.1) $$L(s, \pi) = \prod_v L(s, \pi_v),$$

convergent in a right half-plane. If $v$ is an archimedean place, then one knows (cf. [La1]) how to associate a semisimple $n$-dimensional $\mathbb{C}$-representation $\sigma(\pi_v)$ of the Weil group $W_{F_v}$, and $L(\pi_v, s)$ identifies with $L(\sigma_v, s)$. On the other hand, if $v$ is a finite place where $\pi_v$ is unramified, there is a corresponding semisimple

(Langlands) conjugacy class $A_v(\pi)$ in $\operatorname{GL}(n, \mathbb{C})$ such that

$$(2.2.2) \qquad L(s, \pi_v) = \det(1 - A_v T)^{-1}|_{T = q_v^{-s}}.$$

We may find a diagonal representative $\operatorname{diag}(\alpha_{1,v}(\pi), \ldots, \alpha_{n,v}(\pi))$ for $A_v(\pi)$, which is unique up to permutation of the diagonal entries. Let $[\alpha_{1,v}(\pi), \ldots, \alpha_{n,v}(\pi)]$ denote the resulting unordered $n$-tuple. Since $W_{F,v}^{\operatorname{ab}} \simeq F_v^*$, $A_v(\pi)$ clearly defines an abelian $n$-dimensional representation $\sigma(\pi_v)$ of $W_{F,v}$.

THEOREM 2.2.3 ([JS1]). *Let $n$ be $\geq 1$, and $\pi$ be a nontrivial cuspidal representation of $\operatorname{GL}(n, \mathbb{A}_F)$. Then $L(s, \pi)$ is entire. Moreover, for any finite set $S$ of places of $F$, the incomplete L-function $L^S(s, \pi) = \prod_{v \notin S} L(s, \pi_v)$ is holomorphic in $\Re(s) > 0$.*

When $n = 1$ such a $\pi$ is simply a unitary idele class character, and the result is due to Hecke. Also, when $\pi$ is trivial, $L(s, \pi) = \zeta_F(s)$.

By the theory of Eisenstein series, one has a sum operation $\boxplus$ ([La3]), which results in the following:

THEOREM 2.2.4 ([JS2]). *Given any $m$-tuple of cuspidal representations $\pi_1, \ldots, \pi_m$ of $\operatorname{GL}(n_1, \mathbb{A}_F), \ldots, \operatorname{GL}(n_m, \mathbb{A}_F)$ respectively, there exists a unitary, irreducible automorphic representation $\pi_1 \boxplus \cdots \boxplus \pi_m$ of $\operatorname{GL}n, \mathbb{A}_F)$, $n = n_1 + \cdots + n_m$, which is unique up to equivalence, such that for any finite set $S$ of places,*

$$L^S(s, \boxplus_{j=1}^m \pi_j) = \prod_{j=1}^m L^S(s, \pi_j).$$

Call such a (Langlands) sum $\pi \simeq \boxplus_{j=1}^m \pi_j$, with each $\pi_j$ cuspidal, an *isobaric* representation. Denote by $\operatorname{ram}(\pi)$ the finite set of finite places where $\pi$ is ramified, and let $\mathfrak{N}(\pi)$ be its conductor ([JPSS2]).

For every integer $n \geq 1$, set:

$$(2.2.5) \qquad \mathcal{A}(n, F) = \{\pi : \text{isobaric representation of } \operatorname{GL}(n, \mathbb{A}_F)\}/\simeq,$$

and

$$\mathcal{A}_0(n, F) = \{\pi \in \mathcal{A}(n, F) | \pi \text{ cuspidal}\}.$$

Put $\mathcal{A}(F) = \cup_{n \geq 1} \mathcal{A}(n, F)$ and $\mathcal{A}_0(F) = \cup_{n \geq 1} \mathcal{A}_0(n, F)$.

*Remark.* One can also define the analogs of $\mathcal{A}(n, F)$ for local fields $F$, where the "cuspidal" subset $\mathcal{A}_0(n, F)$ consists of essentially square-integrable representations of $\operatorname{GL}(n, F)$. See [La3] and [Ra1] for details.

Let $\pi, \pi'$ be isobaric automorphic representations in $\mathcal{A}(n, F), \mathcal{A}(n', F)$ respectively. Then there exist an associated Euler product $L(s, \pi \times \pi')$ ([JS2,4], [JPSS], [COPS2], [Sh1,3]), which converges in $\{\Re(s) > 1\}$, and admits a meromorphic continuation to the whole $s$-plane with a functional equation. When $v$ is archimedean or a finite place outside ram$(\pi)$, one has

$$(2.2.6) \qquad L_v(s, \pi \times \pi') \;=\; L(s, \sigma(\pi_v) \otimes \sigma(\pi'_v)).$$

When $n = 1$, $L(s, \pi \times \pi') = L(s, \pi\pi')$, and when $n = 2$ and $F = \mathbb{Q}$, this function is the Rankin-Selberg $L$-function, extended to arbitrary global fields by Jacquet ([J2]).

THEOREM 2.2.7 (Jacquet-Shalika [JS2]). *Let $\pi \in \mathfrak{A}_0(n, F)$, $\pi' \in \mathfrak{A}_0(n', F)$, and $S$ be a finite set of places. Then $L^S(s, \pi \times \pi')$ is holomorphic in $\{\Re(s) > 0\}$ unless $\pi \simeq \pi'^{\vee}$, in which case it has a unique pole (of order 1) at $s = 1$.*

We will also need the general construction of automorphic $L$-functions attached to any $\pi \in \mathcal{A}(n, F)$. Let $S$ be any finite set of places containing the archimedean and ramified places (for $\pi$). If $r$ is any algebraic representation of $\mathrm{GL}(n, \mathbb{C})$, which is the connected component of the $L$-group of $\mathrm{GL}(n)$, we put

$$L^S(s, \pi, r) \;=\; \prod_{v \notin S} L(s, \pi, r),$$

where

$$L(s, \pi, r) \;=\; \det(I - (Nv)^{-s} r(A_v(\pi)))^{-1}.$$

Of particular interest is when $r$ is the exterior square $\Lambda^2$, or the symmetric $k^{\text{th}}$ power $\mathrm{Sym}^k$, for some $k \geq 1$, of the standard representation. When $r$ is simply the standard representation, $L^S(s, \pi, r)$ is evidently just $L^S(s, \pi)$.

One knows (cf. [JS1], [BuF], [Sh1]) that $L^S(s, \pi, \Lambda^2)$ converges absolutely in $\Re(s) > 1$ and admits a meromorphic continuation with a functional equation of the usual type relating $s$ and $1 - s$. One also knows, by [BuG], [Sh1], the analogous properties of $L^S(s, \pi, \mathrm{sym}^2)$. For $n = 2$, one knows more by [GJ], namely that there is an (isobaric) automorphic representation $\mathrm{sym}^2(\pi)$ such that

$$L^S(s, \mathrm{sym}^2(\pi)) \;=\; L^S(s, \pi, \mathrm{sym}^2).$$

For $n = 2$, one also knows the meromorphic continuation and functional equation (cf. [Sh2]) of $L^S(s, \pi, \mathrm{sym}^k)$ for all $k \leq 5$.

2.3. *Base change and automorphic induction.* We will now review the results of Arthur and Clozel [AC] in a form which is suitable for the applications found in this article. See also [Lab2].

In [AC, Chap. (§§3–6)], one finds a construction, for any cyclic extension $K/F$ with $[K : F]$ a prime $\ell$, of maps

$$b_{K/F} : \mathcal{A}(n, F) \to \mathcal{A}(n, K), \; \pi \to \pi_K \quad \text{(base change)}$$

and
$$I_{K/F} : \mathcal{A}(n, K) \to \mathcal{A}(n\ell, F), \pi \to I(\pi) \quad \text{(automorphic induction)},$$

such that at every place $v$ of $F$ which is finite and unramified for the representations and $K/F$ (or archimedean) and a place $w$ of $K$ above $v$, we have (respectively)
$$\operatorname{res}^{F_v}_{K_w}(\sigma(\pi_v)) \simeq \sigma((\pi_K)_w),$$

and
$$\sigma(I(\pi)_v) \simeq \operatorname{ind}^{F_v}_{K_w}(\pi_w).$$

Given any local field $E$, if $\beta$ is an unramified representation of $\operatorname{GL}(r, E)$ (or if $E$ is archimedean), we write $\sigma(\beta)$ to signify the associated $r$-dimensional representation of $W'_E$ (or $W_E$).

There are also local analogs of these maps, customarily called "local base change" and "local automorphic induction". The former is constructed and discussed in great detail in Chapter 1 of [AC]. The later is discussed briefly in [C$\ell$1] and the relevant assertions are consequences of the results of [AC]. An alternate construction of the local automorphic induction for essentially tempered representations, which works also in characteristic $p$, is given in the paper [HH] of Henniart and Herb. The properties we will need are summarized as follows:

PROPOSITION 2.3.1. *Let $K/F$ be a cyclic extension of number fields or local fields of degree $\ell$, a prime. Let $\theta$ be a generator of $\operatorname{Gal}(K/F)$, and let $\chi$ denote the character of the idele class group (resp. multiplicative group) of $F$ in the global (resp. local) case associated to $K$. Then*

1. *The image of $b_{K/F}$ consists precisely of those $\beta \in \mathcal{A}(n, K)$ such that $\beta \simeq \beta \circ \theta$.*

2. *The image of $I_{K/F}$ consists precisely of those $\pi \in \mathcal{A}(n\ell, F)$ such that $\pi \simeq \pi \otimes \chi$.*

3. *For every $\pi \in \mathcal{A}(m, F)$ and $\beta \in \mathcal{A}(n, K)$, we have the* adjointness *property:*
$$L(s, \pi \times I(\beta)) = L(s, \pi_K \times \beta),$$

   *and*
$$\varepsilon(s, \pi \times I(\beta))\varepsilon(s, 1_K)^{nm} = \varepsilon(s, \pi_K \times \beta) \prod_{j=0}^{\ell-1} \varepsilon(s, \chi^j)^{nm}.$$

4. *Suppose $\beta$ is cuspidal (resp. supercuspidal) in $\mathcal{A}(n, K)$, for $K$ global (resp. local). Then*
$$I(\beta)_K \simeq \boxplus_{j=0}^{\ell-1} \beta \circ \theta^j.$$

*Moreover, $I(\beta)$ is cuspidal (resp. supercuspidal) if and only if $\beta$ is* not *isomorphic to $\beta \circ \theta$.*

5. *Suppose $\pi$ is cuspidal (resp. supercuspidal) in $\mathcal{A}(n, F)$, for $F$ global (resp. local). Then*
$$I(\pi_K) \simeq \boxplus_{j=0}^{\ell-1} \pi \otimes \chi^j.$$

*Moreover, $\pi_K$ is cuspidal (resp. supercuspidal) if and only if $\pi$ is* not *isomorphic to $\pi \otimes \chi$.*

For a proof of 1, 2 and the second half of 5, see [AC], [Cl1] and [HH]. For the remaining assertions, see [PR, pp. 7–8], where this proposition is stated as Proposition 3.1, and where $\mathcal{A}(n, F)$ is denoted $\text{Isob}(n, F)$.

## 3. Construction of $\boxtimes : \mathcal{A}(\text{GL}(2)) \times \mathcal{A}(\text{GL}(2)) \to \mathcal{A}(\text{GL}(4))$

The object of this chapter is to prove the following:

THEOREM M. *Let $\pi, \pi'$ be in $\mathcal{A}(2, F)$. Then:*
*Existence.  There exists an isobaric automorphic representation $\pi \boxtimes \pi'$ of $\text{GL}(4, \mathbb{A}_F)$ satisfying (at every finite place $v$)*

$(L_v)$ $\qquad\qquad L(s, (\pi \boxtimes \pi')_v) = L(s, \pi_v \times \pi'_v),$

*and*

$(\varepsilon_v)$ $\qquad\qquad \varepsilon(s, (\pi \boxtimes \pi')_v) = \varepsilon(s, \pi_v \times \pi'_v).$

*Also,*

$(L_\infty)$ $\qquad\qquad L(s, (\pi \boxtimes \pi')_\infty) = L(s, \pi_\infty \times \pi'_\infty).$

*Cuspidality Criterion. Suppose $\pi, \pi'$ are both cuspidal. If neither of them is associated to a character of a quadratic extension, then $\pi \boxtimes \pi'$ is cuspidal if and only if the following holds*:

(C) *$\pi'$ is not equivalent to $\pi \otimes \chi$, for any idele class character $\chi$ of $F$.*

*If $\pi' = I_K^F(\mu)$, for a character $\mu$ of a quadratic extension $K$, then $\pi \boxtimes \pi'$ is cuspidal if and only if the base change $\pi_K$ is cuspidal and not isomorphic to $\pi_K \otimes (\mu \circ \theta)\mu^{-1}$, where $\theta$ denotes the nontrivial automorphism of $K/F$.*

We will call an isobaric representation $\pi \boxtimes \pi'$ satisfying $(L_v)$ and $(\varepsilon_v)$, for all $v$, a *strong lifting* to $\text{GL}(4)/F$, or the *automorphic tensor product*, of the pair $(\pi, \pi')$. If it only satisfies $(L_v)$ for almost all $v$, and not necessarily $(\varepsilon_v)$ or $(L_\infty)$, we will call it a *weak lifting*.

We will briefly explain the arithmetic motivation for the cuspidality criterion above. Suppose $\pi, \pi'$ correspond to irreducible, continuous 2-dimensional representations $\sigma, \sigma'$ over $\mathbb{C}$ of $\mathrm{Gal}(\overline{\mathbb{Q}}/F)$. Then, since the operation $\boxtimes$ on $\mathcal{A}(2, F) \times \mathcal{A}(2, F)$ is supposed to correspond to the tensor product on the Galois side, $\pi \boxtimes \pi'$ should be cuspidal if and only if $\sigma \otimes \sigma'$ is irreducible. When $\sigma, \sigma'$ are nondihedral, reducibility happens if and only if one is a twist of the other by a character. That the criterion should be the same for all nondihedral cuspidal pairs $(\pi, \pi')$ is motivated by the hope (see [La3], [Cl2], and [Ra1]) that there is a group $\mathcal{L}_\mathcal{F}$, whose connected component $\mathcal{L}'_\mathcal{F}$ is pro-reductive with $\mathcal{L}_\mathcal{F}/\mathcal{L}'_\mathcal{F} \simeq \mathrm{Gal}(\overline{\mathbb{Q}}/F)$, and whose irreducible $n$-dimensional representations parametrize cuspidal automorphic representations of $\mathrm{GL}(n, \mathbb{A}_F)$.

3.1. *Relevant objects and the strategy.* We begin with a simple statement:

LEMMA 3.1.1. *Theorem* M *holds in the following three special cases*:

(I) *At least one of* $\{\pi, \pi'\}$ *is not cuspidal.*

(II) *At least one of* $\{\pi, \pi'\}$ *is automorphically induced by a character* $\mu$ *of* (*the idele class group of*) *a quadratic extension* $K$ *of* $F$.

(III) $\pi'$ *is a twist of* $\pi$; *i.e., there exists a character* $\chi$ *of* $C_F$ *such that*
$$\pi' \simeq \pi \otimes \chi.$$

*Proof.* Suppose $\pi'$ is not cuspidal. Then as we have seen in Section 2, there exist idele class characters $\mu_1, \mu_2$ of $F$ such that $\pi' = \mu_1 \boxplus \mu_2$. We set
$$\pi \boxtimes \pi' = (\pi \boxtimes \mu_1) \boxplus (\pi \boxtimes \mu_2),$$
which is not cuspidal. In this case one knows ([JS2]) that the identities $(L_v)$, $(\varepsilon_v)$, and $(L_\infty)$ are all satisfied everywhere. Thus we are done in case (I).

Suppose $\pi'$ is of the form $I_K^F(\mu)$, for a character $\mu$ of $C_K$, for some $K$. Then we set
$$\pi \boxtimes \pi' = I_K^F(\pi_K \otimes \mu),$$
where $\pi_k$ denotes the base change of $\pi$ to $K$. The identity $(L_v)$ at the unramified places $v$ is a direct consequence of the identities of Proposition 2.3.1. The fact that $(L_v)$ holds at every place, and that $(L_\infty)$ and $(\varepsilon_v)$ also hold, will be shown as part of a general assertion in the next section (see Prop. 3.2.1).

It is left to consider (III). In this case, we set
$$\pi \boxtimes \pi' = (\mathrm{sym}^2(\pi) \otimes \chi) \boxplus \omega\chi,$$
where $\omega$ is the central character of $\pi$. The asserted identities are then consequences of the work of Gelbart and Jacquet ([GJ]) and the strong multiplicity one theorem ([JS2]). □

We will say that $(\pi, \pi')$ is of *general type* if we are not in either of these three special cases.

Given any pair $(\Pi, \eta)$ of irreducible, generic, admissible representations of $\mathrm{GL}(n, \mathbb{A}_F)$, $\mathrm{GL}(m, \mathbb{A}_F)$ respectively, which are not necessarily automorphic, we may set

$$(3.1.2) \qquad L(s, \Pi \times \eta) = \prod_v L(s, \Pi_v \times \eta_v)$$

and

$$\varepsilon(s, \Pi \times \eta) = \prod_v \varepsilon(s, \Pi_v \times \eta_v).$$

The following result is a crucial ingredient of our approach.

THEOREM 3.1.3 (Cogdell-Piatetski-Shapiro ([CoPS])). *Let $T$ be a fixed finite set of finite places of $F$. Let $\Pi$ be an irreducible unitary, admissible, generic representation of $\mathrm{GL}(4, \mathbb{A}_F)$ which satisfies the following*:

*For every $\eta \in \mathcal{A}_0(n, F)$, $n \leq 2$, with $\eta_v$ unramified at every $v$ in $T$, we have*:

(MC) $\quad L(s, \Pi \times \eta)$ and $L(s, \Pi^\vee \times \eta^\vee)$ *converge absolutely in large $\Re(s)$,*

*and they admit meromorphic continuations to the whole s-plane.*

(E) $\qquad L(s, \Pi \times \eta)$ and $L(s, \Pi^\vee \times \eta^\vee)$ *are entire.*

*There is a functional equation*

(FE) $\qquad L(1 - s, \Pi^\vee \times \eta^\vee) = \varepsilon(s, \Pi \times \eta) L(s, \Pi \times \eta).$

(BV) $\qquad L(s, \Pi \times \eta)$ *is bounded in vertical strips.*

*Then $\Pi$ is nearly automorphic; i.e., there exists an automorphic representation $\Pi_1$ of $\mathrm{GL}(4, \mathbb{A}_F)$ such that $\Pi_v \simeq \Pi_{1,v}$ for almost all $v$.*

Consequently, Theorem M can be attacked for $(\pi, \pi')$ of general type if we can solve the following two problems (for a finite set $T$ of finite places):

(P1) Define an irreducible, generic admissible representation $\Pi$ of $\mathrm{GL}(4, \mathbb{A}_F)$ such that *at every place $v$*,

$$(3.1.4) \qquad L(s, \Pi_v) = L(s, \pi_v \times \pi'_v)$$

and

$$\varepsilon(s, \Pi_v) = \varepsilon(s, \pi_v \times \pi'_v).$$

(P2) Given any $\eta \in \mathcal{A}_0(m, F)$, $m \leq 2$, with $\eta_v$ unramified at every $v$ in $T$, the hypotheses (MC), (E), (FE) and (BV) hold for $L(s, \Pi \times \eta)$ and $L(s, \Pi^\vee \times \eta^\vee)$.

(P1) cannot be solved as posed without knowing the local Langlands correspondence for GL(4), which in any case we would like to avoid using. We can, however, define a candidate $\Pi_v$ at almost all places, in particular at the archimedean places and at finite places $v$ where $\pi_v$ and $\pi'_v$ are not both supercuspidal. (Even the supercuspidal cases can be dealt with easily in odd residual characteristics.) See Section 3.7. Given any $\pi$ we can find, using [Ku], a chain $F = K_0 \subset K_1 \subset \cdots \subset K_m = K$ with each $K_j/K_{j-1}$ cyclic of prime degree, in fact an infinite family $X = X(\pi)$ of such chains, such that the base change $\pi_K$ of $\pi$ to $K$, which exists by [AC], has the property that at *any* place $w$ of $K$, $\pi_{K,w}$ is not supercuspidal. Hence (P1) can be solved over $K$, and if we can also solve (P2) over a "sufficiently large" subset $Y$ of $X$, we can construct $\pi \boxtimes \pi'$ over each $K$ in $Y$, and then try to find a common descent to $F$ satisfying the desired properties (listed in Theorem M). The part dealing with descent follows an approach taken for GL(2) in our earlier joint work [BR] with D. Blasius. (There is a mistake in [BR], but it does not affect the relevant section there discussing the descent criterion; in any case, we give a complete argument here for GL($n$) (see Prop. 3.6.1).) The best possible situation will be one in which almost every finite place $v$ of $F$ splits completely in some $K \in Y$. There are complications in finding a common descent if $(\pi_K, \pi'_K)$ is not of general type for too many $K$. For this and other reasons, it will be best for us to use an inductive argument in cyclic layers of prime degree. We refer to Sections 3.6 and 3.7 for details.

After replacing $F$ by a suitable larger field, let us suppose that we have solved (P1), and that $(\pi, \pi')$ is of general type with $\pi_v$ not supercuspidal. Fix any $T, \eta$ as above. If $m = 1$, then all the hypotheses of Theorem 3.1.3 are just known assertions of the Rankin-Selberg theory. So we may assume that $m = 2$. (P2) can be solved if we can find a "triple product $L$-function" $L(s, \pi \times \pi' \times \eta)$ which satisfies all these hypotheses and satisfies (at every $v$)

(3.1.5) $$L(s, \pi_v \times \pi'_v \times \eta_v) = L(s, \Pi_v \times \eta_v).$$

There are in fact three candidates for such an $L$-function, but unfortunately, none of them possesses *a priori all* the desired properties. The first one is in some sense the most natural one. To define it, let $\sigma_v$ (resp. $\sigma'_v$) denote, at each $v$, the 2-dimensional representation of $W'_{F_v}$ associated to $\pi_v$ (resp. $\pi'_v$) by the local Langlands correspondence for GL(2) ([Ku]). (When $v$ is archimedean, $W'_{F_v}$ signifies just usual the Weil group.) Put

(3.1.6) $$L(s, \pi \times \pi' \times \eta) = \prod_v L(s, \sigma_v \otimes \sigma'_v \otimes \tau_v),$$

where $\tau_v$ denotes, at each $v$, the representation of $W'_{F_v}$ attached to $\eta_v$. We do not know *a priori* any of the desired properties for this Euler product.

One defines the epsilon factor $\varepsilon(\pi \times \pi' \times \eta)$ as the product over all $v$ of $\varepsilon(s, \sigma_v \otimes \sigma'_v \otimes \tau_v)$; these local factors will also be denoted $\varepsilon(s, \pi_v \times \pi'_v \times \eta_v)$.

The second candidate, which we denote by $L_1(s, \pi \times \pi' \times \eta)$, is the one defined by Piatetski-Shapiro and Rallis in [PSR2], generalizing the construction of Garrett ([G]) over $\mathbb{Q}$ when all three representations correspond to holomorphic modular forms (see also Ikeda ([Ik1]). This is defined by using an integral representation, and is known to satisfy (MC) and (FE). It also satisfies (E) by the arguments of Ikeda (see [Ik2, Thm. 2.7], and also the discussion following Theorem 3.3.11 of this paper). But one does *not yet* have (BV), which is problematic, which we propose to rectify in this paper. One knows that at the finite unramified places ([PSR1]) and the infinite places ([Ik3]) that

$$(3.1.7) \qquad L_1(s, \pi_v \times \pi'_v \times \eta_v) = L(s, \pi_v \times \pi'_v \times \eta_v).$$

See Section 3.3 for a more complete discussion.

The third candidate, which we denote by $L_2(s, \pi \times \pi' \times \eta)$, was constructed in the work of Shahidi (cf. [Sh1, p. 582]), developing the ideas of Langlands in the monograph ([La2]). (See also [PSR2, §0] for an exposition.) To be precise, let $G$ denote the group $\mathrm{Spin}(4,4)$. It has a parabolic $P = MU$ with Levi component $M$ which is a quotient (over $F$) of $\mathrm{GL}(1) \times \mathrm{SL}(2) \times \mathrm{SL}(2) \times \mathrm{SL}(2)$ by $\{1, -1\}$, where $-1$ signifies $H_1 H_3 H_4(-1)$, with the $H_j$ denoting the standard simple coroots. Take the representation $\beta := \omega \otimes \pi_1 \otimes \pi'_1 \otimes \tau_1$ on $M(\mathbb{A}_F)$, with $\pi_1$ (resp. $\pi'_1$, resp. $\tau_1$) being *any* irreducible component of the restriction of $\pi$ (resp. $\pi'$, resp. $\tau$) to $\mathrm{SL}(2, \mathbb{A}_F)$ and $\omega$ being the product of the central characters of $\pi$, $\pi'$ and $\tau$. Extend $\beta$ in the usual way to $P(\mathbb{A}_F)$ by letting $U(\mathbb{A}_F)$ act trivially. Then the Langlands-Eisenstein series on $G(\mathbb{A}_F)$ associated to the representation of $G(\mathbb{A}_F)$ induced by $\beta$ defines, and gives analytic information on, $L_2(s, \pi \times \pi' \times \eta)$. One knows that this $L$-function satisfies (MC) and (FE) (cf. [Sh1]), and also (E) by the very recent work of H. Kim and F. Shahidi([K-Sh]). If $v$ is archimedean or *unramified* or *tempered*, we have

$$(3.1.8) \qquad L_2(s, \pi_v \times \pi'_v \times \eta_v) = L(s, \pi_v \times \pi_v \times \eta_v).$$

There is also a definition of epsilon factors $\varepsilon_2(s, \pi_v \times \pi'_v \times \eta_v)$. The main use of this approach for us comes from knowing that $L_2(s, \pi \times \pi' \times \tau)$ is nonzero on the line $s = 1$.

We take $T$ to be a (finite) set containing all the places $v$ where $\pi_v$ or $\pi'_v$ is in the discrete series, and solve (P2) for the relevant types of $(\pi, \pi')$. In such a case, one first notes, by an extension of the local results of [Ik1], [Ik2] by global arguments, that the local gamma factors (called $\varepsilon'$-factors by some) of the first and third candidates agree at *all* the places, and that their $L$-factors also agree at all the places $v$ where at least one of the representations $\pi_v, \pi'_v$ is not supercuspidal. It is not clear that the $L$-factors agree at all the places $v$. But luckily, it turns out not to matter by a base changing trick. We refer to Sections 3.3 and 3.4 for details.

Next comes the establishment of (BV). The idea here is to make use of Arthur's truncation $\Lambda^T E(f_s)$ of the Eisenstein series $E(f_s)$ on $\mathrm{GSp}(6)/F$ associated to the character $\det^{s+2}$ of the (GL(2) occurring in the Levi of the) Siegel parabolic $P$. (Here $T$ is a sufficiently regular parameter.) The integral representation for $L_1(s, \pi_1 \times \pi_2 \times \pi_3)$ is given by the integral, over $\mathrm{GL}(2)^3$, of $E(f_s)$ against a function $\varphi$ in the space of $\pi_1 \otimes \pi_2 \otimes \pi_3$. By standard arguments using the functional equation and boundedness in $\Re(s) >> 0$, it suffices to show that this integral is of bounded order in vertical strips $\Xi$. Since $\varphi$ decreases rapidly at infinity, it suffices to check that for $s$ in $\Xi$, $E(f_s)(g)$ is bounded in absolute value by a function of $s$ of bounded order times $||g||^N$ for a uniform $N$ as $N$ goes to infinity.

Now, using the fact that $E(f_s)$ is an eigenfunction, though not square-integrable globally, of the Laplacian $\Delta$, we use standard Sobolev estimates and bound $|E(f_s)(g)|$ by its $L^2$-norm in a small neighborhood $V$ of controlled size. The $L^2$-norm over $V$ of $E(f_s) - \wedge^T E(f_s)$ can be estimated by use of the constant terms; these are essentially given by abelian $L$-functions, which are known by Hecke to be of bounded order. It then suffices to estimate the $L^2$ norm of $\wedge^T E(f_s)$ over the whole locally symmetric manifold, which can be understood by the works of Langlands and Arthur.

Finally, to be able to apply the converse theorem, we must also check (3.1.5). The details are in Sections 3.4 and 3.5.

3.2. *Weak to strong lifting, and the cuspidality criterion.* In this section we prove the following:

PROPOSITION 3.2.1.   *Let $\pi, \pi'$ be in $\mathcal{A}(2, F)$. Suppose we have constructed a weak lifting; i.e., an isobaric automorphic representation $\pi \boxtimes \pi'$ of $\mathrm{GL}(4, \mathbb{A}_F)$ satisfying the identity $(L_v)$ (of Theorem M) at almost all $v$. Then we have $(L_v)$ and $(\varepsilon_v)$ at all the finite places $v$, and $(L_\infty)$ as well. In addition, the cuspidality criterion (of Theorem M) holds.*

*Proof.* Let $S$ be the (finite) set of places outside which $(L_v)$ holds. We may assume (see (I) of §3.1) that $\pi$ and $\pi'$ are both cuspidal and nondihedral. Note also that the central character of $\pi \boxtimes \pi'$ is simply $(\omega \omega')^2$, where as before, $\omega, \omega'$ denote the respective central characters of $\pi, \pi'$. This is so because the two idele class characters agree almost everywhere. (In fact, by a theorem of Hecke, it suffices to know that they agree at a set of primes of density $> 1/2$.)

First some notation: If $f(s), g(s)$ are two meromorphic functions of $s$ such that their quotient is invertible, we will write $f(s) \equiv g(s)$. At any place $v$, given a character $\nu$ of $F_v^*$, we can write it as $\nu_0 |.|^z$, for a unitary character $\nu_0$ and a complex number $z$. The real part of $z$ is uniquely defined; we will call it the *exponent of $\nu$*, and denote it $e(\nu)$.

Choose a finite order character $\mu$ of $C_F$ such that $\mu_v$ is sufficiently ramified at every finite place $u$ in $S$ so as to make the Euler $u$-factors of $\pi \boxtimes \pi'$, $\pi \times \pi'$, and their contragredients, equal 1. This is possible by using the results of [JPSS]. Then, comparing the global functional equations of both $L$-functions, we get

(3.2.2)
$$\prod_{w|\infty} L(s, (\pi \boxtimes \pi')_w) L(1-s, \pi_w^\vee \times \pi_w'^\vee) \equiv \prod_{w|\infty} L(1-s, (\pi \boxtimes \pi')_w^\vee) L(s, \pi_w \times \pi_w').$$

For *any* place $v$, archimedean or otherwise, for any $n \geq 1$, and for any cuspidal automorphic representation $\Pi$ of $\mathrm{GL}(n, \mathbb{A}_F)$ (such as $\pi \boxtimes \pi'$), it is known that $L(s, \Pi_v)$ is holomorphic in $\Re(s) > \frac{1}{2} - t$, for some $t = t(\Pi, v) > 0$ (see [BaR, Prop. 2.1, part B]). Consequently, $L(s, (\pi \boxtimes \pi')_v)$ has no pole in common with $L(1 - s, (\pi \boxtimes \pi')_v^\vee)$.

Next we claim that, since $\pi, \pi'$ are (cuspidal) automorphic representations of $\mathrm{GL}(2, \mathbb{A}_F)$, $L(s, \pi_v \times \pi_v')$ is also holomorphic in $\Re(s) > \frac{1}{2} - t$, for some $t > 0$. If $\pi_v$ and $\pi_v'$ are both in the discrete series, then one has holomorphy even in $\Re(s) > 0$ (see [BaR, Lemma 2.3]). If $\pi_v$ is a principal series representation defined by (quasi)characters $\nu, \nu'$, then by [GJ], $e(\nu) < 1/4$ and $e(\nu') < 1/4$. (Strictly speaking, when $v$ is archimedean, one finds in [GJ] only the assertion that these exponents are $\leq 1/4$, but one can eliminate the possibility of exponent $1/4$ by using a simple version of the argument of [LRS1].) At the finite places, one can even replace $1/4$ by $1/5$ by using [Sh2]; see also [LRS2]. Over $\mathbb{Q}$, one can do still better and reduce the bound to $5/28$ ([BuDHI]). But we will not need these finer results. In any case, the claim follows easily.

Consequently, $L(s, \pi_v \times \pi_v')$ has no pole in common with $L(1-s, \pi_v^\vee \times \pi_v'^\vee)$ either. Applying this, along with the earlier remark on the relative primality of $L(s, (\pi \boxtimes \pi')_v)$ and $L(1 - s, (\pi \boxtimes \pi')_v^\vee)$, we get $(L_\infty)$ from (3.2.2).

Next pick any finite $u_0 \in S$ and choose a $\mu$ which is 1 at $u_0$, but is highly ramified at all other finite $u$ in $S$. Again, by comparing the functional equations, we get an identity such as (3.2.2), with the factors at infinity replaced by the corresponding ones at $u_0$. Using the relative primality results above, now with $v = u_0$, we get $(L_{u_0})$. Since $u_0$ was arbitrary, we get $(L_v)$ now at every place $v$.

For the identity of epsilon factors, we fix $u_0$ in $S$ and note that by [JPSS2], we can choose a global character $\mu$ which is 1 at $u_0$ and sufficiently ramified at any other place $u$ in $S$ such that the epsilon factor of $\pi_u \times \pi_u' \otimes \mu_u$ depends only on $\mu_u$ and the square of the product of the central characters $\omega_u, \omega_u'$, and the dependence is simple. Similarly, the epsilon factor of $(\pi \boxtimes \pi')_u$ has the same dependence on $\mu_u$ and the central character of $(\pi \boxtimes \pi')_u$, which we noticed above to be $(\omega_u \omega_u')^2$. The analogous statements hold for the contragredients, and this results in the identity (for all $u \in S - \{u_0\}$)
$$\frac{\varepsilon(s, (\pi \boxtimes \pi')_u)L(s, (\pi \boxtimes \pi')_u)}{L(1-s, (\pi \boxtimes \pi')_u^\vee)} = \frac{\varepsilon(s, \pi_u \times \pi'_u)L(s, \pi_u \times \pi'_u)}{L(1-s, \pi_u^\vee \times \pi_u'^\vee)}.$$
Comparing the global functional equations again, we get $(\varepsilon_{u_0})$ as we already know that the $L$-factors agree. This finishes the proof of the first part of the proposition.

It is left to prove the *cuspidality criterion*. First suppose $\pi \boxtimes \pi'$ is cuspidal. Recall that by (I) of Section 3.1, $\pi \boxtimes \pi'$ is not cuspidal if $\pi$ or $\pi'$ is not. So $\pi$ and $\pi'$ must both be cuspidal. Suppose $\pi' \simeq \pi \otimes \chi$, for an idele class character $\chi$ of $F$. Then, by [JS2] and $(L_v)$ for almost all $v$, $L^S(s, \pi \boxtimes \pi' \otimes (\chi\omega)^{-1})$ must have a pole at $s = 1$, with $S$ denoting the set of ramified and archimedean places, and $\omega$ the central character of $\pi$. (We are using the fact that $\pi \otimes \omega^{-1}$ is the contragredient $\pi^\vee$ of $\pi$.) Then $\pi \boxtimes \pi'$ cannot be cuspidal, leading to a contradiction.

Conversely, suppose that $\pi, \pi'$ are cuspidal and (C) holds, but that $\pi \boxtimes \pi'$ is not cuspidal. We will show that it leads to a contradiction.

Since $\pi \boxtimes \pi'$ is isobaric, we can decompose it uniquely as
$$\pi \boxtimes \pi' = \boxplus_{j=1}^r \eta_j,$$
where each $\eta_j$ is a cuspidal automorphic representation of $\mathrm{GL}(n_j, \mathbb{A}_F)$, with $\sum_{j=1}^r n_j$ being 4. Suppose some $n_j$, say $n_1$, is 1. Then $\eta_1$ is an idele class character, and we have, for any large enough finite set $S$ of places containing the archimedean ones,
$$L^S(s, \pi \times \pi' \otimes \eta_1^{-1}) = \zeta_F^S(s) \prod_{j \neq 1} L^S(s, \eta_j \otimes \eta_1^{-1}).$$
One knows ([JS2]) that, since each $L^S(s, \eta_j \otimes \eta_1^{-1})$ is cuspidal, it cannot vanish at $s = 1$. Thus the pole (at $s = 1$) of $\zeta_F^S(s)$ induces one of the function on the left, which is not allowed by (C). Thus $n_j$ must be $> 1$ for each $j$. Then we must have
$$n_1 = n_2 = 2.$$

Now the key idea is to compute the exterior square $L$-function of $\pi \boxtimes \pi'$ in two different ways. On the one hand, since it is of the form $\eta_1 \boxplus \eta_2$, with each $\eta_j$ a cuspidal of $\mathrm{GL}(2)/F$, we get
$$L^S(s, \pi \boxtimes \pi', \Lambda^2) = L^S(s, \eta_1 \times \eta_2)L^S(s, \omega_1)L^S(s, \omega_2),$$
where $\omega_j$ is the central character of $\eta_j$. This can be seen easily at the unramified places $v$. Indeed, if $\tau_{j,v}$ is the 2-dimensional representation of $W'_{F_v}$ associated to $\eta_{j,v}$, then the above identity is induced by the following:
$$\Lambda^2(\tau_{1,v} \otimes \tau_{2,v}) \simeq (\tau_1 \otimes \tau_2) \oplus \det(\tau_{1,v}) \oplus \det(\tau_{2,v}),$$
which is easy to verify.

Consequently, $L^S(s, \pi \boxtimes \pi', \Lambda^2)$ is divisible by (two) abelian $L$-functions, namely $L^S(s, \omega_1)$ and $L^S(s, \omega_2)$.

On the other hand, denoting by $\sigma_v$ (resp. $\sigma'_v$) the 2-dimensional representation of $W'_{F_v}$ associated to $\pi_v$ (resp. $\pi'_v$), we also have the identity

$$\Lambda^2(\sigma_v \otimes \sigma'_v) \simeq \Lambda^2(\sigma_v) \otimes \mathrm{sym}^2(\sigma'_v) \oplus \mathrm{sym}^2(\sigma_v) \otimes \Lambda^2(\sigma'_v).$$

This implies the following equality of $L$-functions:

$$L^S(s, \pi \boxtimes \pi', \Lambda^2) = L^S(s, \mathrm{sym}^2(\pi) \otimes \omega) L^S(s, \mathrm{sym}^2(\pi') \otimes \omega'),$$

where $\mathrm{sym}^2(\pi)$ (resp. $\mathrm{sym}^2(\pi')$) is the automorphic representation of $\mathrm{GL}(3, \mathbb{A}_F)$ associated to $\pi$ (resp. $\pi'$) by Gelbart-Jacquet ([GJ]), and $\omega$ (resp. $\omega'$) is the central character of $\pi$ (resp. $\pi'$).

Suppose $\pi$ is not automorphically induced by an idele class character of a quadratic extension. Then the main theorem of [GJ] says that $\mathrm{sym}^2(\pi)$ is cuspidal; so $L^S(s, \mathrm{sym}^2(\pi) \otimes \omega')$ cannot be divisible by an abelian $L$-function. Consequently, if neither $\pi$ nor $\pi'$ is automorphically induced, we get a contradiction, and so $\pi \boxtimes \pi'$ must be cuspidal.

It remains to consider the case when $\pi'$ is of the form $I_K^F(\mu)$, for $\mu$ an idele class character of a quadratic extension $K$. Then by (II) of Section 3.1, we know that $\pi \boxtimes \pi'$ must be isomorphic to $I_K^F(\beta)$, where $\beta = \pi_K \otimes \mu$. By hypothesis, $\pi_K$ is cuspidal and not isomorphic to $\pi_K \otimes (\mu \circ \theta)\mu^{-1}$, where $\theta$ is the nontrivial automorphism of $K/F$. Then $\beta$ is a cuspidal, not isomorphic to $\beta \circ \theta$. So by the properties of automorphic induction (Prop. 2.3.3), $\pi \boxtimes \pi'$ is cuspidal. □

3.3. *Triple product L-functions: local factors and holomorphy.* Let $\pi, \pi', \pi''$ be unitary, irreducible, cuspidal automorphic representations of $\mathrm{GL}(2, \mathbb{A}_F)$ of central characters $\omega, \omega', \omega''$ respectively. Put $\omega = \omega\omega'\omega''$. At each place $v$, let $\sigma_v, \sigma'_v, \sigma''_v$ be the 2-dimensional representations of $W'_{F_v}$ associated to $\pi_v, \pi'_v, \pi''_v$ respectively, by the local Langlands correspondence for $\mathrm{GL}(2)$ ([Ku]).

Recall that $L(s, \pi \times \pi' \times \pi'')$ is the (correct) triple product $L$-function defined by the Euler product $\prod_v L(s, \sigma_v \otimes \sigma'_v \otimes \sigma''_v)$, while $L_1(s, \pi \times \pi' \times \pi'')$ is the one given as the greatest common denominator of the integral representation of Garrett and Piatetski-Shapiro and Rallis, and $L_2(s, \pi \times \pi' \times \pi'')$ is the one defined by the Langlands-Shahidi method.

The following basic fact, due to Ikeda in the archimedean case ([Ik3, Thm. 1.1.0]), and Rallis and Piatetski-Shapiro in the non-archimedean case ([PS-R2]), will be very important to us.

PROPOSITION 3.3.1. *Let $v$ be any place of $F$ where $\pi_v, \pi'_v, \pi''_v$ are all unramified. Then*

$$L(s, \pi_v \times \pi'_v \times \pi''_v) = L_1(s, \pi_v \times \pi'_v \times \pi''_v).$$

By a slight refinement of known results we will also establish the following:

PROPOSITION 3.3.2. *Let $v$ be any place of $F$, where at least one of the local representations $\pi_v, \pi'_v, \pi''_v$ is in the principal series (possibly complementary). Then*

(a) $$L(s, \pi_v \times \pi'_v \times \pi''_v) = L_1(s, \pi_v \times \pi'_v \times \pi''_v),$$

*and*

(b) $$\varepsilon(s, \pi_v \times \pi'_v \times \pi''_v) = \varepsilon_1(s, \pi_v \times \pi'_v \times \pi''_v).$$

*Remark.* When $v$ is finite, if two of the local representations are non-supercuspidal, then such an identity already follows from the works of Ikeda ([Ik2,3]). In the special, but important, case when all of $\pi_v, \pi'_v, \pi''_v$ are Steinberg, this was established earlier by Gross and Kudla ([GK]). When $v$ is archimedean, this is treated in [Ik3]. In this case, the $L_1$-factor is defined only up to an invertible holomorphic function, and so the content of this proposition is that we can normalize it appropriately so as to have the stated identities. Note also that one knows by the work of Shahidi ([Sh5]) the identity $L(s, \pi_v \times \pi'_v \times \pi''_v) = L_2(s, \pi_v \times \pi'_v \times \pi''_v)$ for any archimedean $v$.

We have stated here only what we need, and a stronger assertion will be proved later in Section 4.4 *after* establishing Theorem M.

*Proof.* By hypothesis, there exist quasi-characters $\mu, \nu$ of $F_v^*$ such that one of the local representations, say $\pi_v$, is the principal series defined by $(\mu, \nu)$. The starting point for us is the following:

LEMMA 3.3.3. *The set of poles of $L_1(s, \pi_v \times \pi'_v \times \pi''_v)$ is contained in the set of poles of $L(s, \pi_v \times \pi'_v \times \pi''_v)$.*

*Proof of Lemma* 3.3.3. We first recall the following:

PROPOSITION 3.3.4. *For any quasi-character $\lambda$ of $F_v^*$,*

$$L(s, \lambda \otimes \pi'_v \otimes \pi''_v) = L(s, \lambda \otimes \sigma'_v \otimes \sigma''_v)$$

*and*

$$\varepsilon(s, \lambda \otimes \pi'_v \otimes \pi''_v) = \varepsilon(s, \lambda \otimes \sigma'_v \otimes \sigma''_v).$$

For a proof, see the author's paper with D. Prasad ([PR, Prop. 4.2]), where the identity is established via a global argument.

On the other hand, it is immediate that

(3.3.5) $$L(s, \sigma_v \otimes \sigma'_v \otimes \sigma''_v) = L(s, \mu \otimes \sigma'_v \otimes \sigma''_v) L(s, \nu \otimes \sigma'_v \otimes \sigma''_v),$$

and

$$\varepsilon(s, \sigma_v \otimes \sigma'_v \otimes \sigma''_v) = \varepsilon(s, \mu \otimes \sigma'_v \otimes \sigma''_v) \varepsilon(s, \nu \otimes \sigma'_v \otimes \sigma''_v).$$

One also has the following result of Ikeda:

PROPOSITION 3.3.6 ([Ik3, Thm. 1.8]). *The quotient*

$$\frac{L_1(s, \pi_v \times \pi'_v \times \pi''_v)}{L(s, \mu \otimes \pi'_v \otimes \pi''_v) L(s, \nu \otimes \pi'_v \otimes \pi''_v)}$$

*is entire.*

Set

$$\gamma(s, \pi_v \times \pi'_v \times \pi''_v) = \varepsilon(s, \pi_v \times \pi'_v \times \pi''_v) \frac{L(1-s, \pi_v^\vee \times (\pi'_v)^\vee \times (\pi''_v)^\vee)}{L(s, \pi_v \times \pi'_v \times \pi''_v)}.$$

Define $\gamma_j(s, \pi_v \times \pi'_v \times \pi''_v)$ analogously, for $j = 1, 2$, using the $L_j$ and $\varepsilon_j$-factors instead. (The gamma factor is called the $\varepsilon'$-factor by some, e.g., in [Ik1,2].)

In view of Propositions 3.3.4, 3.3.6 and identity (3.3.5), Lemma 3.3.3 will follow once we establish

PROPOSITION 3.3.7.

$$\gamma(s, \pi_v \times \pi'_v \times \pi''_v) = \gamma_1(s, \pi_v \times \pi'_v \times \pi''_v).$$

Using [Sh4], one can show such an identity with $\gamma_1$ replaced by $\gamma_2$.

*Proof of Proposition* 3.3.7. By Theorem 3 of [Ik2], we know that

$$(3.3.8) \qquad \gamma_1(s, \pi_v \times \pi'_v \times \pi''_v) = \gamma_1(s, \mu \otimes \pi'_v \times \pi''_v) \gamma_1(s, \nu \otimes \pi'_v \times \pi''_v).$$

By applying the identities (3.3.4), we then see that

$$(3.3.9) \qquad \gamma_1(s, \pi_v \times \pi'_v \times \pi''_v) = \gamma_1(s, \mu \otimes \sigma'_v \otimes \sigma''_v) \gamma_1(s, \nu \otimes \sigma'_v \otimes \sigma''_v).$$

In view of (3.3.5), the right-hand side of (3.3.9) equals $\gamma(s, \pi_v \times \pi'_v \times \pi''_v)$ as claimed. □

As noted above, this also proves Lemma 3.3.3.

*Proof of Proposition* 3.3.2. There is an index $\lambda(\tau)$ of a unitary irreducible $\tau$ of $\mathrm{GL}(2, F_v)$ which measures the failure of $\tau$ to be tempered; it is $t$, resp. 0, if $\tau$ is a complementary series attached to $(\mu|.|^t, \mu|.|^{-t})$ with $t > 0$ and $\mu$ unitary, resp. $\tau$ is tempered, and set (as in [Ik1]) $\lambda(\pi_v, \pi'_v, \pi''_v) = \lambda(\pi_v) + \lambda(\pi'_v) + \lambda(\pi''_v)$. Then it is easy to see that

$$(3.3.10) \qquad L(s, \sigma_v \otimes \sigma'_v \otimes \sigma''_v) \text{ is holomorphic in } \Re(s) > \lambda(\pi_v, \pi'_v, \pi''_v).$$

(This is the Galois analog of Lemma 2.1 of [Ik1].) We know by Shahidi ([Sh2]) that $\lambda(\tau)$ is always $< \frac{1}{5}$ for $\tau \in \{\pi_v, \pi'_v, \pi''_v\}$. Consequently, $\lambda(\pi_v, \pi'_v, \pi''_v)$ is always $< 1/2$ unless all the three representations are nontempered. In the former case, i.e., if $\lambda(\pi_v, \pi'_v, \pi''_v) < 1/2$, then $L(s, \sigma_v \otimes \sigma'_v \otimes \sigma''_v)$ and $L(1 - s, \sigma_v^\vee \otimes (\sigma_v')^\vee \otimes (\sigma_v'')^\vee)$ have no common poles. The same holds for the $L_1$-functions by Lemma 3.3.3, and the assertion follows in this case.

It remains to treat the case when $\lambda(\pi_v, \pi_v', \pi_v'') \geq 1/2$ (which should not happen!). Suppose $\pi_v$ is in the principal series associated to the characters $\mu, \nu$. The unitarity implies that the set $\{\overline{\mu}, \overline{\nu}\}$ equals $\{\mu^{-1}, \nu^{-1}\}$, and so one of the following happens: (i) $\overline{\mu} = \mu^{-1}$; (ii) $\overline{\mu} = \nu^{-1}$. In case (i), $\mu$ and $\nu$ are unitary, so we are in case (ii). Now $\pi_v$ is the complementary series representation associated to the pair of characters $(\mu, \overline{\mu}^{-1})$. We may write $\mu = \chi |.|^t$, for some unitary (possibly ramified) character $\chi$ and real number $t$. Then $\nu = (\overline{\chi}^{-1})|.|^{-t} = \chi |.|^{-t}$, and so $\pi_v$ is $\pi_{1,v} \otimes \chi$ with $\pi_{1,v}$ unramified. Similarly, $\pi_v'$ (resp. $\pi_v''$) is $\pi_{1,v}' \otimes \chi'$ (resp. $\pi_{1,v}'' \otimes \chi''$), with $\pi_{1,v}', \pi_{1,v}''$ unramified and $\chi', \chi''$ unitary. Then the integral representation for triple product $L$-functions gives the identity (with $\xi = \chi\chi'\chi''$):

$$L_1(s, \pi_v \times \pi_v' \times \pi_v'') = L_1(s, \pi_{1,v} \times \pi_{1,v}' \times (\pi_{1,v}'' \otimes \xi)).$$

Consequently, if all three local representations are nontempered, we may assume, by replacing $(\pi_v, \pi_v', \pi_v'')$ by $(\pi_{1,v}, \pi_{1,v}', \pi_{1,v}'' \otimes \xi)$, that at most one of them is ramified. If $\xi$ is ramified, then both the $L$- and the $L_1$-functions of $(\pi_v, \pi_v', \pi_v'')$ at $s$ do not share any pole with the corresponding contragredient functions at $1 - s$, and the assertion of Proposition 3.3.2 follows. When $\xi$ is unramified, we can appeal to Proposition 3.3.1. □

THEOREM 3.3.11. *Let $\pi, \pi', \pi''$ be cuspidal automorphic representations of* $\mathrm{GL}(2, \mathbb{A}_F)$ *satisfying the following*:

- *At least one of $\{\pi, \pi', \pi''\}$ is nondihedral; and*

- *At every finite place $v$, at least one of $\{\pi_v, \pi_v', \pi_v''\}$ is in the principal series.*

Then $L(s, \pi \times \pi' \times \pi'')$ is entire.

In view of Proposition 3.3.2, it suffices to show that $L_1(s, \pi \times \pi' \times \pi'')$ is entire, and this follows from the work of Ikeda. But a remark is in order. In [Ik1], one finds an assertion (see Theorem 2.7 there) that $L_1(s, \pi \times \pi' \times \pi'')$ is entire as long as one of the three representations is not dihedral, *without any local conditions*. This is as it should be, but a small correction needs to be made in his proof. To be precise, in the case when $\omega^2 = 1, \omega \neq 1$ with associated quadratic extension $K/F$ (see [Ik1, pp. 231–234]), he first shows that $L_1(s, \pi_K \times \pi_K' \times \pi_K'')$ is entire and then proceeds to assert that *any pole of $L_1(s, \pi \times \pi' \times \pi'')$ on $\Re(s) = 1$ is a pole of $L_1(s, \pi_K \times \pi_K' \times \pi_K'')$*. We do not see how to verify this claim; the problem is with the bad factors. The way he tries to derive this is to assert the following two identities:

(Ik1) $\qquad L_1(s, \pi \times \pi' \times \pi'') = L_2(s, \pi \times \pi' \times \pi''),$

and

(Ik2) $\quad L_2(s, \pi_K \times \pi'_K \times \pi''_K) = L_2(s, \pi \times \pi' \times \pi'') L_2(s, \pi \times \pi' \times \pi'' \otimes \omega),$

and then appeal to Shahidi's result ([Sh6]):

(3.3.12) $\qquad L_2(s, \pi \times \pi' \times \pi'' \otimes \omega) \neq 0 \quad \text{if} \quad \Re(s) = 1.$

We are unable to verify either of the identities (Ik1), (Ik2), though they both hold locally *almost everywhere*.

This can be easily fixed as follows. First note that at any place $v$, (3.3.8) holds with $L$ replaced by $L_2$. Since $3/5$ is the maximum possible value of $\lambda(\pi, \pi', \pi'')$ (by [Sh1]), the local factors (of either of three $L$-functions) is invertible on $\{\Re(s) = 1\}$. So (3.3.12) holds with $L_2$ replaced by the incomplete $L$-function $L_2^S$, taken over the places outside a finite set $S$. Choose $S$ to be such that [Ik1], [Ik2] both hold with $L_1, L_2$ replaced by $L_1^S, L_2^S$ respectively. (Over $K$, we remove the factors at places above $S$.) It follows that $L_1(s, \pi \times \pi' \times \pi'')$ has a pole if and only if $L^S(s, \pi \times \pi' \times \pi'')$ does, which happens if and only if $L^S(s, \pi_K \times \pi'_K \times \pi''_K)$ and (hence) $L(s, \pi_K \times \pi'_K \times \pi''_K)$ do. The rest of Ikeda's argument goes through verbatim.

3.4. *Boundedness in vertical strips.* A very important prerequisite to the application of converse theorems is the knowledge that the relevant $L$-functions are, in addition to being entire, bounded in vertical strips. This is what we prove here for a class of triple product $L$-functions. After the proof of the main theorem in Section 3.8, the same will hold for all such $L$-functions.

THEOREM 3.4.1. *Let $F$ be a totally complex number field, and let $\pi_1, \pi_2, \pi_3$ be unitary, cuspidal automorphic representations of $\mathrm{GL}_2(\mathbb{A}_F)$. Assume that at any finite place $v$, at least one $\pi_{j,v}$ is in the principal series. Then the entire function $L(s, \pi_1 \times \pi_2 \times \pi_3)$ is bounded in vertical strips of finite width.*

Once we prove Theorem M completely in Section 3.7, this theorem will become valid with no hypotheses on $F$ and on the $\pi_j$, because the $L$-function will become a Rankin-Selberg $L$-function on $\mathrm{GL}(4) \times \mathrm{GL}(2)$.

*Proof.* In view of Propositions 3.3.1 and 3.3.2, we only need to prove the assertion of the theorem for $L_1(s, \pi_1 \times \pi_2 \times \pi_3)$, which has an integral representation involving a Siegel Eisenstein series $\mathcal{E}(f)$ on $\mathrm{GSp}_6/F$. Our key idea is to appeal to the fact that $\mathcal{E}(f)$ is an eigenfunction of the Laplacian and then to make use of Arthur's truncation of $\mathcal{E}(s)$.

We now need to set up notation and preliminaries. Let

$$G = \mathrm{GSp}(6)/F = \{g \in \mathrm{GL}_6/F \mid {}^t g J g = \lambda J, \text{ for some } \lambda \in \mathbb{G}_m\},$$

where $J = \begin{pmatrix} 0 & -I_3 \\ I_3 & 0 \end{pmatrix}$, with $I_n$ denoting, for any $n \geq 1$, the identity $n \times n$-matrix. The root system relative to the diagonal torus $A$ in $G$ is of rank 3 and is described by

$$\Phi = \{\pm e_i \pm e_j \mid 1 \leq i, j \leq 3\},$$

where $e_1, e_2, e_3$ are standard basis vectors of Lie $T_{\mathbb{C}}^*$. Pick the following set of simple roots:

$$\Delta = \{\alpha = e_1 - e_2, \beta = e_2 - e_3, \gamma = 2e_3\}.$$

Denote by $P_0 = AU_0$ the corresponding minimal parabolic subgroup of $G$, whose elements are matrices of the form $\begin{pmatrix} A & B \\ 0 & D \end{pmatrix}$, $A, B, D \in M_3$, with $A$ (resp. $D$) upper (resp. lower) triangular. For each subset $\theta$ of $\Delta$, let $A_\theta$ denote $(\cap_{\delta \in \theta} \ker(\delta))^0 \subset A$, $M_\theta$ be the centralizer of $A_\theta$ in $G$, and $P_\theta = M_\theta U_\theta$ be the corresponding *standard parabolic* subgroup of $G$, split over $F$, containing $P_0$. The *Weyl group* $W$ is a semidirect product $(\mathbb{Z}/2)^3 \rtimes S_3$, whose elements are made up of permutations of $\{e_1, e_2, e_3\}$ (possibly) composed with flips (sign changes) $e_i \mapsto -e_i$, for $i \in \{1, 2, 3\}$.

Let $P = P_{\{\alpha,\beta\}} = MU$ be the (maximal) Siegel parabolic subgroup with modular function $\delta p$. We have the matrix representation

$$P = \left\{ p = \begin{pmatrix} \lambda A & 0 \\ \hline 0 & {}^t A^{-1} \end{pmatrix} \begin{pmatrix} I & X \\ \hline 0 & I \end{pmatrix} \Big| \lambda \in \mathbb{G}_m, A \in \mathrm{GL}_3 X = {}^t X \right\},$$

with $\delta_p(p) = |\det A|^6$. Let $K = \prod_v K_v$ be the standard maximal compact mod center subgroup of $G(\mathbb{A}_F)$.

Let $\omega_i$ be the central character of $\pi_i$, for each $i \leq 3$. Put $\omega = \omega_1 \omega_2 \omega_3$ and consider the space $J(\omega, s)$ of right $K$-finite functions $f_s$ on $G(\mathbb{A}_F)$ satisfying for all $p \in P(\mathbb{A}_F)$ and $g \in G(\mathbb{A}_f))$:

$$(3.4.2) \qquad f_s(pg) = \omega(\lambda) |\lambda|^{3s+\frac{3}{2}} \omega(\det A) |\det A|^{2s+1} f_s(g).$$

We will henceforth take $f_s$ to be a *decomposable, entire and $K$-finite section*, which is *good* in the sense of [Ik1], [PS-R]. The corresponding local space at $v$ will be denoted $J_v(\omega, s)$. The relevant Eisenstein series is then given by the formula (for all $g \in G(\mathbb{A}_F)$)

$$(3.4.3) \qquad E(f_s)(g) = \sum_{\gamma \in P(F) \backslash G(F)} f_s(\gamma g).$$

(Note that $J(\omega, s)$ corresponds to $I(\omega, 2s - 1)$ in Proposition 1.6 of [Ik1].)

Let $\varphi = \varphi_1 \otimes \varphi_2 \otimes \varphi_3$ be a cusp form in the space of $\pi_1 \otimes \pi_2 \otimes \pi_3$ on $H(\mathbb{A}_F)$, where

$$H = \{(g_1, g_2, g_3) \in \mathrm{GL}_2^3 \mid \det g_1 = \det g_2 = \det g_3\}.$$

Denote the corresponding Whittaker function by $W = W_1 \otimes W_2 \otimes W_3$. Then one knows (see [Ik1], [PS-R2]) that in $\mathrm{Re}(s) \gg 0$ one has the identity (with $f_s = \otimes_v f_{v,s}$):

$$(3.4.4) \quad \langle E(f_s), \varphi \rangle_H := \int_{C(\mathbb{A}_F)H(F)\backslash H(\mathbb{A}_F)} E(h, f_s)\varphi(h)dh = \prod_v \Psi(f_{v,s}; W_v),$$

where each $\Psi(f_{v,s}; W_v)$ is a local integral having $L_1(s, \pi_{1,v} \times \pi_{2,v} \times \pi_{3,v})$ as its greatest common denominator. By [PS-R2] one knows that at every unramified $v$, $\Psi(f_{v,s}; W_v)$ coincides with $L_1(s, \pi_{1,v} \times \pi_{2,v} \times \pi_{3,v})$ for a suitable choice of $f_{v,s}$. Hence we get

$$L_1(s, \pi_1 \times \pi_2 \times \pi_3) \prod_{v \in S} \frac{\Psi(f_{v,s}; W_v)}{L_1(s, \pi_{1,v} \times \pi_{2,v} \times \pi_{3,v})} = \langle E(f_s), \varphi \rangle_H,$$

where $S$ is a finite set of places containing the ramified and archimedean places.

Recall that $L_1(s, \pi_1 \times \pi_2 \times \pi_3)$ has an Euler product in $\Re(s) > 1$ and a functional equation. It is clear that it is bounded for $\Re(s) \gg 0$. So by the Phragmen-Lindelöf theorem, Theorem 3.4.1 will be proved if we establish that this $L$-function is of bounded order in vertical strips of finite width.

LEMMA 3.4.5. *For each $v$, the function $\frac{\Psi(f_{v,s}; W_v)}{L_1(s, \pi_{1,v} \times \pi_{2,v} \times \pi_{3,v})}$ is entire and of bounded order, for a suitable choice of $f_{v,s}$.*

*Proof of Lemma* 3.4.5. Put $g(s, \pi_{1,v} \times \pi_{2,v} \times \pi_{3,v}) = \frac{\Psi(f_{v,s}; W_v)}{L_1(s, \pi_{1,v} \times \pi_{2,v} \times \pi_{3,v})}$. The entireness of this function follows from the fact that $L_1(s, \pi_{1,v} \times \pi_{2,v})$ is the greatest common denominator of the local integrals $\Psi(f_{v,s}; W_v)$ as the $(f_{v,s}, W_v)$ vary. The fact that $g(s, \pi_{1,v} \times \pi_{2,v} \times \pi_{3,v})$ is of bounded order is obvious at any $v$-adic place as the local integral is a rational function of $Nv^{-s}$. So let $v$ be archimedean. Since we have taken $F$ to be totally complex, $F_v$ is $\mathbb{C}$. There is a local functional equation

$$g(s, \pi_{1,v} \times \pi_{2,v} \times \pi_{3,v}) = \varepsilon_1(s, \pi_{1,v} \times \pi_{2,v} \times \pi_{3,v})g(1-s, \overline{\pi}_{1,v} \times \overline{\pi}_{2,v} \times \overline{\pi}_{3,v}^{\vee}),$$

where the intervening $\varepsilon$-factor is, thanks to Proposition 3.3.2 (b), an exponential function. Now, it suffices to prove that, for a suitable choice of $f_{v,s}$, the local integral $\Psi(f_{v,s}; W_v)$ is of bounded order in $\Re(s) \geq 1/2$.

Let $\lambda$ denote $\lambda(\pi_{1,v}) + \lambda(\pi_{2,v}) + \lambda(\pi_{3,v})$, where, for each $j$, $\lambda(\pi_{j,v})$ denotes, as in the proof of Proposition 3.3.2 above, the index of $\pi_{j,v}$ measuring its failure to be tempered. It is known ([Ik3]) that the integral $\Psi(f_{v,s}; W_v)$ is absolutely convergent in $\Re(s) > \lambda$, defining a holomorphic function there. Noting that $|f_{v,s}|$ transforms under the left action of $P(F_v)$ by a character depending only on $\Re(s)$, we see that $\Psi(f_{v,s}; W_v)$ is bounded in any vertical strip $\{a \leq \Re(s) \leq b\}$ if $a > \lambda$.

Now suppose that at least one $\pi_{j,v}$ is tempered. By Gelbart-Jacquet ([GJ]), one has $\lambda(\pi_{i,v}) < 1/4$ for any $i$. (For $F = \mathbb{Q}$, one even has the archimedean estimate: $\lambda(\pi_{i,v}) < 5/28$ by [LRS1].) So $\lambda < 1/2$ in this case, and by the remarks above, $\Psi(f_{v,s}; W_v)$ is bounded in any vertical strip $\{1/2 \leq \Re(s) \leq b\}$ of finite width.

It remains to consider when *every* $\pi_j$ is nontempered. Arguing as in the proof of Proposition 3.3.2, we may assume that $\pi_{1,v}$ and $\pi_{2,v}$ are unramified and that $\pi_{3,v}$ is the twist of an unramified representation by a unitary character $\xi$ which is ramified if nontrivial. If $\xi = 1$, i.e., when each $\pi_{j,v}$ is unramified, it is proved in [Ik3] that there exists a function $f_{v,s}$ such that $\Psi(f_{v,s}; W_v)$ equals $L_1(s, \pi_{1,v} \times \pi_{2,v} \times \pi_{3,v})$. So we are done in this case by the standard properties of the gamma function. Now suppose $\xi$ is nontrivial. We can then use the $\xi$-twisted form of the unramified $f_{v,s}$ employed by Ikeda, in which case $\Psi(f_{v,s}; W_v)$ will be essentially given by the integral expression preceding Lemma 2.3 in [Ik3], the change being that the integrand will get multiplied by $\xi(\sqrt{xyz}R^{-1}z)\tilde{\xi}(k)$, where $\tilde{\xi}$ is the character of the maximal compact $K_v$ whose restriction to $K_v \cap P(F_v)$ is $\xi$. (Here the notation is as in [Ik3].) Now we appeal to the fact that the Whittaker function $W_{1,v}$ (resp. $W_{2,v}$, resp. $W_{3,v}$) satisfies a second order differential equation on $\mathrm{diag}(x,1)$ (resp. $\mathrm{diag}(y,1)$, resp. $\mathrm{diag}(z,1)$). Using integration by parts, we can then express an elementary function times $\Psi(f_{v,s}; W_v)$ in terms of a finite sum of integrals, each of which converges absolutely in $\Re(s) > \lambda - 1/2$ and is bounded in $\{a \leq \Re(s) \leq b\}$ if $a > \lambda - 1/2$. Since $\lambda$ is $< 3/4$ by [GJ], we are done. $\square$

Since an entire function which is the ratio of two entire functions of bounded order is itself of bounded order, it suffices then to establish that the entire function $\langle E(f_s), \varphi \rangle_H$ is of bounded order. (Such a reduction was employed in [RuS] for an analogous situation.)

Since $\varphi$ vanishes rapidly at infinity, it will suffice to establish the following:

PROPOSITION 3.4.6. *Fix a vertical strip $\Xi$ in $\mathbb{C}$ of finite width, and a Siegel set $\mathcal{S}$ in $G(\mathbb{A}_F)$. Then there exist $N \in \mathbb{Z}$ and a function $\Phi(s)$ of of bounded order, depending only on the width of $\Xi$ and $\mathcal{S}$, such that*

$$|E(f_s)(g)| \leq |\Phi(s)| \, ||g||^N,$$

*for all $g$ in $\mathcal{S}$.*

By a function of bounded order we mean a quotient $\xi_1(s)/\xi_2(s)$ such that, for some $t_0 \geq 0$, $C > 0$, and $M \in \mathbb{Z}_+$, we have for all $s = x + iy \in \Xi$ with $y \geq t_0$, the bound $|\xi_j(s)| \leq Ce^{|y|^M}$, for $j = 1, 2$.

*Proof of Proposition* 3.4.6. Let us choose a sufficiently small compact, open subgroup $K$ of $G(\mathbb{A}_{F,f})$ such that $E(f_s)$ is right invariant under $K$, and such that $Y := G(F)Z_\infty^+ \backslash G(\mathbb{A}_F)/K$ is a smooth manifold. Fix an $\epsilon > 0$, and also a chart on $Y$ defining the differentiable structure. We will work with the Siegel set modulo $K$, though we will suppress this in our notation.

For every $g$ in $\mathcal{S}$, choose a neighborhood $V = V_g$ of $g$ with nonempty interior, which is small enough to lie inside the unit ball of radius $\epsilon$ in local coordinates.

Recall that $E(f_s)$ is an eigenfunction of the Laplacian $\Delta$ on all of $Y$, as it is defined by the Casimir operator, and that it is moreover *locally* in $L^p$, for any $p \geq 1$.

LEMMA 3.4.7. *Let $g, V$ be as above. Then there is a polynomial $P(s)$ independent of $g$ such that*

$$||E(f_s)(g)|| \leq |P(s)|.||E(f,s)||_{2,V}.$$

*Proof of Lemma* 3.4.7. Let $V' = V'_g$ be any neighborhood of $g$ with nonempty interior such that its closure is contained in the interior of $V$. It suffices to prove the asserted bound for $||E(f,s)||_{\infty,V'}$. (We write $||.||_{p,X}$ for the $L^p$-norm over $X$.) As remarked above, $E(f_s)$ is an eigenfunction, say, with eigenvalue $\lambda_s$, of the second order elliptic differential operator $\Delta$. It can be seen that $E(f_s)$ is in the localized Sobolev space $\mathcal{H}_n^{\mathrm{loc}}(V)$ for any $n$. And since $\bar{V}' \subset V$, $E(f_s)$ equals a function in $\mathcal{H}_n(V)$ on $V'$.

LEMMA 3.4.8. *There exists a polynomial $Q(X)$ in $\mathbb{C}[X]$, depending only on $\epsilon$ and $\lambda_s$ such that*

$$||E(f,s)||_{\infty,V'} \leq |Q(\lambda_s)|.||E(f,s)||_{2,V}.$$

Since $\lambda_s$ is a polynomial, in fact quadratic, in $s$, Lemma 3.4.7 follows from this. When $\lambda_s$ is real, which happens in our case if and only if $s$ lies on the unitary cross, the results of Section 4 of [So], in particular the identity (4.2) and Lemma 4.1, easily imply Lemma 3.4.8. The methods of [So] extend to the case of complex eigenvalues as well, but T. Wolff has pointed out to us how this lemma can be deduced more directly from the local regularity properties of elliptic differential operators; see Lemma 3.4.9 below. (The difficult part of [So] involves precisely bounding the degree of $Q(X)$, but it is not important for our purposes.) For any $r$, denote by $||.||_{(2,r);X}$ the norm associated to the $r^{\mathrm{th}}$ $L^2$-derivative on $X$ so that $||u||_{(2,r);X}$ equals $\sum_{|\nu| \leq r} ||\partial^\nu u||_{2,X}$, with the $\partial^\nu$ denoting distribution derivatives.

LEMMA 3.4.9. *Let $\Omega$ be a subset of $\mathbb{R}^N$ contained in the unit ball of radius $\epsilon$, and let $\Omega'$ be a subset of $\Omega$ with nonempty interior such that $\bar{\Omega}' \subset \Omega$. Then,*

(1) *For any integer $r > N/2$, there exists a constant $C_1 > 0$ depending only on $\epsilon$ and $\Delta$ such that, for all $u$ in $\mathcal{H}_r(\Omega)$,*

$$||u||_{\infty,\Omega} \leq C_1 ||u||_{(2,r);\Omega}.$$

(2) *For any integer $i \geq 2$, there is a constant $C_i > 0$ depending only on $\epsilon$ and $\Delta$, such that for any $u \in \mathcal{H}_i(\Omega)$:*

$$||u||_{(2,i);\Omega'} \leq C_i \left( ||u||_{2,\Omega'} + ||\Delta u||_{(2,i-2);\Omega} \right).$$

*Proof of Lemma 3.4.9.* If we write $\mathcal{H}^0(\Omega)$ for the subspace of $\mathcal{H}(\Omega)$ obtained by completing $C_c^\infty(\Omega)$, then the first (resp. second) assertion of this lemma is proved on pages 148–151 (resp. 262–267) of [GiT] (resp. [Fo]). Moreover, given any $u$ in $\mathcal{H}_r(\Omega)$, for any $r$, and a subset $\Omega'$ of nonempty interior with its closure contained in the interior of $\Omega$, we can find a smooth cut-off function $\psi$ such that $v = \psi u$ is in $\mathcal{H}_r^0(\Omega)$ and $v = u$ on $\Omega'$ (see, for example, [Fo, pp. 275–276] ), and we can do this uniformly for all $r$. Since we can bound $||\psi u||_{(2,r),\Omega}$ by a universal constant times $||u||_{(2,r),\Omega}$, we get the first part of the lemma, and for the second part we only need to bound $||\Delta(\psi u)||_{(2,i-2),\Omega}$ by a constant times $||\Delta(u)||_{(2,i-2),\Omega} + ||u||_{(2,i-1),\Omega}$. We will suppress $\Omega$ henceforth. Recall that $||\Delta(\psi u)||_{(2,i-2),\Omega}$ is (by definition) $\sum_{|\nu|\leq i-2} ||D^\nu \Delta(\psi u)||$, which is $\sum_{|\nu|\leq i-2} ||D^\nu((\Delta\psi)u + (\nabla\psi)(\nabla u) + \psi(\Delta u))||_2$; this can be bounded by

$$\sum_{|\nu|\leq i-2} \sum_{\beta+\gamma=\nu} (||D^\beta(\Delta\psi)D^\gamma(u)||_2 + \ldots).$$

It is easy to estimate derivatives of $\psi$, and moreover, all but one of the terms involves at most $i-1$ derivatives of $u$ and is hence bounded by a constant times $||u||_{(2,i-1)}$. The remaining term involves $i-2$ derivatives of $\Delta u$, and this is bounded by a constant times $||\Delta u||_{(2,i-2)}$. The constants involved are controllable in terms of $\epsilon$. □

*Proof of Lemma 3.4.8.* Fix a $g_0$ and identify an open neighborhood of $g_0$ containing $V$ with an open neighborhood of $\mathbb{R}^N$, $N = \dim(Y)$, contained in the unit ball of radius $\epsilon$. By shrinking $V$, we may assume that $E(f,s)$ is in $\mathcal{H}_i(V)$ for all $i$. By the first part of Lemma 3.4.9, we are left to estimate $||E(f_s)||_{(2,r),V}$ for an $r > N/2$, which we will take to be of the form $2m$. Choose neighborhoods $V_j$, $0 \leq j \leq m-3$, of $g$ such that we have the nested inclusions

$$V' = V_{m-3} \subset \ldots V_1 \subset V_0 = V,$$

such that the closure of each $V_i$ is in the interior of $V_{i+1}$. Apply the second part of Lemma 3.4.9 succesively with $u = E(f_s)$ and

$$(i, \Omega', \Omega) = (2(m-j), V_{m-j}, V_{m-j+1})$$

with $j$ ranging over $\{0, \ldots, m-3\}$, repeatedly making use of the fact that $\Delta E(f_s)$ is $\lambda_s E(f_s)$. The assertion of Lemma 3.4.8 follows at $g_0$. Now if $g$ is any other point, we can translate the neighborhood $V$ to $g$ and see that the $G_\infty$-invariance of $\Delta$ gives the corresponding result with the same constant. □

Let $T \in \operatorname{Lie} A_\mathbb{R}$ be a sufficiently regular parameter, i.e., with $\langle \nu, T \rangle$ being sufficiently large for each $\nu$ in $\Delta$. Denote by $\wedge^T E(f_s)$ the Arthur truncation of $E(f_s)$ (cf. [A1] or [MW1, p. 35]). One knows that $\wedge^T E(f_s)$ decreases rapidly at infinity, and is hence in $L^p(Y)$ for all $p$. Clearly, $\|\psi\|_{2,C_g}$ is bounded above by $\|\psi\|_{2,G}$ for any square-integrable function $\psi$ on $Y$. (Abusing notation, we write $\|.\|_{p,G}$ instead of $\|.\|_{p,Y}$.) So, in view of Lemma 3.4.8, Proposition 3.4.6 (and the theorem) will be proved once we establish the following:

PROPOSITION 3.4.10. *Fix $\Xi, \mathcal{S}$ as above. Then there exist $N \in \mathbb{Z}$ and functions $\Phi_1(s), \Phi_2(s)$ of bounded order such that, for all $g$ in $\mathcal{S}$,*

(i) $$\|\wedge^T E(f_s)\|_{2,G} \leq |\Phi_1(s)|$$

*and, up to shrinking $V = V_g$ to a smaller neighborhood of $g$ with nonempty interior,*

(ii) $$\|E(f,s) - \wedge^T E(f,s)\|_{\infty,V} \leq |\Phi_2(s)| \|g\|^N.$$

Our first object will be to prove the bound (i) of this proposition. For the sake of an inductive argument later, we will in fact need to analyze the $L^2$-norm of the truncation of a class of *relative* Eisenstein series $E^{P'}(h_s^w)$ (see (3.4.15) below) attached to any standard parabolic subgroup $P' = M'U'$ and a Weyl element $w$ representing a double coset in $W_{M'}\backslash W/W_M$. We set

$$h_s^w = M_w(f_s),$$

where $M_w$ is the intertwining operator defined by

$$M_w(f_s)(g) = \int_{(U_0 \cap w\bar{U}_0 w^{-1})(\mathbb{A}_F)} f_s(w^{-1}ug)\, du.$$

Here $\bar{U}_0$ denotes the opposite unipotent subgroup of the maximal unipotent subgroup $U_0$.

*Remark* 3.4.11. Since $f_s$ is good, it factorizes as $\prod_v f_{v,s}$ and at any $v$ where $\omega$ is unramified, $f_{v,s}$ is the standard function $\phi_{v,s}$ and $M_w(\phi_{v,s}) = c_v^w(s)\phi_{v,s}^w$, where $c^w(s) = \prod_v c_v^w(s)$ is a ratio of *abelian L-functions*; see [Ik1], [PS-R1]. Moreover, the possible poles are contained (cf. [Ik1]) in the set $\{\frac{3}{2}, 1, 0, -\frac{1}{2}\}$.

Recall that every element $w$ in the Weyl group $W$ is determined by what it does to $\{e_1, e_2, e_3\}$, and that each $e_i$ gets sent to $\pm e_j$, for some $j$. We will adopt the following notation: If $w$ sends $e_1$ to $e_3$, $e_2$ to $-e_2$ and $e_3$ to $-e_1$, for example, we will denote $w$ by $[3, \bar{2}, \bar{1}]$. Put

(3.4.12) $$\Sigma_M = \{w_1, w_2, w_3, w_4, w_5, w_6, w_7, w_8\}$$

with $w_1 = \text{id} = [1, 2, 3]$, $w_2 = [\bar{2}, \bar{3}, 1]$, $w_3 = [\bar{1}, \bar{2}, \bar{3}]$, $w_4 = [1, 2, \bar{3}]$, $w_5 = [1, 3, \bar{2}]$, $w_6 = [1, \bar{3}, \bar{2}]$, $w_7 = [2, \bar{1}, 3]$ and $w_8 = [2, \bar{3}, 1]$.

For each (standard) parabolic $P' = M'U' \neq G$, define a subset $\Sigma_{M,M'}$ of $\Sigma_M$ by the following table:

(3.4.13)

| $\Delta_{M'}$ | $\Sigma_{M,M'}$ |
|---|---|
| $\{\alpha, \beta\}$ | $\{w_1, w_2, w_4, w_6\}$ |
| $\{\alpha, \gamma\}$ | $\{w_1, w_2, w_5\}$ |
| $\{\beta\gamma\}$ | $\{w_1, w_7\}$ |
| $\{\alpha\}$ | $\{w_1, w_2, w_3, w_4, w_5, w_6\}$ |
| $\{\beta\}$ | $\{w_1, w_3, w_4, w_6, w_7, w_8\}$ |
| $\{\gamma\}$ | $\{w_1, w_2, w_5, w_7\}$ |
| $\phi$ | $\Sigma_M$ |

LEMMA 3.4.14. *Let $P' = U'M'$ be a standard parabolic. Then $\Sigma_{M,M'}$ is a set of representatives for $W_{M'}\backslash W/W_M$, where $W_{M'}$ (resp. $W_M$) denotes the Weyl group of the roots relative to $M'$ (resp. $M$). Moreover, for any $w \in \Sigma_{M,M'}$, $Q'_w := wPw^{-1} \cap M'$ is a parabolic subgroup of $M'$, $P'_w := Q'_w \cdot U'$ (but not $wPw^{-1} \cap P'$) is a (standard) parabolic subgroup of $G$ contained in $P'$, and if $U'_w = wPw^{-1} \cap U'$, then $U' = U'_w \cdot V'_w$, where $V'_w = w\bar{U}_0 w^{-1} \cap U_0$.*

*Proof of Lemma* 3.4.14. Put

$$W^*_{M,M'} = \left\{ w \in W \;\middle|\; \begin{array}{l} w(\nu) > 0, \quad \text{for all } \nu \in \Delta_M; \\ w^{-1}(\nu) > 0, \quad \text{for all } \nu \in \Delta_{M'} \end{array} \right\}.$$

Then it is well-known that $W^*_{M,M'}$ is a set of representatives in $W$ for $W_{M'}\backslash W/W_M$. Explicitly, one has, when $\Delta_{M'}$ is empty,

$$W^*_{M,M'} = \left\{ \begin{array}{c} w_1, \bar{w}_2 := [3, \bar{2}, \bar{1}], \, \bar{w}_3 := [\bar{3}, \bar{2}, \bar{1}], \, w_4, w_5, \\ w_6, \bar{w}_7 := [2, \bar{3}, \bar{1}], \, w_8 \end{array} \right\}.$$

The only difference between this set and the $\Sigma_M$ defined above is that $w_2$, $w_3$, $w_7$ are replaced by $\bar{w}_2, \bar{w}_3, \bar{w}_7$. But it is easy to see that for $j = 2, 3, 7$, $\bar{w}_j = w_j w'_j$, with $w'_j \in W_M$. Thus, when $\Delta_{M'}$ is empty, $\Sigma_{M,M'}$ is also a set of representatives for $W_{M'}\backslash W/W_M$ in $W$. The other cases of $M'$ can be checked by comparing with the table on page 21 for $\Sigma_{M,M'}$. The remaining assertions of the lemma follow by a straightforward computation.

The relative Eisenstein series of interest to us is attached to $(P', w)$ and is given by

$$E^{P'}(h_s^w)(g) = \sum_{\gamma \in Q'_w(F) \backslash M'(F)} h_s^w(\gamma g). \tag{3.4.15}$$

Denote by $\wedge^{T,P'} E^{P'}(h_s^w)$ the truncation of $E^{P'}(h_s^w)$ as defined on page 97 of [A1]. Let $\langle \wedge^{T,P'} E^{P'}(h_s^w), \wedge^{T,P'} E^{P'}(h_s^w) \rangle_{P',T}$ denote the scalar product defined on page 45 of [A2], which specializes to the usual scalar product when $P' = G$.

PROPOSITION 3.4.16. $\langle \wedge^{T,P'} E^{P'}(h_s^w), \wedge^{T,P'} E^{P'}(h_s^w) \rangle_{P',T}$ is of bounded order.

*Proof of Proposition* 3.4.16. To mesh better with the notation of Langlands and Arthur we will shift the parameter $s$ by $1/2$ so that the unitary axis is the line $\Re(s) = 0$ (instead of $\Re(s) = 1/2$), but by abuse of notation we will continue to denote the Eisenstein series by $E^{P'}(h_s^w)$.

For $s, s'$ in $\mathbb{C}$, set

$$h(s, s') = \langle \wedge^{T,P'} E^{P'}(h_s^w), \wedge^{T,P'} E^{P'}(h_{-\overline{s'}}^w) \rangle_{P',T}, \tag{3.4.17}$$

which is a meromorphic function in $(s, s')$ whose poles are in a union of vertical lines $s = c_i$ and horizontal lines $s' = c'_j$, depending on the polar set of $E^{P'}(h_s^w)$, which is in turn governed (cf. [La4]) by the constant term along $P_0$. The intersection of this polar set of $h(s, s')$ with $\{\overline{s} + s' = 0\}$ has only a *finite* number of points in $\Xi$. So, for large enough $t_0$, the function $s \to h(s, -\overline{s})$ is well-defined and finite on $\{s = x + iy \in \Xi \mid |y| \geq t_0\}$.

By Langlands we can write $E^{P'}(h_s^w)$ in terms of residues of *cuspidal* Eisenstein series. In the notation of [A2, p. 43], it is a finite sum of terms of the form

$$\text{Res}_{\Lambda \to \Lambda_S}(F_B(\Lambda), \Lambda + \beta s), \tag{3.4.18}$$

where $S$ is a suitable flag of affine subspaces $\mathfrak{t} = \Lambda_\mathfrak{t} + \text{Lie}\,(A_{P'_w})^*$ in $\text{Lie}\,(A_B)^*$, $\beta$ in $\text{Lie}\,(A_{P'_w})^*$, $B = M_B U_B$ a relevant parabolic subgroup of $P'_w$, and $F_B$ a meromorphic function of $\Lambda$ on $\text{Lie}\,(A_B)^*$ with singularities lying in hyperplanes.

Note that $h_s^w$ is associated to the representation induced from $P'_w$ by a 1-dimensional representation, which is residual on $M'_w$ associated to a representation induced by a character $\xi_\Lambda^w$ of the minimal parabolic $P_0$, which is exponential in $\Lambda$ and of bounded order in each coordinate direction. The construction of Langlands ([La4]) is natural relative to induction in stages, and we see that $E^{P'}(h_s^w)$ must be the residue of the cuspidal Eisenstein series defined by $\xi_\Lambda^w$. Hence the only relevant $B$ occurring in (3.4.18) is $P_0$. In this case the function $F_B(\Lambda)$ is right $K$-finite in the space of the representation induced by $\xi_\Lambda^w$. By the Iwasawa decomposition we see then that $F_B(\Lambda)$ is of exponential type.

By Lemma 3.1 of [A2], we see that $h(s, s')$ is a finite sum of double residues

$$(3.4.19) \quad \mathrm{Res}_{\Lambda \to \Lambda_S} \mathrm{Res}_{\Lambda' \to \Lambda_{S'}} \omega^{T,P'}(\Lambda + \beta s, -\Lambda' + \beta s', F_B(\Lambda), F_B(\Lambda')),$$

where (cf. [A2, p. 46])

$$\omega^{T,P'}(\lambda, \lambda', \phi, \phi') = \sum_{P_1 \subset P'} \sum_{t,t' \in W^{P'}(B,P_1)} \frac{\langle M_{t,\lambda}(\phi), M_{t',-\overline{\lambda'}}(\phi') \rangle e^{(t\lambda - t'\lambda')(T)}}{a(P_1, P') \prod_{\nu \in \Delta_{P_1}^{P'}}(t\lambda - t'\lambda')(\nu^\vee)}.$$

Here $a(P_1, P')$ is a nonzero constant, and $W^{P'}(B, P_1)$ is the set of isomorphisms, possibly empty, of Lie $A_B$ onto Lie $A_{P_1}$ leaving Lie $A_{P'}$ pointwise fixed which arises from the restricted Weyl group $W(G, A_0)$. (In our particular case, $P_1$ has to be $B = P_0$.) We refer to [A2] for a definition of $\Delta_{P_1}^{P'}$.

The $c$-functions associated to the intertwining operators $M_{t,\lambda}$ (of (3.4.19)) acting on functions $\phi(\lambda)$ coming from induced representations from the Borel subgroup are, just like the $c_s^w$ associated to $M_w$ discussed in Remark 3.4.11, ratios of abelian $L$-functions ([La2]), hence bounded of order 1. Combining this with that fact that $F_B(\Lambda)$ is of pure exponential type, we see that the expression $\omega^{T,P'}(\Lambda + \beta s, -\Lambda' + \beta s, F_B(\Lambda), F_B(\Lambda'))$ is of bounded order in $\Lambda$ and $\Lambda'$. Its singularities are in hyperplanes, and evaluating the residues for $\Lambda \to \Lambda_S$ and $\Lambda' \to \Lambda_{S'}$, we see that the double residue of (3.4.19) must be of bounded order in $(s, s')$ along any direction. This then implies the result we want about $h(s, -\overline{s})$ in $\{s = x + iy \in \Xi \mid |y| \geq t_0\}$. $\square$

*Proof of Proposition* 3.4.10. Part (i) of this proposition clearly follows by application of Proposition 3.4.16 in the special case $P' = G$. It remains to prove (ii). First we need to derive a manageable expression for $E(f_s) - \wedge^T E(f_s)$. This takes some work and we start with some preliminaries.

For each standard parabolic $P' = M'U'$, $M' = \mathrm{Cent}_G(A')$, put

$$M_0' = \left\{ m \in M'(\mathbb{A}_F) \;\middle|\; \begin{array}{l} |\chi(m)| = 1, \\ \text{for all } \chi \in \mathrm{Hom}(M', \mathbb{G}_m) \end{array} \right\}.$$

Then there exists a map

$$\log_{M'} : M'(\mathbb{A}_F) \longrightarrow \mathrm{Lie}(A')_{\mathbb{C}}$$

defined in order to satisfy $\langle \chi, \log_{M'}(m) \rangle = \log |\chi(m)|$, for all $m \in M'(\mathbb{A}_F)$ and $\chi \in \mathrm{Hom}(M', \mathbb{G}_m)$. Then $M_0'(\mathbb{A}_F)$ is precisely $\ker(\log_{M'})$. On the other hand, the Iwasawa decomposition $G(\mathbb{A}_F) = U'(\mathbb{A}_F) M'(\mathbb{A}_F) K$ leads to a map

$$m_{P'} : G(\mathbb{A}_F) \longrightarrow M_0'(\mathbb{A}_F) \backslash M'(\mathbb{A}_F),$$

given by sending $g$ to the coset $M'_0(\mathbb{A}_F)m$, if $g = umk$, with $u \in U'(\mathbb{A}_F)$, $m \in M'(\mathbb{A}_F)$ and $k \in K$. Thus $\log_{M'} \circ m_{P'}$ is a map from $G(\mathbb{A}_F)$ into $\text{Lie}(A')_{\mathbb{C}}$.

Note that $A_\Delta$ is the center $Z$ of $G$. We have a decomposition $\text{Lie } A'_{\mathbb{R}} = \text{Lie } Z_{\mathbb{R}} \oplus \sigma^G_{M'}$, where $\sigma^G_{M'}$ is the subspace generated by the set $\Delta_{M'}$ of simple roots relative to $M'$. Denote by $\hat{\tau}_{P'}$ the characteristic function of the sum of $\text{Lie } Z_{\mathbb{R}}$ and the interior of the cone of $\sigma^G_{M'}$ generated by $\Delta_{M'}$. Then one has (see [MW1, p. 35]) (for all $G \in G(\mathbb{A}_F)$):

$$(3.4.20) \quad (E(f_s) - \wedge^T E(f_s))(g) = \sum_{P_0 \subseteq P' \subsetneq G} (-1)^{(\dim A' - \dim A_0)}$$

$$\times \sum_{\gamma \in P'(F)\backslash G(F)} \hat{\tau}_{P'}(\log_{M'}(m_{P'}(\gamma g) - T') E(f_s)_{P'}(\gamma g),$$

where $T'$ is the projection of $T$ in $\text{Lie } A'_{\mathbb{R}}$, and $E(f_s)_{P'}$ is the *constant term* of $E(f_s)$ along $P'$, namely:

$$(3.4.21) \qquad E(f_s)_{P'}(g) = \int_{U'(F)\backslash U'(\mathbb{A}_F)} E(f_s)(ug)\, du.$$

PROPOSITION 3.4.22. (a) *For each $P'$, the constant term of $E(f_s)$ along $P'$ is given by*

$$E(f_s)_{P'}(g) = \sum_{w \in \Sigma_{M,M'}} E^{P'}(h_s^w)(g)$$

*for all $g$ in $G(\mathbb{A}_F)$.*

(b) *For $T$ large enough (and regular),*

$$E(f_s) - \wedge^T E(f_s) = \sum_{P_0 \subseteq P' \subsetneq G} (-1)^{(\dim A' - \dim A_0)}$$

$$\times \sum_{w \in \Sigma_{M,M'}} \hat{\tau}_{P'}(\log_{M'}(m_{P'}(g) - T') E^{P'}(M_w(f_s))(g).$$

*Remark.* For $P' = P_0$, the expression given in part (a) is well-known by the work of Rallis and Piatetski-Shapiro ([PS-R1]); see also [Ik1]. For *cuspidal* Eisenstein series (which is *not* what we have here), there is a general formula for any reductive group $G$; see [MW1, §II.1.7], for example.

*Proof of Proposition* 3.4.22. (a) Combining (3.4.3) and Lemma 3.4.14, and remembering that $U'_w$ denotes $wPw^{-1} \cap P'$, we get

$$E(f_s)_{P'}(g) = \sum_{w \in \Sigma_{M,M'}} \sum_{\gamma \in Q'_w(F)\backslash M'(F)} \int_{U'_w(F)\backslash U'(\mathbb{A}_F)} f_s(w^{-1}u\gamma g)\, du.$$

Now we use the decomposition $U' = U'_w \cdot V'_w$ of Lemma 3.4.14, and recall the $P(F)$-invariance of $f_s$ together with the fact that the intertwining operator $M_w$ is defined (see (3.4.12)) by integration over $V'_w(\mathbb{A}_F)$. We see that, up to normalizing the measures so that volume of the compact group $(wUw^{-1} \cap U')(F) \backslash (wUw^{-1} \cap U')(\mathbb{A}_F)$ is 1,

$$E(f_s)_{P'}(g) = \sum_{w \in \Sigma_{M,M'}} \sum_{\gamma \in Q'_w(F) \backslash M'(F)} M_w(f_s)(\gamma g)\, du.$$

This proves (a).

(b) By definition, $Q'_w \backslash M' = P'_w \backslash P'$. Since $m_{P'}$ is left invariant under elements of $P'(F)$, the assertion is a consequence, in view of the identity of (a), of the following:

LEMMA 3.4.23. *For $g$ lying in a fixed Siegel domain, for $T$ sufficiently large (and regular):*

$$\hat{\tau}_{P'}(\log_{M'}(\gamma g)) - T') E^{P'}_w(f_s(\gamma g)) = 0 \quad \text{if} \quad \gamma \in G(F) - P'(F).$$

*Proof of Lemma* 3.4.23. We may restrict ourselves to a Siegel domain of the form $\mathfrak{S}_c = \Omega A'_c K$, where $\Omega$ is a compact subset of $P'(\mathbb{A}_F)$, $c$ a positive real number, and $A'_c$ the set of elements $t$ in $A'(\mathbb{A}_F)$ such that $|\nu(t)| > c$ for each $\nu$ in $\Delta_{M'}$. As usual, one chooses $c$ to be small enough and $\Omega$ large enough so that $G(\mathbb{A}_F) = P'(F) \mathfrak{S}_c$.

Let $\gamma \in G(F) - P'(F)$. Then by the Bruhat decomposition, $\gamma = \bar{u}w$, for some Weyl element $w$, $\gamma \in P'(F)$ and $\bar{u} \in \bar{U}'_w(F) := wU'(F)w^{-1}$.

By the Iwasawa decomposition relative to $P'$, we can write any $g$ in $G(\mathbb{A}_F)$ as a product $uxtk$, with $u \in U'(\mathbb{A}_F)$, $x \in M'_0(\mathbb{A}_F)$, $t \in A'(\mathbb{A}_F)$ and $k \in K$. For any root $\nu$ (in $\Phi_{M'}$) let us write $H_\nu(g)$ for $|\nu(t)|$, which is well-defined. By the Iwasawa decomposition relative to $P_0$, we can also write $g = nt_0 k_1$ with $n \in U_0(\mathbb{A}_F)$, $t_0 \in A_0(\mathbb{A})$ and $k \in K$. Then $H_\nu(g)$ identifies with $\nu$ evaluated at the projection of $t_0$ in $A'(\mathbb{A}_F)$. For $w \in W$, we will write $\nu^w(t) = \nu(t^w) := \nu(wtw^{-1})$.

*Claim* 3.4.24. Let $\bar{u} \in \bar{U}'_w(F)$. Then there exists $\nu \in \Phi^+_{M'}$ such that $\nu^w = \nu^{-1}$ on $A'$ and $H_\nu(\bar{u}) \leq 1$.

Suppose the claim holds. Put $\gamma = y\bar{u}w \in G(F)$. Then for any $g = uzk \in G(\mathbb{A}_F)$, with $z \in U'(\mathbb{A}_F)M'(\mathbb{A}_F)$, $u \in U'(\mathbb{A}_F)$, $k \in K$, we have

$$H_\nu(\gamma g) = H_\nu(\bar{u}wuzk) = H\nu(\bar{u}u^w z^w) = H_\nu(z^w) H_\nu((z^w)^{-1}(\bar{u}u^w)z^w).$$

Note that $H_\nu(z^w) = H_{\nu^w}(z) = H_\nu(g)^{-1} < c^{-1}$. On the other hand, since the conjugate of $\bar{u}u^w$ by $z^w$ lies in $\bar{U}'_w(F)$, we have by the claim, $H_\nu((z^w)^{-1}(\bar{u}u^w)z^w) \leq 1$. So $H_\nu(\gamma g) < c^{-1}$. Since $\nu$ is a positive sum of roots in $\Delta_{M'}$, we can take $T$ to be large enough so that $\nu(T') > c^{-1}$. Then we must have $\hat{\tau}_{P'}(\log_{M'}(\gamma g)) - T') = 0$.

It remains to prove the claim. We can define local "heights" $H_{\nu,v}$ in the obvious way to get $H_\nu(g) = \prod_v H_{\nu,v}(g_v)$. So it suffices to find a $\nu$ with $\nu^w = \nu^{-1}$ such that $H_{\nu,v}(\bar{u}_v) \leq 1$ everywhere. It is instructive to begin by looking at the situation for the basic case $G' = \mathrm{SL}_2$, and the argument for $G$ is simply a jazzed up version. Recall that at each $v$, the standard maximal compact subgroup, if $G'(F_v)$ is $G'(O_v)$, of $F$ is $p$-adic, SO(2) if $F_v = \mathbb{R}$, and SU(2) if $F_v = \mathbb{C}$. Let $\bar{u} = \begin{pmatrix} 1 & 0 \\ x & 1 \end{pmatrix}$, with $x \in F^*$. By the Iwasawa decomposition, we can find $a \in F_v^*$ and $y \in F_v$ such that

$$k := \begin{pmatrix} a^{-1} & 0 \\ 0 & a \end{pmatrix} \begin{pmatrix} 1 & y \\ 0 & 1 \end{pmatrix} \begin{pmatrix} 1 & 0 \\ x & 1 \end{pmatrix} = \begin{pmatrix} a^{-1}(1+xy) & a^{-1}y \\ ax & a \end{pmatrix}$$

lies in $K_v$. In the non-archimedean case, if $|x| \leq 1$, we take $a = 1$ (and $y = 0$) and if $|x| > 1$, we are forced to have $|a| \leq |x|^{-1}$. Suppose $F_v = \mathbb{R}$. Then $k^t k$ must be the identity matrix, and this forces $|a| = \frac{1}{\sqrt{x^2+1}}$ (and $y = -x$). The complex case is similar and we get $|a| = \frac{1}{\sqrt{|x|^2+1}}$. Now consider the global case and write $a$ for the idele $(a_v)$ with $a_v$ being (for each $v$) the element chosen above (and denoted just by $a$ for brevity of notation). Then we get

$$|a| \leq \prod_{v|\infty} \frac{1}{\sqrt{1+|x|_v^2}} \prod_{v \text{ finite}} |x|_v^{-1} \leq \prod_{v|\infty} \frac{|x|_v}{\sqrt{1+|x|_v^2}} \leq 1,$$

by Artin's product formula, which gives $\prod_{\text{all } v} |x|_v^{-1} = 1$ (since $x \in F^*$). Take $\nu$ to be the root $\begin{pmatrix} a & 0 \\ 0 & b \end{pmatrix} \mapsto a/b$. Clearly, $\bar{U} = w \begin{pmatrix} 1 & x \\ 0 & 1 \end{pmatrix} w^{-1}$ and $\nu^w = \nu^{-1}$, with $w = \begin{pmatrix} 0 & -1 \\ 1 & 0 \end{pmatrix}$. Moreover, since $|a| \leq 1$, we have $H_\nu(\bar{u}) \leq 1$.

In our case ($G = \mathrm{GSp}(6)$), we will now give the argument for the most basic case $P' = P_0$ and $\bar{U}_w = \bar{U}_0$. The other cases can be easily extracted from this. We write

$$\bar{u} = \begin{pmatrix} 1 & 0 & 0 & & & \\ t_1 & 1 & 0 & & 0 & \\ t_3 & t_2 & 1 & & & \\ \hline u_1 & u_2 & u_3 & 1 & -t_1 & t_1 t_2 - t_3 \\ v_1 & v_2 & v_3 & 0 & 1 & -t_2 \\ w_1 & w_2 & w_3 & 0 & 0 & 1 \end{pmatrix}$$

with $u_i, v_i, w_i, t_i$ in $F$. Let $v$ be a place. By the Iwasawa decomposition (over $F_v$), we can find $t = \mathrm{diag}(a, b, c, a^{-1}, b^{-1}, c^{-1}) \in A_0(F_v)$ and $n \in U_0(F_v)$ such that $k := t^{-1}n^{-1}u$ lies in $K'_v$, the standard maximal compact subgroup of $\mathrm{Sp}(6, F_v)$. If we write

$$n^{-1} = \left( \begin{array}{ccc|ccc} 1 & s_1 & s_3 & z_1 & z_2 & z_3 \\ 0 & 1 & s_2 & y_1 & y_2 & y_3 \\ 0 & 0 & 1 & x_1 & x_2 & x_3 \\ \hline & & & 1 & 0 & 0 \\ & 0 & & -s_1 & 1 & 0 \\ & & & -s'_3 & -s_2 & 1 \end{array} \right), \text{ with } s'_3 = s_3 - s_1 s_2,$$

then we have:

(3.4.25) $$k = \left( \begin{array}{c|c} * & * \\ \hline X & Y \end{array} \right)$$

where

$$X = \begin{pmatrix} au_1 & au_2 & au_3 \\ b(-s_1 u_1 + v_1) & b(-s_1 u_2 + v_2) & b(-s_1 u_3 + v_3) \\ c(-s'_3 u_1 - s_2 v_1 + w_1) & c(-s'_3 u_2 - s_2 v_2 + w_2) & c(-s'_3 u_3 - s_2 v_3 + w_3) \end{pmatrix}$$

and

$$Y = \begin{pmatrix} a & -at_1 & -at'_3 \\ -bs_1 & b(s_1 t_1 + 1) & b(s_1 t'_3 - t_2) \\ -s'_3 c & c(s'_3 t_1 - s_2) & c(s'_3 t'_3 + s_2 t_2 + 1) \end{pmatrix}.$$

In the non-archimedean case, we are then forced to have: $|a|, |b|, |c| \leq 1$. Moreover, suppose some element, call it $d$, of the set $\{u_1, u_2, u_3, t_1, t'_3\}$ is nonzero. Then we can take $a = 1$ if $d$ is integral at $v$ and we are forced to have $|a| \leq |d|^{-1}$ if $d$ is nonintegral at $v$. In the archimedean case, just as in the SL(2) case, we have $a = \frac{1}{\sqrt{1+|x|_v^2}}$. Again the global $a$ has absolute value bounded above by 1. In this case, we are done by taking $\nu = 2e_1$ and $w$ to be the longest root. Suppose every element of this set is zero. then $t_3$ is also zero, and we get

$$\bar{u} = \left( \begin{array}{ccc|ccc} 1 & 0 & 0 & & & \\ 0 & 1 & 0 & & 0 & \\ 0 & t_2 & 1 & & & \\ \hline 0 & 0 & 0 & 1 & 0 & 0 \\ v_1 & v_2 & v_3 & 0 & 1 & -t_2 \\ w_1 & w_2 & w_3 & 0 & 0 & 1 \end{array} \right) k_1,$$

with $k_1$ in $K'_v$. Now we may replace $\bar{u}$ by $\bar{u} k_1^{-1}$. Then we get a condition on $b$ analogous to the one on $a$ above if some element of the set $\{v_1, v_2, v_3, s_1, t_2\}$ is nonzero. In this case we take $\nu$ to be $2e_2$. If every element of this set is also zero, then since $\bar{u}$ is nonzero, some $w_j$ must be nonzero, which leads to a condition on $|c|$. We are done by taking $\nu$ to be $2e_3$ in this case. This finishes the proof of the Claim, Lemma 3.4.23 and Proposition 3.4.22.

*Proof of Proposition* 3.4.10 (*contd.*). In view of part (b) of Proposition 3.4.22, it suffices to show that, for all $P' = M'U'$ and $w$ in $W^*_{M,M'}$, we have the bound

$$\|E^{P'}(h_s^w)\|_{2,V'_g} \leq |\Phi(s)| \|g\|^N, \tag{3.4.26}$$

for some function $\Phi(s)$ of bounded order and integer $N$.

Since $Y$ has finite volume, say $A$, which is bigger than or equal to $\mathrm{vol}(U_g)$, we get the bound

$$\|E^{P'}(h_s^w)\|_{2,V'_g} \leq A \|E^{P'}(h_s^w)\|_{\infty,V'_g}. \tag{3.4.27}$$

Let us call $(P', w)$ *convenient* if $P'_w = P'$. This happens in two cases: (i) when $P' = P_0$ and $w$ is arbitrary, and (ii) when $P' = P$ and $w = w_1$ (identity). Suppose we are in either of these cases. Then $E^{P'}(h_s^w) = M_w(f_s)$. Recall that the effect of the intertwining operator $M_w$ on $f_s$ is given by a ratio of abelian $L$-functions (see Remark 3.4.11), which are functions of bounded order by Hecke. Moreover, the definition of $f_s$ shows that its value at $g$ is bounded by a constant (independent of $s$) times $\|g\|^r$, for some integer $r$ depending on the width of $\Xi$. Making use of (3.4.27), and possibly shrinking $V_g$ slightly, we can find an $N$ for which (3.4.26) holds in such a convenient case.

Now suppose $(P', w)$ is *inconvenient*. Then $E^{P'}(M_w(f_s))$ is left-invariant under $U'(\mathbb{A}_F)M'(F)$. Let $A'$ denote the center of $M'$ and put $K' = K \cap M(\mathbb{A}_{F,f})$. Since $P'$ commutes with $A'$, we get the requisite behavior of $E^{P'}(M_w(f_s))$ along the root directions along $A'$. It remains to study the growth along the root directions which are trivial on $A'$. By the decomposition $G(\mathbb{A}_F) = P'(\mathbb{A}_F)K_\infty K$, the problem reduces, in view of (3.4.27), to bounding the growth, for each $k \in K_\infty/(K_\infty \cap M'_\infty)$, of the function on $Y' := M'(F)(Z'_\infty)^+ \backslash M'(\mathbb{A}_F)/K'$ given by $g \to E^{P'}(h_s^w)(gk)$. Let $W_g$ denote the corresponding neighborhood (of the image of $g$), which we may shrink to be in the open ball of radius $\epsilon$ in the local coordinates, and let $W'_g$ be the subset defined by $U'_g$ such that $\overline{W'_g} \subset W_g$. Note that the relative Eisenstein series is an eigenfunction of the Laplacian $\Delta'$ defined by the Casimir operator on $M'(F_\infty)$ with eigenvalue $\lambda'_s$, which is again a quadratic function of $s$. Applying the earlier argument (see Lemmas 3.4.8 and 3.4.9), and making use of the $K_\infty$-finiteness of $f_s$, we get for a polynomial $P'(s)$ independent of $g$ and $k$:

$$\|E^{P'}(h_s^w)\|_{\infty,W'_g} \leq |P'(s)| \cdot \|E^{P'}(h_s^w)\|_{2,W'_g}. \tag{3.4.28}$$

We can truncate $E^{P'}(h_s^w)$, and in view of Proposition 3.4.16, it is enough to estimate $\|E^{P'}(h_s^w) - \wedge^{T,P'}E^{P'}(h_s^w)\|_{2,W'_g}$. In view of the obvious analog of Proposition 3.4.22, it will suffice to bound the function $\|E^{P''}(H_s^u)\|_{\infty,W'_g}$ for

any parabolic $P''$ contained in $P'$, Weyl element $u$ in $W_{M'_w}\backslash W_{M'}/W_{M''}$ and function $H_s^u = M_u(h_s^w)$. In the convenient case the proof goes as above. In the inconvenient cases, the parabolic taking the place of $P'_w$ is necessarily the minimal parabolic $P_0$. So if we truncate one more time, all the resulting cases will be convenient. □

This finishes the proof of Propositions 3.4.10, 3.4.6, and hence Theorem 3.4.1.

*Remark* 3.4.29. There is clearly an inductive argument buried here which should help us understand the nature of integral representations of $L$-functions for larger groups.

3.5. *Modularity in the good case.* Let $F$ be a totally complex number field. We will call a pair $(\pi, \pi')$ of cuspidal automorphic representations of $\mathrm{GL}(2, \mathbb{A}_F)$ *good* if at *each* finite place $v$, either $\pi_v$ or $\pi'_v$ is not supercuspidal.

THEOREM 3.5.1. *For every good pair $(\pi, \pi')$, there exists an isobaric automorphic representation $\pi \boxtimes \pi'$ of $\mathrm{GL}(4, \mathbb{A}_F)$ such that the identities $(L_v), (\varepsilon_v), (L_\infty)$ asserted in Theorem* M *hold, where $v$ is an arbitrary finite place.*

*Proof.* In view of Lemma 3.1.1 we may, and we will, assume that neither $\pi$ nor $\pi'$ is dihedral, i.e., automorphically induced from a grossencharacter of a quadratic extension, and that $\pi'$ is not a character twist of $\pi$.

Our first object is to define, for each place $v$, an irreducible admissible representation $\Pi_v$ of $\mathrm{GL}(4, F_v)$ canonically associated to the pair $(\pi_v, \pi'_v)$. As before, let $\sigma_v$ (resp. $\sigma'_v$) be the 2-dimensional representation of $W'_{F_v}$ attached to $\pi_v$ (resp. $\pi'_v$), as in [Ku]. First let $v$ be archimedean. Then the local Langlands correspondence exists, and we may (and will) take $\Pi_v$ to be the representation associated to $\sigma_v \otimes \sigma'_v$ by the method of [La 1].

Next consider the non-archimedean case. Here we will, for every $n \geq 1$ with $d$ a divisor of $n$ and a supercuspidal representation $\beta$ of $\mathrm{GL}(d, F_v)$, denote the associated (generalizd) special representation of $\mathrm{GL}(n, F_v)$ by $\mathrm{St}_n(\beta)$. By hypothesis, one of the local representations, say $\pi_v$, is not supercuspidal. If $\pi_v$ is a *principal series representation* attached to quasi-characters $\mu, \nu$ of $F_v^*$, then let us set

$$\Pi_v := (\mu \otimes \pi'_v) \boxplus (\nu \otimes \pi'_v).$$

It remains to consider when $\pi_v$ is the special representation $\mathrm{St}_2(\nu)$ defined by a quasi-character $\nu$. If $\pi'_v$ *is supercuspidal*, then let us set

$$\Pi_v := \mathrm{St}_4(\pi'_v \otimes \nu).$$

If $\pi'_v$ *is a principal series representation* defined by quasi-characters $\mu', \nu'$, then we set

$$\Pi_v := (\mathrm{St}_2(\nu) \otimes \mu') \boxplus (\mathrm{St}_2(\nu) \otimes \nu').$$

82    DINAKAR RAMAKRISHNANFinally, if $\pi'_v$ is also a special representation $\mathrm{St}(\nu')$, for a quasi-character $\nu'$, then we set

$$\Pi_v := \mathrm{St}_3(\nu\nu') \boxplus \nu\nu'.$$

If we put $\Pi := \otimes_v \Pi_v$, then $\Pi$ is an irreducible, admissible, generic representation of $\mathrm{GL}(4, \mathbb{A}_F)$. We want to apply the converse theorem for $\mathrm{GL}(4)$ (cf. [CoPS], which is recalled in Theorem 3.1.3). Let $T$ be a finite set of finite places containing those $v$ where both $\pi_v$ and $\pi'_v$ are in the discrete series, and let $\eta$ be an arbitrary cuspidal automorphic representation of $\mathrm{GL}(m, \mathbb{A}_F)$ with $m = 1$ or $2$, such that $\eta_v$ is unramified at every $v$ in $T$. We are interested in the formally defined Rankin-Selberg $L$-function $L(s, \Pi \times \eta)$, and more precisely in the local integrals which show up in the proof of the converse theorem in [CoPS]. By [JPSS] we know that, at each finite place $v$, the $L$-factor $L(s, \Pi_v \times \eta_v)$ is the greatest common denominator of the corresponding local integrals $\Psi(s, W_v, W'_v, \Phi_v)$; when $v$ is archimedean, we know (cf. [JS4]) that the poles of $\Psi(s, W_v, W'_v, \Phi_v)$ are contained in the poles of the $L$-factor. Moreover, we have a local functional equation involving the correct $L$- and $\varepsilon$-factors. See also [MW2]. (When $v$ is archimedean, it is not known if the $L$-factor is realized in terms of the local integrals, unless one extends the function spaces [JS4], but this will not play any role for us.) Using this in conjunction with Theorems 3.3.11 and 3.4.1, we see that we can apply Theorem 3.1.3 and conclude the near-automorphy of $\Pi$ once we note the following:

LEMMA 3.5.2. *At every place $v$ of $F$, there exist the identities*

$$L(s, \Pi_v \times \eta_v) = L(s, \pi_v \times \pi'_v \times \eta_v)$$

*and*

$$\varepsilon(s, \Pi_v \times \eta_v) = \varepsilon(s, \pi_v \times \pi'_v \times \eta_v).$$

*Proof.* For $m = 1$, this is simply a rewording of Proposition 3.3.4. So assume that $m = 2$, and denote, for each $v$, the 2-dimensional representation of $W'_{F_v}$ associated to $\eta_v$ by $\tau_v$. Suppose $v$ is non-archimedean. If $\pi_v$ is a principal series representation defined by quasi-characters $\mu, \nu$, then $\sigma_v$ is simply $\mu \oplus \nu$, and $L(\pi_v \times \pi'_v \times \eta_v)$ (resp. $\varepsilon(\pi_v \times \pi'_v \times \eta_v)$) coincides with

$$L(s, \mu \otimes \sigma'_v \otimes \tau_v) L(s, \nu \otimes \sigma'_v \otimes \tau_v)$$

(resp. $\varepsilon(s, \mu \otimes \sigma'_v \otimes \tau_v) \varepsilon(s, \nu \otimes \sigma'_v \otimes \tau_v)$). Again the assertion follows from Proposition 3.3.4, thanks to the definition of $\Pi_v$. It works similarly when $\pi'_v$ is in the principal series. So we may assume that both $\pi_v$ and $\pi'_v$ are in the discrete series. Then by the choice of $(T, \eta)$, $\eta_v$ is an unramified principal series. If $\pi_v$ and $\pi'_v$ are both supercuspidal, the assertion follows once again by Proposition 3.3.4. It remains to consider when (i) $\pi_v$ is special and $\pi'_v$ is supercuspidal, and (ii) $\pi$ and $\pi'$ are both special. The lemma follows in both cases, due to our choice of $\Pi_v$, from Proposition 81 and Theorem 8.2 of [JPSS]. Finally let $v$ be archimedean. But here, one knows by [La1] that $L(s, \Pi_v \times \eta_v)$ (resp. $\varepsilon(s, \Pi_v \times \eta_v)$) equals $L(s, \sigma_v \otimes \sigma'_v \otimes \tau_v)$ (resp. $\varepsilon(s, \sigma_v \otimes \sigma'_v \otimes \tau_v)$) and we are done.

Now we know that $\Pi$ is nearly automorphic, i.e., that there exists an automorphic representation $\Pi_1$ of $\mathrm{GL}(4, \mathbb{A}_F)$ such that $\Pi_v$ is isomorphic to $\Pi_{1,v}$ for almost all $v$.

LEMMA 3.5.3. *There is an isobaric automorphic representation $\Pi_2$ of $\mathrm{GL}(4, \mathbb{A}_F)$ such that $\Pi_v \simeq \Pi_{2,v}$ for almost all $v$.*

*Proof.* Indeed, this is obvious if $\Pi_1$ is cuspidal. So assume not. Then $\Pi_1$ is a Jordan-Hölder component, thanks to Proposition 2 of [La6], of a representation induced from a cuspidal automorphic representation $\beta$ of (the Levi component) of a standard parabolic subgroup $P$ of $\mathrm{GL}(4)$. There is by definition an isobaric constituent $\Pi_2$ of this same induced representation, whose local components are at the unramified places the same as those of $\Pi_1$ (and those of $\Pi$). Moreover, by the reverse implication of Proposition 2 of [La6], $\Pi_2$ is also automorphic. $\square$

Now Theorem 3.5.1 follows by applying Proposition 3.2.1 by taking $\pi \boxtimes \pi'$ to be $\Pi_2$.

*Remark* 3.5.4. It is natural to wonder, since $\pi$ is generic, if it is exactly the same as $\Pi_2$. We will show later (in §4.3) that this is indeed the case. We will then have a uniqueness statement for $\pi \boxtimes \pi'$.

3.6. *A descent criterion.* The object of this section is to prove the following simple extension of the proposition in Section 4.2 of [BR]. (There is a mistake in [BR], but this does not affect Sections 1–4. In any case, the proof given below is completely self-contained.)

PROPOSITION 3.6.1. *Fix $n, p \in \mathbb{N}$ with $p$ prime. Let $F$ be a number field, $\{K_j \mid j \in \mathbb{N}\}$ a family of cyclic extensions of $F$ with $[K_j : F] = p$, and for each $j \in \mathbb{N}$, $\pi_j$ a cuspidal automorphic representation of $\mathrm{GL}(n, \mathbb{A}_{K_j})$. Suppose that, for all $j, r \in \mathbb{N}$, the base changes of $\pi_j, \pi_r$ to the compositum $K_j K_r$ satisfy*

(DC) $$(\pi_j)_{K_j K_r} \simeq (\pi_r)_{K_j K_r}.$$

*Then there exists a unique cuspidal automorphic representation $\pi$ of $\mathrm{GL}(n, \mathbb{A}_F)$ such that*
$$(\pi)_{K_j} \simeq \pi_j,$$
*for all but a finite number of $j$.*

*Proof.* First we need the following:

LEMMA 3.6.2. *Let $\pi$ be a cuspidal automorphic representation of* $\mathrm{GL}(n, \mathbb{A}_F)$. *Then there exist at most a finite number of idele class characters $\chi$ such that*

$$\pi \simeq \pi \otimes \chi.$$

*Proof of Lemma* 3.6.2. Suppose $\pi$ is isomorphic to its self-twist by an idele class character $\chi$. First we note that if $\omega$ is the central character of $\pi$, then the central character of $\pi \otimes \chi$ is simply $\omega \chi^n$. Thus $\chi$ must have order dividing $n$.

Next we claim that $\chi$ must be unramified at *any* finite place $v$ where $\pi$ is unramified. Indeed, for any such $v$, $L(s, \pi_v)$ is by definition, of the form $\prod_{j=1}^{n} L(s, \mu_j)$, for some unramified characters $\mu_j$ of $F_v^*$; this $L$-factor is not 1. On the other hand, if $\chi_v$ is ramified, $L(s, \pi_v \otimes \chi_v)$ is none other than $\prod_{j=1}^{n} L(s, \mu_j \chi)$, which equals 1. Hence the claim.

Consequently, if $S$ is the (finite) support of the conductor of $\pi$, then the number of possible self-twists is bounded by the number of idele class characters $\chi$ of $F$ such that the (finite) support of the conductor of $\chi$ is in $S$. Every such character cuts out, by class field theory, a cyclic extension $F(\chi)/F$ of degree $m \leq n$, which is ramified only at the places in $S$. But the number of such extensions is, for each $m \leq n$, known to be finite by a well-known variant of a classical theorem of Hermite. The assertion now follows. □

Now we continue with the proof of the proposition. For each $j$, let $\theta_j$ be a generator of $\mathrm{Gal}(K_j/F)$, and $\delta_j$ a character of $F$ cutting out $K_j$ (by class field theory). Note that, for each $i \geq 1$, the pull-back to $K_i$ of $\delta_j$ by the norm map $N_i$ from $K_i$ to $F$ cuts out the compositum $K_i K_j$.

We claim that

(3.6.3) $$\pi_j \circ \theta_j \simeq \pi_j \quad (\text{for all } j).$$

For all $j, r \geq 1$, let $\theta_{j,r}$ denote the automorphism of $K_j K_r$ such that

(i) $\theta_{j,r}|_{K_j} = \theta_j$, and

(ii) $\theta_{j,r}|_{K_r} = 1$ (where 1 denotes the identity automorphism).

It is easy to see that the base change of $\pi_j \circ \theta_j$ to $K_j K_r$ is simply $(\pi_j)_{K_j K_r} \circ \theta_{j,r}$. Applying (DC), we then have

$$(\pi_j \circ \theta_j)_{K_j K_r} \simeq (\pi_r)_{K_j K_r} \circ \theta_{j,r} \simeq (\pi_r)_{K_j K_r} \simeq (\pi_j)_{K_j K_r},$$

since $\theta_{j,r}$ is trivial on $K_r$. Since $K_j K_r$ is a cyclic extension of $K_j$ of prime degree, we have by Arthur-Clozel (Prop. 2.3.1)

$$\pi_j \circ \theta_j \simeq \pi_j \otimes (\delta_r \circ N_j)^{m_r},$$

for some $m_r \in \{0, 1, \ldots, p-1\}$. For every fixed $r \geq 1$, and for all $k \neq r$, we then have the self-twist identity

$$\pi_j \simeq \pi_j \otimes (\delta_r \circ N_j)^{m_r}(\delta_k \circ N_j)^{-m_k}.$$

Note that $\delta_r \circ N_j$ and $\delta_k \circ N_j$ must be distinct unless their ratio is a power of $\delta_j$. So the lemma above forces $m_r$ to be 0 for all but a finite number of $r$. Then, since (3.6.3) is independent of $r$, the claim follows.

As a result, by applying Proposition 2.3.1 again, we see that there exists, for each $j \geq 1$, a cuspidal automorphic representation $\pi(j)$ of $\mathrm{GL}(n, \mathbb{A}_F)$ such that

$$\pi_j \simeq (\pi(j))_{K_j}.$$

Such a $\pi(j)$ is of course unique only up to twisting by a power of $\delta_j$.

It is important to note that, for any $r \neq j$, we have the following compatibility for base change in (cyclic) stages:

$$(3.6.4) \qquad ((\pi(j))_{K_j})_{K_j K_r} \simeq ((\pi(j))_{K_r})_{K_j K_r}.$$

We see this as follows. Let $v$ be a finite place of $K_j K_r$ which is unramified for the data. Denote by $u$ (resp. $w$, resp. $w'$) the place of $F$ (resp. $K_j$, resp. $K_r$) below $v$. If $\sigma_u$ denotes the representation of $W'_{F_u}$ associated to $\pi(j)_u$, then

$$\mathrm{res}^{(K_j)_w}_{(K_j K_r)_v}(\mathrm{res}^{F_u}_{(K_j)_w}(\sigma_u)) \simeq \mathrm{res}^{(K_r)_{w'}}_{(K_j K_r)_v}(\mathrm{res}^{F_u}_{(K_r)'_w}(\sigma_u)).$$

Then (2.3.0) implies the local identity (for all such $v$)

$$((\pi(j)_u)_{(K_j)_w})_{(K_j K_r)_v} \simeq ((\pi(j)_u)_{(K_r)'_w})_{(K_j K_r)_v}.$$

The global isomorphism (3.6.4) follows by the strong multiplicity one theorem.

We can then rewrite (DC) as saying, for all $j, r \geq 1$,

$$((\pi(j))_{K_j})_{K_j K_r} \simeq ((\pi(r))_{K_j})_{K_j K_r}.$$

Applying Proposition 2.3.1, we get

$$(\pi(j))_{K_j} \simeq (\pi(r))_{K_j} \otimes (\delta_r \circ N_j)^{m(r,j)},$$

for some integer $m(r, j)$. We can replace $\pi(r)$ by $\pi(r) \otimes \delta_r^{-m(r,j)}$ and get

$$(\pi(j))_{K_j} \simeq (\pi(r))_{K_j}.$$

Then, by replacing $\pi(j)$ by a power of $\delta_j$, we can arrange for $\pi(j)$ and $\pi(r)$ to be isomorphic. In sum, we have produced, for *every pair* $(j, r)$, a *common descent*, say $\pi(j, r)$, of $\{\pi_j, \pi_r\}$.

Fix $a, b \in \{1, \ldots, p-1\}$, and consider the possible isomorphism

$$(3.6.5) \qquad \pi(j, r) \simeq \pi(j, r) \otimes \delta_j^a \delta_r^{-b}.$$

We claim that this cannot happen outside a finite set $S_{a,b}$ of pairs $(j, r)$. To see this fix a pair $(i, \ell)$ and consider the relationship of $\pi(i, \ell)$ to $\pi(j, r)$. Since $\pi(i,\ell)$ and $\pi(j,\ell)$ have the same base change to $K_\ell$, they must differ by twisting by a power of $\delta_\ell$. Similarly, $\pi(j,\ell)$ and $\pi(j,r)$ differ by a character twist as they have the same base change to $K_r$. Put together, this shows that $\pi(i,\ell)$ and $\pi(j,r)$ are twists of each other. Then (3.6.5) would imply that

$$\pi(i,\ell) \simeq \pi(i,\ell) \otimes \delta_j^a \delta_r^{-b}.$$

The claim now follows since (by the lemma above) $\pi(i,\ell)$ admits only a finite number of self-twists and since the characters $\delta_j^a \delta_r^{-b}$ are all distinct for distinct pairs $(j,r)$ (as $a,b$ are fixed).

Now choose a pair $(j,r)$ *not belonging* to $S_{a,b}$ for *any* $(a,b) \in (\mathbb{Z}/p)^* \times (\mathbb{Z}/p)^*$, and set

$$\pi = \pi(j,r).$$

We assert that, for all but a finite number of indices $m$,

$$\pi_{K_m} \simeq \pi_m.$$

It suffices to show that, for any large enough $m$, $\pi = \pi(j,r)$ is isomorphic to either $\pi(j,m)$ or $\pi(m,r)$. Suppose neither is satisfied. Then there exist $a, b \in \{1, \ldots, p-1\}$ such that

$$\pi(j,m) \simeq \pi(j,r) \otimes \delta_j^a \quad \text{and} \quad \pi(m,r) \simeq \pi(j,r) \otimes \delta_r^b.$$

We also have $\pi(j,m) \simeq \pi(m,r) \otimes \delta_m^c$, for some $c \in \{0, 1, \ldots, p-1\}$. Putting these together, we get the self-twisting identity

$$\pi(j,r) \simeq \pi(j,r) \otimes \delta_j^a \delta_r^{-b} \delta_m^{-c}.$$

By our choice of $(j,r)$, $c$ cannot be 0. For each nonzero $c$, the set of indices $m$ for which such an identity can hold is finite (again by the lemma). Hence we get a contradiction for large enough $m$, which implies that $a$ or $b$ should be 0. The proposition is now proved. □

3.7. *Modularity in the general case.* In this section we will complete the proof of Theorem M without any hypothesis on $(\pi, \pi')$. We begin by noting the following (well-known):

LEMMA 3.7.1. *Let $v$ be a finite place where $\pi_v$ is supercuspidal. Then there is a finite solvable, normal extension $E$ of $F_v$ such that the base change $\pi_{v,E}$ is in the principal series.*

*Proof.* Let $\sigma_v$ denote as before the 2-dimensional representation of $W'_{F_v}$ associated to $\pi_v$. By [Ku], $\sigma_v$ is an irreducible of $W_{F_v}$. The structure of the Weil group implies that $\sigma_v$ must in fact be of the form $\tau_v \otimes \nu$, where $\tau_v$ is an irreducible $\mathbb{C}$-representation of $W_{F_v}$ *factoring through a finite quotient*, and $\nu$ a quasi-character. So there is a finite (necessarily solvable) normal extension $E/F_v$ such that the restriction of $\tau_v$ to the Weil group of $E(v)$ is trivial. (Take $E$ to be the fixed field of $\overline{F}_v$ under the kernel of $\tau_v$.) Then by the functoriality of base change, $(\pi_v)_{E(v)}$ must be in the principal series, in fact of the form $(1 \boxplus 1) \otimes (\nu \circ N_{E(v)/F})$. □

It should be noted that by using Krasner's lemma, we can find a solvable Galois extension $k'/k$ of number fields with local extension $E/F_v$ such that $\mathrm{Gal}(k'/k)$ is isomorphic to $\mathrm{Gal}(E/F_v)$ (see [PR, §4, Lemma 3], for example). But we *cannot* hope to be able to take $k$ to be $F$! But we *can* find (see below) a finite chain of cyclic extensions of $F$ with good properties. However, the field on top need not be normal over $F$.

*Proof of Theorem* M. We may, by Proposition 3.1.1, assume that we are not in one of the special situations (I), (II), (III) treated there.

Given any finite solvable group $G$, define its *length* $\ell(G)$ to be the length $\ell$ of any chain $\{1\} = G_0 \subset G_1 \subset \cdots \subset G_\ell = G$, with each $G_i$ normal in $G_{i+1}$ of prime index. Let $S = S(\pi)$ denote the (possibly empty) finite set of finite places where either $\pi_v$ is not supercuspidal. At each $v \in S$, let $\ell(\pi_v)$ denote the length of a minimal Galois extension $L(v)/F_v$ such that the base change $(\pi_v)_{L(v)}$ is in the principal series. Let $\ell(\pi)$ denote the maximum of $\{\ell(\pi_v) \mid v \in S\}$, and let $S'$ denote the subset of $S$ where this maximum is attained. Further, for each $v$ in $S'$, let $p(v)$ denote the maximum over all $L(v)$, of the degree, required to be a prime or 1, of the largest cyclic extension $K(v)$ of $F_v$ contained in $L(v)$. Let $p = p(\pi)$ be the maximum of $p(v)$ over all $v$ in $S'$, and let $S''$ denote the subset of $S'$ where $p(v) = p$ (and $\ell(\pi_v) = \ell(\pi)$). Note that $p$ is a *prime* unless $\pi$ is *good* over $F$, i.e., has no supercuspidal components, in which case $p = 1$.

Let us set

(3.7.2) $$r(\pi) = (\ell(\pi), p(\pi)).$$

We will order these pairs as follows: $(\ell, p) < (\ell', p)$ if either $\ell < \ell'$ or if $\ell = \ell'$ *and* $p < p'$.

Suppose $r = r(\pi)$ is $(0, 1)$. Then, for any quadratic extension $K'$ which is totally complex, the assertion holds by Theorem 3.5.1. Also, given any finite place $u$, we can find a totally complex, quadratic extension $K'(u)$ such that $u$ splits in $K'(u)$. Applying the descent criterion (Prop. 3.6.1), we can find a common descent having the requisite properties.

Now let $r > (0, 1)$, and assume by induction that the theorem is proved (over all number fields $K$) for pairs $(\pi_1, \pi_2)$ of cuspidal automorphic representations of $\mathrm{GL}(2, \mathbb{A}_K)$ with $r(\pi_1, \pi_2) < r$.

Fix, at every place $v$ in $S''$, a character $\chi_v$ of $F_v^*$ cutting out a $K(v)$ (as above) of degree $p$.

Enumerate the set of finite places where $\pi$ is unramified, and list them as $\{v_1, v_2, \ldots, v_j, \ldots\}$.

Fix an index $j \geq 1$ (for the moment), and let $S(j)$ be the union of $v_j$ with $S''$. Let $\chi_{v_0}$ denote the trivial character of $F_{v_0}^*$.

Now by the Grunewald-Wang theorem (see [AT, Chap. 10, Thm. 5]), we can find a global character $\chi(j)$ of $C_F$ of order $\ell$ whose local restrictions are given by $\chi_v$ at every $v$ in $S(j)$. (Note that the set $S_0$ which occurs in *loc. cit.* is empty in our case as each local degree is $p$ or $1$ and cannot be divisible by 4, thus allowing us to find $\chi(j)$ of order $p$, not just $2p$.) Let $K_j$ be the $p$-extension of $F$ cut out by $\chi(j)$.

By construction we have, for every $j \geq 1$,

$$(3.7.3) \qquad r(\pi) < r.$$

Thus, by induction, the theorem holds for $\pi_{K_j}$ for each $j$. Note that if the automorphic representation $\pi_{K_j}$ is not cuspidal for some $j$, it must be dihedral (with $p = 2$), and by the remark earlier, we may assume that we are not in this case. Put

$$(3.7.4) \qquad \Pi_j = \pi_{K_j} \boxtimes \pi'_{K_j}.$$

The nondihedrality of $\pi$ also gives the following:

LEMMA 3.7.5. *There exist at most a finite number of indices $j$ such that $\Pi_j$ is not cuspidal.*

*Proof.* It suffices to prove that, outside a finite set of indices, we have

1. $\pi_{K_j}$ is not dihedral;

2. $\pi_{K_j}$ is not the twist of $\pi'_{K_j}$ by an idele class character of $K_j$.

Suppose $\pi_{K_j}$ is dihedral. Then, by Gelbart-Jacquet ([GJ]), the symmetric square of $\pi_{K_j}$ is Eisensteinian, while by assumption, $\text{sym}^2(\pi)$ is not. This forces $p$ to be 3 and $\text{sym}^2(\pi)$ to be automorphically induced by a character $\mu$ of $K_j$. In other words, the symmetric square of $\pi$ admits a self-twist relative to the character of $F$ corresponding to $K_j$. By Lemma 3.6.2, this can only happen for a finite number of $j$.

We may then assume that $\pi_{K_j}$ is not dihedral and consider the second assertion. Suppose $\pi_{K_j}$ is the twist of $\pi'_{K_j}$ by an idele class character. Then, for a character $\mu$ of $C_{K_j}$, $L^{T'}(s, \Pi_j \otimes \mu)$ must have a pole at $s = 1$, where $T'$ is any finite set of places of $K_j$. By the inductivity of $L$-functions, the same holds for the ($L$-function of the) automorphic induction $I_{K_j}^F(\Pi_j)$. Let us choose $T'$ to be the inverse image of the set $T$ of places of $F$ where $\pi$ or $\pi'$ or $K_j$ is ramified. Then it is easy to see that we have the identity:

$$L^T(s, I_{K_j}^F(\Pi_j \otimes \mu)) = L^T(s, \pi \times (\pi' \boxtimes I_{K_j}^F(\mu))).$$

Since $T$ contains all the ramified places, the right-hand side $L$-function is an incomplete form of *any* of the triple product $L$-functions attached to $(\pi, \pi', I_{K_j}^F(\mu))$. By [Ik1], it cannot have a pole unless $\pi$ is dihedral, and we are done. □

Consequently, after shrinking the index set by removing an appropriate finite subset and renumbering it, we may assume that

$$\Pi_j \text{ is cuspidal for every } j.$$

Next we fix a pair $(j, r)$ of indices and consider the descent criterion (DC) of Proposition 3.6.1. Let $w$ be a finite place where $((\Pi_j)_{K_jK_r})_w$, $((\Pi_r)_{K_jK_r})_w$ and $K_iK_j$ are all unramified. Then, by construction, both these local representations correspond to the restriction (to the Weil group of $(K_jK_r)_w$) of $\sigma_v \otimes \sigma'_v$, where $v$ signifies the place of $F$ below $w$. (Recall that $\sigma_v, \sigma'_v$ are associated to $\pi_v, \pi'_v$ respectively.) This leads to the identity

$$L^Y(s, (\Pi_j)_{K_jK_r}) = L^Y(s, (\Pi_r)_{K_jK_r}),$$

for a large enough finite set $X$ of places, which gives (DC) by the strong multiplicity one theorem.

Applying Proposition 3.6.1, we then get a unique cuspidal descent $\Pi$ on $\mathrm{GL}(4)/F$ such that, for all but a finite number of indices,

$$\Pi_{K_j} \simeq \Pi_j.$$

Finally, by construction, each (unramified) finite place $v_j$ splits completely in $K_j$; let $w_j$ be a divisor of $v_j$ in $K_j$. This implies (by the definition of base change) that, for almost all $j$,

$$L(s, \Pi_{v_j}) = L(s, (\Pi_j)_{w_j}) = L(s, (\pi_{K_j})_w \times (\pi'_{K_j})_w) = L(s, \pi_v \times \pi'_v),$$

and similarly for the $\varepsilon$-factors. Thus $\Pi$ is a weak (tensor product) lifting of $(\pi, \pi')$. By Proposition 3.2.1, it is also the strong lifting. □

## 4. Applications

4.1. *A multiplicity one theorem for* $\mathrm{SL}(2)$. Denote by

$$L_0^2(\mathrm{SL}(2, F) \backslash \mathrm{SL}(2, \mathbb{A}_F))$$

the subspace of cusp forms in $L^2(\mathrm{SL}(2, F) \backslash \mathrm{SL}(2, \mathbb{A}_F))$. It is acted upon by $\mathrm{SL}(2, \mathbb{A}_F)$ by right translation, and this action is unitary relative to the natural scalar product on $L^2(\mathrm{SL}(2, F) \backslash \mathrm{SL}(2, \mathbb{A}_F))$. The main result of this section is:

THEOREM 4.1.1.  *The representation of* $\mathrm{SL}(2, \mathbb{A}_F)$ *on*

$$L_0^2(\mathrm{SL}(2, F) \backslash \mathrm{SL}(2, \mathbb{A}_F))$$

*has multiplicity one.*

This was conjectured by Labesse and Langlands (see [LL] and [Lab1]). The theorem says that the commutant of this representation is abelian, and equivalently, since the space of cusp forms is completely reducible, every irreducible unitary representation of $\mathrm{SL}(2, \mathbb{A}_F)$ occurring in (the decomposition of) this space appears with multiplicity one.

The stable trace formula for $\mathrm{SL}(2)$ was analyzed deeply, and in detail, in [LL]. We refer to this paper for all the background material. The results there show in particular that multiplicity one for $\mathrm{SL}(2)$ is a consequence of the following:

THEOREM 4.1.2. *Let $\pi, \pi'$ be unitary, cuspidal automorphic representations of $\mathrm{GL}(2, \mathbb{A}_F)$. Suppose we have, for almost all $v$,*

(LL(v)) $$\mathrm{Ad}(\pi_v) \simeq \mathrm{Ad}(\pi'_v).$$

*Then there exists an idele class character $\chi$, which is unique if $\pi$ is not automorphically induced by a character of a quadratic extension, such that*

$$\pi' \simeq \pi \otimes \chi.$$

*If $\pi$ and $\pi'$ have the same central character, then $\chi$ is quadratic.*

Here $\mathrm{Ad}(\pi)$ denotes the automorphic representation $\mathrm{sym}^2(\pi) \otimes \omega^{-1}$ of $\mathrm{GL}(3, \mathbb{A}_F)$, where $\omega$ is the central character of $\pi$.

As noted in the introduction, in the special case when $\pi$ and $\pi'$ are defined by holomorphic eigenforms, with $F$ necessarily totally real, one can establish this result [Ra2] by making use of the associated $\ell$-adic representations. But this does not work in general, and the main point here is to give a unified proof using the modularity of $\pi \boxtimes \pi'$.

Before beginning the proof of this theorem, we would like to indicate the following concrete result:

COROLLARY 4.1.3. *Let $\pi, \pi'$ be unitary, cuspidal automorphic representations of $\mathrm{GL}(2, \mathbb{A}_F)$ of trivial character. Suppose we have at every unramified finite place $v$,*

(Sq(v)) $$a_v(\pi)^2 = a_v(\pi')^2,$$

*where $a_v(\pi) = \mathrm{tr}(A_v(\pi))$ is the $v^{\mathrm{th}}$ Hecke eigenvalue of $\pi$. Then there is a unique quadratic idele class character $\chi$ such that*

$$a_v(\pi') = \chi_v(\varpi_v) a_v(\pi'),$$

*for almost all $v$, where $\varpi_v$ denotes the uniformizer at $v$. If in addition, the conductors $N, N'$ of $\pi, \pi'$ are square-free, then $\chi = 1$.*

*Proof that Theorem* 4.1.2 *implies Corollary* 4.1.3. Let $\pi, \pi'$ be as in the corollary, and let $v$ be a finite place where both $\pi$ and $\pi'$ are unramified. Since by hypothesis the central characters are trivial, we may write $A_v(\pi) = \mathrm{diag}(\alpha_v, \alpha_v^{-1})$ and $A_v(\pi') = \mathrm{diag}(\alpha'_v, \alpha'^{-1}_v)$. Then the Langlands class of $\mathrm{Ad}(\pi_v)$ (resp. $\mathrm{Ad}(\pi'_v)$) is given by $\{\alpha_v^2, 1, \alpha_v^{-2}\}$ (resp. $\{\alpha'^2_v, 1, \alpha'^{-2}_v\}$). Consequently, we have (in $\Re(s) > 1$)

$$\log(L(s, \mathrm{Ad}(\pi_v))) = 1 + \sum_{m \geq 1} \frac{\alpha_v^{2m} + 1 + \alpha_v^{-2m}}{m(Nv)^{ms}}$$

and

$$\log(L(s, \mathrm{Ad}(\pi'_v))) = 1 + \sum_{m \geq 1} \frac{\alpha'^{2m}_v + 1 + \alpha'^{-2m}_v}{m(Nv)^{ms}}.$$

We claim that, for every $m \geq 1$,

$$\alpha_v^{2m} + \alpha_v^{-2m} = \alpha'^{2m}_v + \alpha'^{-2m}_v.$$

By induction on $m$, it suffices (by the binomial formula) to verify

$$(\alpha_v^2 + \alpha_v^{-2})^m = (\alpha'^2_v + \alpha'^{-2}_v)^m,$$

which follows immediately from the hypothesis (Sq(v)). This proves the claim.

Since (LL(v)) then holds for almost all $v$, the corollary follows by applying Theorem 4.1.2.

*Proof of Theorem* 4.1.2. First note that, since (LL(v)) holds for almost all $v$, the automorphic representations $\mathrm{Ad}(\pi)$ and $\mathrm{Ad}(\pi')$ are isomorphic by the strong multiplicity one theorem.

Suppose $\pi$ is dihedral, i.e., of the form $I_K^F(\mu)$, for a character $\mu$ of $C_K$, for some quadratic extension $K$ of $F$. Its central character then identifies with $\mu_0 \delta$, where $\mu_0$ is the restriction of $\mu$ to $C_F$ and $\delta$ is the quadratic character of $C_F$ associated to $K/F$ by class field theory. Denote by $\theta$ the nontrivial automorphism of $K/F$. We claim that

(ad) $\qquad \mathrm{Ad}(\pi) = \delta \boxplus I_K^F(\mu/(\mu \circ \theta)).$

Since $(\mu_0)_K = \mu(\mu \circ \theta)$, this is equivalent to the identification of $\mathrm{sym}^2(\pi)$ with $\mu_0 \boxplus I_K^F(\mu^2)$, which is easy to verify at the unramified places. So the claim follows by the strong multiplicity one theorem ([JS2]).

A particular consequence of the claim is that $L(s, \mathrm{Ad}(\pi) \otimes \delta)$ has a pole at $s = 1$. This forces $\pi'$ also to be dihedral, for otherwise, $L(s, \mathrm{Ad}(\pi') \otimes \delta)$ will be entire ([GJ]). Then $\delta$ must in fact occur in the isobaric decomposition of $\mathrm{Ad}(\pi'_K)$, which in turn is a summand of $\pi' \boxtimes \pi'^\vee$. Then $L(s, \pi' \boxtimes \pi'^\vee \otimes \delta)$ also has a pole at $s = 1$. This implies that

$$\pi' \simeq \pi' \otimes \delta.$$

In other words, $\pi'$ is of the form $I_K^F(\mu')$, for a character $\mu'$ of $C_K$.

Now, base changing $\mathrm{Ad}(\pi) = \mathrm{Ad}(\pi')$ to $K$, and using (ad), we get the identity

$$\mu/(\mu \circ \theta) \boxplus (\mu \circ \theta)/\mu \simeq \mu'/(\mu' \circ \theta) \boxplus (\mu' \circ \theta)/\mu'.$$

Replacing $\mu'$ by $\mu' \circ \theta$ if necessary, we can then deduce that

$$\mu/\mu' = (\mu/\mu') \circ \theta.$$

In other words, there exists a character $\chi$ of $C_F$ such that $\mu' = \mu \chi_K$. Then

$$\pi' = I_K^F(\mu') \simeq I_K^F(\mu) \otimes \chi \simeq \pi \otimes \chi,$$

as desired. To verify the middle isomorphism, one first checks it locally almost everywhere and then applies the strong multiplicity one theorem.

So we may assume henceforth, in the main case of interest, that neither $\pi$ nor $\pi'$ is dihedral. Put

$$\Pi = \pi \boxtimes \pi'.$$

We know that this is an isobaric automorphic representation of $\mathrm{GL}(4, \mathbb{A}_F)$. The main point here is to prove the following:

LEMMA 4.1.4.  *$\Pi$ is not cuspidal.*

*Proof of Lemma* 4.1.4.  Let $S$ be a finite set of places containing the ramified places (for $\pi$ and $\pi'$) and the archimedean ones.

Suppose $\Pi$ is cuspidal. Then, by [JS2], we know that the Rankin-Selberg $L$-function $L^S(s, \Pi \times \Pi^\vee)$ must have a *simple pole* at $s = 1$.

On the other hand, at any finite place $v$, if we denote (as usual) by $\sigma_v$ (resp. $\sigma_v'$) the representation of $W_{F_v}'$, we get, as $\sigma_v \otimes \sigma_v^\vee \simeq \mathrm{Ad}(\sigma_v) \oplus 1$, the following identity:

$$(\sigma_v \otimes \sigma_v') \otimes (\sigma_v \otimes \sigma_v')^\vee \simeq 1 \oplus \mathrm{Ad}(\sigma_v) \oplus \mathrm{Ad}(\sigma_v') \oplus \mathrm{Ad}(\sigma_v) \otimes \mathrm{Ad}(\sigma_v').$$

At the unramified places $v$, we know that $\sigma_v \otimes \sigma_v'$ corresponds to $\Pi_v$. So this translates to the $L$-function identity

(Id)  $L^S(s, \Pi \times \Pi^\vee) = \zeta_F^S(s) L^S(s, \mathrm{Ad}(\pi)) L^S(s, \mathrm{Ad}(\pi')) L^S(s, \mathrm{Ad}(\pi) \times \mathrm{Ad}(\pi')).$

But we know that $\mathrm{Ad}(\pi)$ is isomorphic to $\mathrm{Ad}(\pi')$. Moreover, since $\mathrm{Ad}(\pi) = \mathrm{sym}^2(\pi) \otimes \omega^{-1}$,

$$\mathrm{Ad}(\pi)^\vee \simeq \mathrm{sym}^2(\pi^\vee) \otimes \omega \simeq \mathrm{Ad}(\pi),$$

as $\mathrm{sym}^2(\pi^\vee)$ is $\mathrm{sym}^2(\pi) \otimes \omega^{-2}$.

Consequently, $L^S(s, \mathrm{Ad}(\pi) \times \mathrm{Ad}(\pi'))$ has a pole at $s = 1$. Also, since $\pi$ and $\pi'$ are nondihedral, $L^S(s, \mathrm{Ad}(\pi))$ and $L^S(s, \mathrm{Ad}(\pi'))$ are entire (cf. [GJ]). Then the identity (Id) implies that $L^S(s, \Pi \times \Pi^\vee)$ has at least a double pole at $s = 1$, leading to the desired contradiction.  □

Thus $\pi, \pi'$ are nondihedral cusp forms on $\mathrm{GL}(2)/F$ with $\pi \boxtimes \pi'$ not cuspidal. Then by the cuspidality criterion of Theorem M (see §3), which was proved in Section 3.2, we see that $\pi'$ must be isomorphic to $\pi \otimes \chi$ for an idele class character $\chi$, proving Theorems 4.1.2 and 4.1.1. □

4.2. *Some new functional equations.* Let $\pi_1, \pi_2, \pi_3, \pi_4$ be cuspidal automorphic representations of $\mathrm{GL}(2, \mathbb{A}_F)$, and let $S$ be a finite set of places containing the archimedean and ramified places. Let us set

$$L^S(s, \pi_1 \times \pi_2 \times \pi_3 \times \pi_4) = \prod_{v \notin S} \det(I - (Nv)^{-s} A_v(\pi_1) \otimes A_v(\pi_2) \otimes A_v(\pi_3) \otimes A_v(\pi_4))^{-1},$$

$$L^S(s, \mathrm{sym}^2(\pi_1) \times \pi_2 \times \pi_3) = \prod_{v \notin S} \det(I - (Nv)^{-s} \mathrm{sym}^2(A_v(\pi_1)) \otimes A_v(\pi_2) \otimes A_v(\pi_3))^{-1},$$

and

$$L^S(s, \mathrm{sym}^3(\pi_1) \times \pi_2) = \prod_{v \notin S} \det(I - (Nv)^{-s} \mathrm{sym}^3(A_v(\pi_1)) \otimes A_v(\pi_2))^{-1}.$$

These functions of $s$ converge absolutely in some right half-plane.

THEOREM 4.2.1. *Each of the three (incomplete) Euler products above admit meromorphic continuations to the whole $s$-plane. In addition, each of them can be completed to a convergent Euler product over all places in $\Re(s) > 1$, admitting meromorphic continuation and a functional equation of the usual form (as conjectured by Langlands). They are nonzero on the line $s = 1$ with no pole except possibly at the point $s = 1$. The first two $L$-functions are in fact analytic everywhere except for possible poles at $s = 0, 1$.*

*Remark.* The meromorphic continuation and functional equation of the third of the $L$-functions above is already known by the work of Shahidi using the theory of Eisenstein series ([Sh2]). But his approach does not seem to be able to give invertibility in $\Re(s) \geq 1$ except for a possible pole at $s = 1$.

*Proof.* Since by Theorem M, there exist $\pi_1 \boxtimes \pi_2$ and $\pi_3 \boxtimes \pi_4$ in $\mathcal{A}(4, F)$, we may set

$$L(s, \pi_1 \times \pi_2 \times \pi_3 \times \pi_4) = L(s, (\pi_1 \boxtimes \pi_2) \times (\pi_3 \boxtimes \pi_4)).$$

By the definition of $\boxtimes$, the unramified local factors of the $L$-function on the right agree with those of $L^S(s, \pi_1 \times \pi_2 \times \pi_3 \times \pi_4)$. This way we get a complete Euler product in $\Re(s) > 1$, and the remaining assertions follow immediately for the first $L$-function from the corresponding ones, due to Jacquet, Piatetski-Shapiro and Shalika ([JPSS]), Mœglin-Waldspurger ([MW2]), Shahidi [Sh6], for the Rankin-Selberg $L$-functions on $\mathrm{GL}(4) \times \mathrm{GL}(4)$.

We complete the second and third $L$-functions by respectively setting

$$L(s, \mathrm{sym}^2(\pi_1) \times \pi_2 \times \pi_3) \;=\; L(s, \mathrm{sym}^2(\pi_1) \times (\pi_2 \boxtimes \pi_3))$$

and

$$L(s, \mathrm{sym}^3(\pi_1) \times \pi_2) \;=\; \frac{L(s, \mathrm{sym}^2(\pi_1) \times (\pi_1 \boxtimes \pi_2))}{L(s, (\pi_1 \boxtimes \pi_2) \otimes \omega_1)},$$

where $\omega_1$ is the central character of $\pi_1$.

The verification of the fact that the unramified factors of these $L$-functions are the right ones is an easy exercise, which is left to the reader. Since the Rankin-Selberg $L$-functions on $\mathrm{GL}(n) \times \mathrm{GL}(m)$ are nonzero in $\{\Re(s) \geq 1\}$ and are analytic everywhere except for possible poles at $s = 0, 1$, all the assertions of Theorem 4.2.1 follow except for knowing that the third $L$-function is nonzero at $s = 1$. So we will be done if we prove the following:

LEMMA 4.2.2.

$$-\mathrm{ord}_{s=1} L(s, \pi_1 \boxtimes \pi_2 \otimes \omega_1) \;\leq\; -\mathrm{ord}_{s=1} L(s, \mathrm{sym}^2(\pi_1) \times (\pi_1 \boxtimes \pi_2)).$$

Suppose $L(s, \pi_1 \boxtimes \pi_2 \otimes \omega_1)$ has a pole at $s = 1$. Since $\pi_1, \pi_2$ are cuspidal, the pole has to be a simple pole, and besides, this can happen if and only if $\pi_2$ is isomorphic to the contragredient of $\pi_1 \otimes \omega_1$. It then follows, since $\pi_1 \boxtimes \pi_1^\vee = \mathrm{Ad}(\pi_1) \boxplus 1$, that

(Fac)    $L(s, \mathrm{sym}^2(\pi_1) \times (\pi_1 \boxtimes \pi_2)) \;=\; L(s, \mathrm{Ad}(\pi_1) \times \mathrm{Ad}(\pi_1)) L(s, \mathrm{Ad}(\pi_1)).$

We saw earlier that $\mathrm{Ad}(\pi_1)$ is self-dual. Hence $L(s, \mathrm{Ad}(\pi_1) \times \mathrm{Ad}(\pi_1))$ has a pole at $s = 1$. We claim that $L(s, \mathrm{Ad}(\pi_1))$ is holomorphic at $s = 1$, which is clear if $\pi$ is nondihedral, as $\mathrm{Ad}(\pi_1)$ is then cuspidal by [GJ]. So suppose $\pi$ is of the form $I_K^F(\mu)$ for a character of (the idele classes of) a quadratic extension $K$ of $F$, with nontrivial automorphism $\theta$ over $F$. Then, by the identity (ad) appearing in the proof of Theorem 4.1.2, $L(s, \mathrm{Ad}(\pi_1))$ can have a pole at $s = 1$ if and only if $\mu/(\mu \circ \theta)$ is the pull-back by norm to $K$ of either the trivial character of $C_F$ or the quadratic character $\delta$ (attached to $K$). In either case, we will have $\mu = \mu \circ \theta$, contradicting the cuspidality of $\pi_1$.

Consequently, the right- (and hence the left-) hand side of (Fac) has a pole at $s = 1$. This proves the lemma and completes the proof of the theorem. □

*Remark* 4.2.3. It is now clear how we can also define, in an analogous way, the epsilon factors, both locally and globally, for the $L$-functions under consideration.

4.3. *Root numbers and representations of orthogonal type.* Let $k$ be a non-archimedean local field of characteristic zero. Our first aim is to prove the following, which is not clear even if one knew the local Langlands conjecture for $GL(4)$.

PROPOSITION 4.3.1. *Let $n \leq 4$, and let $\{\eta_j | j \leq n\}$ be a set of supercuspidals of $GL(2, k)$, with corresponding 2-dimensional irreducibles $\{\tau_j | j \leq n\}$. Then we have*

$$L(s, \eta_1 \times \ldots \times \eta_n) = L(s, \tau_1 \otimes \ldots \otimes \tau_n)$$

*and*

$$\varepsilon(s, \eta_1 \times \ldots \times \eta_n) = \varepsilon(s, \tau_1 \otimes \ldots \otimes \tau_n).$$

*Proof of Proposition* 4.3.1. For $n = 1$, this is contained in [Ku]. It is also well-known for $n = 2$; see [PR, Prop. 4.2], for example. So assume $n \geq 3$. If either of the $\eta_j$ is the local automorphic induction of a (quasi)character $\mu$ of $E^*$, for some quadratic extension $E/k$, then the proposition is the same as Proposition 4.2 of *loc. cit.* Otherwise, we can modify the argument there and appeal to the global Rankin-Selberg product to get what we want. We give the details for completeness and we start by appealing to [PR, Lemma 3, §4], and find a number field $F$ with $k = F_u$ for some place $u$, and irreducible 2-dimensional representations $\sigma_j$ of $\mathrm{Gal}(\overline{\mathbb{Q}}/F)$ with solvable image such that their respective restrictions to the decomposition group at $u$ identify with $\tau_j$. Solvability allows us to apply the theorem of Langlands ([La5]) and Tunnell ([Tu]), and find corresponding cuspidal automorphic representations $\pi_j$ with $\pi_{j,u} \simeq \eta_j$. Then we choose a finite order character $\mu$ of $C_F$ which is highly ramified at all the ramified places except $u$ where it is 1. Comparing the functional equations of $L(s, \pi_1 \times \ldots \times \pi_n)$ and $L(s, \sigma_1 \otimes \ldots \otimes \sigma_n)$, we proceed as in Proposition 3.2.1 and deduce the assertions above. □

*Remark* 4.3.2. In Chapter 3 we associated, to each pair $(\pi, \pi')$ of cuspidal automorphic representations of $GL(2, \mathbb{A}_F)$, an irreducible admissible representation $\Pi$ of $GL(4, \mathbb{A}_F)$ and an isobaric automorphic representation $\pi \boxtimes \pi'$ such that $\Pi_v \simeq (\pi \boxtimes \pi')_v$ for almost all $v$. Using Proposition 4.3.3 we can conclude that such an identity holds at *every* $v$.

To clarify this remark, we will digress a bit about the local Langlands conjecture for $GL(4)$, which predicts a canonical bijection, preserving local factors, between the set $\mathcal{A}^0(4, k)$ of irreducible, supercuspidal representations $\eta$ of $GL(4, k)$ and the set $\mathcal{R}^0(4, k)$ of irreducible, continuous, 4-dimensional $\mathbb{C}$-representations $\tau$ of $W_k$. There is a numerical bijection between the sets ([He4], compatible with twisting by unramified characters. A theorem of M. Harris ([Ha]) gives a (geometrically defined) bijection $\tau \to \eta$, which gives moreover an equality of epsilon factors of pairs

$$\varepsilon(s, \tau \otimes \tau') = \varepsilon(s, \eta \times \eta'),$$

for all irreducibles $\tau'$, with corresponding $\eta'$, of dimension $m$, for all $m \leq 4$, *provided* the residual characteristic is prime to $2m$. On the other hand, a recent result of J. Chen ([Ch]) asserts that any $\eta$ in $\mathcal{A}^0(4,k)$ is uniquely determined by the collection of epsilon factors $\varepsilon(s, \eta \times \eta')$, as $\eta'$ ranges over $\mathcal{A}^0(m,k)$ for $m = 1$ and $m = 2$. Consequently, as noted in [PR, §4], the local Langlands conjecture follows for GL(4) in the odd residual characteristic case by combining the results of [Ha] and [Ch]. In the case of even residual characteristic, shown by [He2], one knows the bijection for representations of "cyclic type", i.e., those which admit a quartic self-twist. Moreover, Proposition 4.1 of [PR] shows how to extend this and treat also the case of representations admitting just a quadratic self-twist.

Now let $\mathrm{GO}(4, \mathbb{C})$ denote the subgroup of $\mathrm{GL}(4, \mathbb{C})$ consisting of orthogonal similitudes with similitude factor $\mu$. Let $\mathrm{GSO}(4, \mathbb{C})$ denote the kernel of the homomorphism $\mathrm{GO}(4, \mathbb{C}) \to \{\pm 1\}$, $g \to \mu(g)^{-2}\det(g)$. In this section we will be concerned with the subset $\mathcal{R}^1(4,k)$ of $\tau$ in $\mathcal{R}^0(4,k)$ whose image lie in $\mathrm{GSO}(4, \mathbb{C})$. It is well known that there is a short exact sequence

$$1 \to \mathbb{C}^* \to \mathrm{GL}(2, \mathbb{C}) \times \mathrm{GL}(2, \mathbb{C}) \to \mathrm{GSO}(4, \mathbb{C}) \to 1,$$

inducing a surjection

$$\mathrm{Hom}(W_k, \mathrm{GL}(2, \mathbb{C}) \times \mathrm{GL}(2, \mathbb{C})) \to \mathrm{Hom}(W_k, \mathrm{GSO}(4, \mathbb{C})).$$

See, for example, the proof of Lemma 4 in [PR, §5]. Consequently, every $\tau$ in $\mathcal{R}^1(4,k)$ is a tensor product $\tau_1 \otimes \tau_2$, with each $\tau_j$ an irreducible of dimension 2. The pair $(\tau_1, \tau_2)$ is not unique, as it can be replaced by $(\tau_1 \otimes \nu, \tau_2 \otimes \nu^{-1})$, but this is the only ambiguity.

Motivated by this, we define $\mathcal{A}^1(4,k)$ to be the irreducible supercuspidal representations $\eta$ of $\mathrm{GL}(4,k)$ which are of the form $\eta_1 \boxtimes \eta_2$ (see §3.7 for a definition). In what follows, we will make use of the local Langlands correspondence for GL(2) and GL(3), established respectively by Kutzko ([Ku]) and Henniart ([He1]).

PROPOSITION 4.3.3. *There is a unique bijection $\tau \to \eta$ between $\mathcal{R}^1(4,k)$ and $\mathcal{A}^1(4,k)$, compatible with twisting by characters and taking duals, such that the following identity of epsilon factors holds for all $\tau' \in \mathcal{R}^0(m,k)$, with corresponding $\eta' \in \mathcal{A}^0(m,k)$, for $m = 1$ and $m = 2$:*

$$\varepsilon(s, \tau \otimes \tau') = \varepsilon(s, \eta \times \eta').$$

*Proof.* Let $\tau \in \mathcal{R}^1(4,k)$. Then we know by the remark above that it is of the form $\tau_1 \otimes \tau_2$, with $\tau_1, \tau_2$ irreducible 2-dimensionals of $W_k$. Let $\eta_1, \eta_2$ be the respective supercuspidals of $\mathrm{GL}(2,k)$ given by the local Langlands correspondence. Put
$$\eta = \eta_1 \boxtimes \eta_2.$$
Suppose $\eta = \eta_1 \boxtimes \eta_2$ is not supercuspidal. Then in its isobaric sum decomposition, there must occur a supercuspidal, say $\beta$, of $\mathrm{GL}(m,k)$, for some $m < 4$. (Since $\eta$ is a representation of $W_k$, there will be no summand of Steinberg type.) Then $L(s, \beta \times \beta^\vee)$, and hence $L(s, \eta \times \beta^\vee)$ will have a pole at $s = 0$ (see [JPSS], [JS2]).

Now let $\lambda$ denote the irreducible $m$-dimensional representation of $W_k$ attached to $\beta$ by the local Langlands correspondence, which we can do by [Ku] and [He1] as $m \leq 3$. Then, by the proposition above, $L(s, \tau \otimes \lambda^\vee)$ equals $L(s, \eta \times \beta^\vee)$, and hence must have a pole at $s = 0$ as well. This cannot happen as $\tau$ is an irreducible of dimension 4, and $\lambda$ is one of dimension $< 4$. This proves that $\eta$ is supercuspidal, and we get a "reciprocity map"
$$r : \mathcal{R}^1(4,k) \to \mathcal{A}^1(4,k).$$
The $\varepsilon$-factor identity of Proposition 4.3.1, when used in conjunction with J. Chen's theorem ([Ch]), shows the injectivity of $r$. It is also easy to see that $r$ is compatible with taking duals and twisting by characters.

The surjectivity of $r$ is shown by reversing the construction of $\tau \to \eta$. The irreducibility of $\tau$ is checked by use of Proposition 4.3.2. □

Let $\pi_1, \pi_2, \pi_3, \pi_4$ be cuspidal automorphic representations of $\mathrm{GL}(2, \mathbb{A}_F)$. Define its *root number* to be
$$W(\pi_1 \times \pi_2 \times \pi_3 \times \pi_4) = \prod_v W(\pi_{1,v} \times \pi_{2,v} \times \pi_{3,v} \times \pi_{4,v}),$$
where the *local root numbers* are defined by
$$W(\pi_{1,v} \times \pi_{2,v} \times \pi_{3,v} \times \pi_{4,v}) = \varepsilon(\frac{1}{2}, \pi_{1,v} \times \pi_{2,v} \times \pi_{3,v} \times \pi_{4,v}).$$
Note that the functional equation implies that the root number is $\pm 1$ in the self-dual situation.

THEOREM 4.3.4. *Let $\pi_j$, $j \leq 4$, be as above. Assume that the central character of each $\pi_j$ is trivial. Then we have*
$$W(\pi_1 \times \pi_2 \times \pi_3 \times \pi_4) = 1.$$

When all the $\pi_j$'s are isomorphic, say to $\pi$, the root number is
$$W(\mathrm{Sym}^2(\pi) \times \mathrm{Sym}^2(\pi))W(\mathrm{Sym}^2(\pi))^2 W(1),$$
and this theorem follows in that case from Proposition 6.1 of [PR].

*Proof.* This is an immediate consequence of Proposition 4.3.1 in conjunction with the method of [PR] (see also [Ro]). We give the argument for completeness.

It suffices to show that the local root number at *any* place $v$ is 1. Clearly, for each $j \leq 4$, the associated representation $\sigma_{j,v}$ of $W'_{F_v}$ is self-dual with determinant 1 as $\pi_j$ has trivial central character. Put

$$\beta_v = \sigma_{1,v} \otimes \sigma_{2,v} \otimes \sigma_{3,v} \otimes \sigma_{4,v}.$$

By Proposition 4.3.1, we are reduced to showing that $W(\beta_v)$ (with the obvious definition) is 1. Being the tensor product of an even number of symplectic representations, $\beta_v$ is orthogonal, and by Deligne ([De]) the root number of $\tilde{\beta}_v := \beta_v \ominus 4[1]$ is given by its second Stiefel-Whitney class $w_2(\tilde{\beta}_v)$. But since $\mathrm{Sp}(2, \mathbb{C})$ is simply connected, and since the image of $\beta_v$ in $\mathrm{O}(4, \mathbb{C})$ factors through a 4-fold product of $\mathrm{Sp}(4, \mathbb{C})$, this image is simply connected. It follows that $\tilde{\beta}_v$ lifts to the Spin group, which is a 2-fold cover of $\mathrm{SO}(4, \mathbb{C})$. Consequently, $w_2(\tilde{\beta}_v)$ is trivial, from which it follows that $W(\beta_v) = W(\tilde{\beta}_v)W(1)^4$ is 1. □

4.4. *Triple product L-functions revisited.* Let $\pi_1, \pi_2, \pi_3$ be cuspidal automorphic representations of $\mathrm{GL}(2, \mathbb{A}_F)$ of respective central characters $\omega_1, \omega_2, \omega_3$. Recall from Section 3.1 the definitions of the three candidates for the triple product $L$-functions and $\varepsilon$-factors associated to $(\pi_1, \pi_2, \pi_3)$. The object of this section is to extract the following precise statement from Theorem M.

THEOREM 4.4.1. *Let $\pi_j$, $j \leq 3$, be as above. Then, at* every *place $v$,*

$$L_1(s, \pi_{1,v} \times \pi_{2,v}, \pi_{3,v}) \stackrel{(a)}{=} L(s, \pi_{1,v} \times \pi_{2,v}, \pi_{3,v}) \stackrel{(b)}{=} L_2(s, \pi_{1,v} \times \pi_{2,v}, \pi_{3,v}),$$

*and*

$$\varepsilon_1(s, \pi_{1,v} \times \pi_{2,v}, \pi_{3,v}) \stackrel{(a)}{=} \varepsilon(s, \pi_{1,v} \times \pi_{2,v}, \pi_{3,v}) \stackrel{(b)}{=} \varepsilon_2(s, \pi_{1,v} \times \pi_{2,v}, \pi_{3,v}).$$

*Proof.* Since we know this (see §3.3) at archimedean places, we may restrict our attention to finite places; fix such a place $u$. Let $p$ be the residual characteristic of $u$. Recall from the discussion in Section 3.3 that (a) (resp. (b)) is known (i) when all the representations are unramified, *and* (ii) when one of the representations is not supercuspidal (resp. when all the representations are tempered).

One has the following (at any place $v$) by the works of Shahidi:

(4.4.2) $\qquad \gamma(s, \pi_{1,v} \times \pi_{2,v} \times \pi_{3,v}) = \gamma_2(s, \pi_{1,v} \times \pi_{2,v} \times \pi_{3,v}).$

In view of the remark above, the only case to check is when one of the representations, say $\pi_{1,v}$, is a subquotient of a principal series representation defined by quasi-characters $\mu_1, \mu_2$ of $F_v^*$. Then one has by the multiplicativity of local factors ([Sh4]), the identity

$$\gamma_2(s, \pi_{1,v} \times \pi_{2,v} \times \pi_{3,v}) = \gamma_2(s, \mu_1 \otimes \pi_{2,v} \times \pi_{3,v})\gamma_2(s, \mu_\otimes \pi_{2,v} \times \pi_{3,v}).$$

The analog for $\gamma(s, \pi_{1,v} \times \pi_{2,v} \times \pi_{3,v})$ holds by the decomposition of $\sigma_{1,v} \otimes \sigma_{2,v} \otimes \sigma_{3,v}$.

Our first object is to establish a weak analog of (4.4.2) for the $\gamma_1$-factor in place of the $\gamma_2$-factor.

LEMMA 4.4.3. *Let $S(p)$ be the set of places of $F$ above $p$. Then,*

$$\prod_{v \in S(p)} \gamma(s, \pi_{1,v} \times \pi_{2,v} \times \pi_{3,v}) = \prod_{v \in S(p)} \gamma_1(s, \pi_{1,v} \times \pi_{2,v} \times \pi_{3,v}).$$

*Remark* 4.4.4. In [Ik1, p. 229, Cor. to Lemma 2.2], one finds an argument to prove Theorem 4.4.1 (a) in the special case when *each $\pi_{j,u}$ is dihedral*, being associated to a character $\chi_{j,u}$ of a quadratic extension $K_{j,u}$ of $F_u$. But in the course of the argument, there is an assertion that there is a quadratic extension $K_j$ of $F$ with local extension $K_{j,u}/F_u$, together with an idele class character $\chi$ of $K$, such that (1) $\chi$ agrees with $\chi_{j,u}$ on $K_{j,u}$ (which is fine) *and* that (2) the cuspidal automorphic representation $\eta_j$ of $GL(2, \mathbb{A}_F)$ defined by $\chi_j$ has principal series components at *all* the places $v$ outside $u$. We are unable to make sure that (2) can be achieved. (One cannot in general control the behavior of $\chi_j$ at an infinite number of places!) However, the argument can be made to work with a slight refinement as follows. We can find $K_j, \chi_j$ satisfying (1) and also (2'): $\chi_j$ is trivial at every place in $S(p) - \{u\}$; this is possible by the Grunewald-Wang theorem ([AT]). By Lemma 4.4.3 above and Lemma 2.1 of [Ik1], the desired equality then follows.

*Proof of Lemma* 4.4.3. Let $T$ be the set of finite places where at least one of the representations $\{\pi_1, \pi_2, \pi_3\}$ is ramified. Then, comparing the functional equations of $L(s, \pi_1 \times \pi_2 \times \pi_3)$ and $L_1(s, \pi_1 \times \pi_2 \times \pi_3)$, we get the equality

$$\prod_{v \in T} \gamma(s, \pi_{1,v} \times \pi_{2,v} \times \pi_{3,v}) = \prod_{v \in T} \gamma_1(s, \pi_{1,v} \times \pi_{2,v} \times \pi_{3,v}).$$

When $v$ is not in $S(p)$, its norm is a power of a prime $\ell$ distinct from $p$, and the corresponding local factors are inverses of polynomials in $\ell^{-s}$. Since the identity above holds for all $s$, it is easy to separate the contribution from the primes in $S(p)$ and those outside. □

*Proof of Theorem* 4.4.1 (*contd.*). Now we need the following:

PROPOSITION 4.4.5. *Let $v$ be a finite place, and let $L_*$ denote $L, L_1$ or $L_2$. Then, for any quasi-character $\nu$ which is sufficiently ramified,*

$$L_*(s, \pi_{1,v} \times \pi_{2,v} \times \pi_{3,v} \otimes \nu) = 1.$$

*Proof of Proposition* 4.4.5. First we consider when $L_* = L$. Since by definition $L(s, \pi_{1,v} \times \pi_{2,v} \times \pi_{3,v} \otimes \nu)$ is the same as $L(s, \sigma_{1,v} \otimes \sigma_{2,v} \otimes \sigma_{3,v} \otimes \nu)$, it suffices, by the additivity of the $L$-factor, to show that, for any irreducible summand $\tau$ of $\sigma_{1,v} \otimes \sigma_{2,v} \otimes \sigma_{3,v}$, viewed as a representation of $W_{F_v} \times \mathrm{SL}(2, \mathbb{C})$, we have $L(s, \tau \otimes \nu)$ is 1 for $\nu$ of large enough conductor. Any irreducible $\tau$ is of the form $\beta \otimes S^k$, where $\beta$ is an irreducible of $W_{F_v}$ and $S^j$ denotes the symmetric $j^{\mathrm{th}}$ power representation of $\mathrm{SL}(2, \mathbb{C})$. When $\beta$ is 1-dimensional, $L(s, \tau \otimes \nu)$ equals $L(s, \beta\nu|.|^{j/2})$, which is trivial (by Tate's thesis) if $\nu$ is sufficiently ramified. If $\beta$ is higher dimensional, then $L(s, \tau \otimes \nu)$ is a product, for suitable half-integers $t$, of factors of the form $L(s, \beta \otimes \nu|.|^t)$, which are 1 for *any* $\nu$ (because $\beta$ is irreducible of dimension $> 1$).

Next consider the case $L_* = L_i$, with $i = 1, 2$. We may assume that at least one of the $\pi_{j,v}$, say $\pi_{3,v}$, is tempered. Suppose $\pi_{1,v}$ is not supercuspidal. Then by using Proposition 3.3.6 for $L_1$, and (4.4.2) for $L_2$ together with the temperedness of $\pi_{3,v}$, we deduce that (for any $\nu$)

$$\frac{L_1(s, \pi_{1,v} \times \pi_{2,v} \times \pi_{3,v} \otimes \nu)}{L(s, \pi_{1,v} \times \pi_{2,v} \times \pi_{3,v} \otimes \nu)} \quad \text{is} \quad \text{entire.}$$

The assertions (a) and (b) of the theorem about $L$-factors follows in this case by the triviality of the denominator for sufficiently ramified $\nu$. The identities for the $\varepsilon$-factors also follows because the respective $\gamma$-factors are the same.

So we may take all the $\pi_{j,v}$ to be supercuspidal. We will now use a slight variant of the inductive argument utilized in the proof of Section 3.7. Let $\ell(\pi_{1,v})$, $\ell(\pi_1)$, $p(v)$, $p = p(\pi_1)$, $r = r(\pi_1) = (\ell(\pi_1), p(\pi_1))$ be as in Section 3.7, together with the (lexicographic) ordering $<$ introduced there. We will use induction on $r$. If $r = (0, 1)$, we are done, and so take that $(0, 1) < r$ and assume by induction that the proposition holds for $L_1$ at all places for any triple $(\pi, \pi', \pi'')$ of cuspidal automorphic representations of $\mathrm{GL}(2)$ over a number field if $r(\pi) < r$. Fix at every place $v$ in $S''$, a character $\chi_v$ cutting out a $K(v)$ of degree $p$ such that the base change of $\pi_{1,v}$ to $K(v)$ has $\ell$ equal to $\ell(\pi_1)/p$. Choose, using the Grunewald-Wang theorem, a global character $\chi$ of order $p$, cutting out a cyclic extension $K/F$ with local extensions $K(v)/F_v$. Then $r(\pi_{1,K})$ is less than $r$ and we may use induction in conjunction with the functional equations of $L_i(s, \pi_1 \times \pi_2 \times \pi_3 \otimes \nu)$ and $L_i(s, \pi_{1,K} \times \pi_{2,K} \times \pi_{3,K} \otimes \nu)$ to conclude that, for any idele class character $\nu$ of $F$ which is sufficiently ramified at each of the places in $S''$, we have (for $i = 1, 2$)

$$\prod_{v \in S''} \gamma_i(s, \pi_{1,v} \times \pi_{2,v} \times \pi_{3,v} \otimes \nu) = 1.$$

As remarked earlier, we may assume that, for every $j \leq 3$ and every $v \in S''$, the local representation $\pi_{j,v}$ is supercuspidal. We may then apply Lemma 2.1 of [Ik1] to conclude that $L_i(s, \pi_{1,v} \times \pi_{2,v} \times \pi_{3,v} \otimes \nu)$ has no pole in common with $L_i(1-s, \pi_{1,v}^\vee \times \pi_{2,v}^\vee \times \pi_{3,v}^\vee \otimes \nu^{-1})$. The proposition follows for $L_i$. □

*Remark* 4.4.6. One can also prove Proposition 4.4.5 for $L_2$ by arguing as in the proof of Lemma 5 of [PR]. Moreover, Shahidi has informed us that he can now prove a very general statement of this sort for $L_2$.

*Proof of Theorem* 4.4.1 (*contd.*). Fix any finite place $u$ of residual characteristic $p$, and let $S(p)$ be as in Lemma 4.4.3. Choose an idele class character $\nu$ of $F$ such that (i) $\nu_u$ is trivial, and (ii) $\nu_v$ is highly ramified at each $v \neq u$ in $S(p)$. Then, comparing the functional equations of $L_1(s, \pi_1 \times \pi_2 \times \pi_3 \otimes \nu)$ and $L(s, \pi_1 \times \pi_2 \times \pi_3 \otimes \nu)$ and using (4.4.2), Lemma 4.4.3 and Proposition 4.4, we get

$$(4.4.7) \qquad \gamma_1(s, \pi_{1,u} \times \pi_{2,u} \times \pi_{3,u} \otimes \nu) = \gamma(s, \pi_{1,u} \times \pi_{2,u} \times \pi_{3,u} \otimes \nu).$$

Again, we can reduce to the case when each $\pi_{j,u}$ is supercuspidal and so, applying Lemma 2.1 of [Ik1] and the obvious analog for $L_2$, we deduce the assertions (a) and (b) of the theorem for the $L$-factors. The identities for the $\varepsilon$-factors also follow by (4.4.7). □

4.5. *The Tate conjecture for* 4-*fold products of modular curves.* Let $Y$ denote the moduli scheme over $\mathbb{Q}$ parametrizing elliptic curves $E$ with level infinity structure, i.e., equipped with an isomorphism of the Tate module $T_f(E) = \lim_n E[n]$ with $\hat{\mathbb{Z}}^2$. It comes with a natural right action of $\mathrm{GL}(2, \mathbb{A}^f)$, and a smooth compactification $X$ ([KM]). For every compact, open subgroup $K$ of $\mathrm{GL}(2, \mathbb{A}^f)$, let $Y_K$ (resp. $X_K = Y_K \cup Y_K^\infty$) denote $Y/K$ (resp. $X/K$). Then $X_K$ is a smooth projective curve over $\mathbb{Q}$, and one has the identifications

$$Y_K(\mathbb{C}) = \mathrm{GL}(2, \mathbb{Q}) \backslash \mathcal{H}^\pm / K,$$

where $\mathcal{H}^\pm$ denotes $\mathbb{C} - \mathbb{R}$, and

$$Y_K^\infty = \mathbb{P}^1(\hat{\mathbb{Z}})/K.$$

A cofinal system of compact open subgroups is provided by $\{K_1(N) | N \geq 1\}$, where

$$K_1(N) = \{k = (k_{ij}) \in \mathrm{GL}(2, \hat{\mathbb{Z}}) \,|\, k_{21}, k_{22} - 1 \in N\hat{\mathbb{Z}}\}.$$

For every $N \geq 1$, we will also be interested in

$$K_0(N) = \{k = (k_{ij}) \in \mathrm{GL}(2, \hat{\mathbb{Z}}) \,|\, k_{21} \in N\hat{\mathbb{Z}}\}.$$

We will write $X_1(N), X_0(N)$ instead of $X_{K_1(N)}, X_{K_0(N)}$ respectively.

Let $V$ be a smooth projective variety over a number field $k$. For every $j \geq 0$, denote by $C^j(V/k)$ the $\mathbb{Q}$-vector space spanned by codimension $j$ algebraic cycles on $V$ over $k$ modulo homological equivalence. Then, for every prime $\ell$, one has an $\ell$-adic cycle class map

$$C^j(V/k) \longrightarrow H^{2j}_{\mathrm{et}}(V_{\overline{k}}, \mathbb{Q}_\ell)(j),$$

which is injective and lands in the group of codimension $j$ *Tate cycles* over $k$, namely

$$\mathrm{Ta}_\ell^j(V/k) \;=\; H^{2j}_{\mathrm{et}}(V_{\overline{k}}, \mathbb{Q}_\ell)(j)^{\mathrm{Gal}(\overline{k}/k)}.$$

Let $S$ be a finite set of places containing the ramified places for the $\mathrm{Gal}(\overline{\mathbb{Q}}/k)$-module $H^{2j}_{\mathrm{et}}(V_{\overline{k}}, \mathbb{Q}_\ell)$, the archimedean places, and also the places above $\ell$. The associated $L$-function is

$$L^{(2j)}(s, V/k) \;=\; \prod_{v:\,\mathrm{finite}} L_v^{(2j)}(s, V/k),$$

where

$$L_v^{(2j)}(s, V/k) \;=\; \det(I - \mathrm{Fr}_v T | H^{2j}_{\mathrm{et}}(V_{\overline{k}}, \mathbb{Q}_\ell)^{I_v})^{-1}|_{T=(Nv)^{-s}}.$$

Here $\mathrm{Fr}_v$ denotes the geometric Frobenius at $v$ and $I_v$ the inertia group at $v$. In computing $L_v^{(2j)}(s, V/k)$, one takes $\ell$ to be prime to $v$. For each $\ell$, let $S = S(\ell)$ be the finite set of places consisting of the ramified places for the $\mathrm{Gal}(\overline{\mathbb{Q}}/k)$-module $H^{2j}_{\mathrm{et}}(V_{\overline{k}}, \mathbb{Q}_\ell)$, together with the places above $\ell$. One knows by Deligne's proof of the Weil conjectures that the incomplete $L$-function $L^{(2j),S}(s, V/k)$, obtained from $L^{(2j)}(s, V/k)$ by deleting the factors over $S$, converges absolutely in $\Re(s) > j+1$. One also has a definition of $L_v^{(2j)}(s, V/k)$ for any archimedean place $v$ (see [Ra3, §.5.5], for example). Write $L_\infty^{(2j)}(s, V/k)$ for the product $\prod_{v|\infty} L_v^{(2j)}(s, V/k)$, and put

$$\tilde{L}^{(2j)}(s, V/k) \;:=\; L_\infty^{(2j)}(s, V/k) L^{(2j)}(s, V/k).$$

One also has epsilon factors $\varepsilon_v^{(2j)}(s, V/k)$ for each $v$, and the global one

$$\varepsilon^{(2j)}(s, V/k) \;=\; \prod_{\mathrm{all}\,v} \varepsilon_v^{(2j)}(s, V/k),$$

which is an invertible holomorphic function.

THEOREM 4.5.1.  *Fix positive integers $N_j$, $1 \le j \le 4$, and put*

$$V \;=\; X_1 \times X_2 \times X_3 \times X_4,$$

*with $X_j$ being $X_0(N_j)$ or $X_1(N_j)$ for each $j \le 4$. Let $k$ be a finite solvable, normal extension of $\mathbb{Q}$, and let $\ell$ be a prime. Then the following hold*:

(a) *$L^{(4)}(s, V/k)$ admits a meromorphic continuation to the whole $s$-plane, has a convergent Euler product in $\{\mathrm{Re}(s) > 3\}$, and satisfies the exact functional equation*:

$$\tilde{L}^{(4)}(s, V/k) \;=\; \varepsilon^{(4)}(s, V/k) \tilde{L}^{(4)}(5-s, V/k).$$

(b) *$L^{(4)}(s, V/k)$ has no pole anywhere except possibly at the "edge" point $s=3$.*

(c) *If we let $-\mathrm{ord}_{s=3}$ stand for the order of the pole at $s = 3$,*

$$-\mathrm{ord}_{s=3} L^{(4)}(s, V/k) = \dim_{\mathbb{Q}_\ell} \mathrm{Ta}^2_\ell(V/k).$$

(d) *Suppose some $X_j$ is $X_0(N_j)$ with $N_j$ square-free. Then*

$$-\mathrm{ord}_{s=3} L^{(4),S}(s, V/k) = \dim_{\mathbb{Q}} C^2(V/k).$$

*In this case, for any number field $k$, $\dim_{\mathbb{Q}_\ell} \mathrm{Ta}^2_\ell(V/k)$ and $\dim_{\mathbb{Q}} C^2(V/k)$ are equal.*

Parts (a) and (b) are as predicted by the Hasse-Weil hypothesis, while parts (c) and (d) verify certain conjectures of Tate (see [Ta2], [Ra3, Chap. 5]).

*Proof.* First we recall some well-known basic facts about the cohomology of the modular curves $X_1(N)$. For $j \in \{0, 1, 2\}$, put (for any prime $\ell$)

$$W^j_\ell(N) := H^j_{\mathrm{et}}(X_1(N)_{\overline{\mathbb{Q}}}, \mathbb{Q}_\ell).$$

For any field extension $E$ of $\mathbb{Q}$, denote by $\mathcal{H}_E(N)$ the $E$-algebra of Hecke correspondences at level $N$. This algebra is semisimple and acts on $W^j_\ell(N) \otimes_{\mathbb{Q}} E$, commuting with the Galois action. If $\pi = \pi_\infty \otimes \pi_f$ is an irreducible automorphic representation of $\mathrm{GL}(2, \mathbb{A}_\mathbb{Q})$, with $\pi_f$ rational over $\overline{\mathbb{Q}}$, then the $K_1(N)$-invariants in the space of $\pi_f$ over $\overline{\mathbb{Q}}$ is naturally an $\mathcal{H}_{\overline{\mathbb{Q}}}(N)$-module. Conversely, any irreducible $\mathcal{H}_{\overline{\mathbb{Q}}}(N)$-module occurring in $W^j_\ell(N)$ is a summand of some $\pi_f^{K_1(N)}$. Denote by $\mathrm{Coh}^j(N)$ the set of irreducible automorphic representations $\pi$ contributing to the cohomology in degree $j$ with $\pi_f^{K_1(N)}$ nonzero. Then one knows that $\pi$ is 1-dimensional for $j = 0, 2$, and that for $j = 1$ it is cuspidal with $\pi_\infty$ in the lowest discrete series, corresponding, in the classical language, to a holomorphic newform of weight 2 and level $N$. Decomposing according to the Hecke algebra, one gets

$$(4.5.2) \qquad W^j_\ell(N) \otimes_{\mathbb{Q}_\ell} \overline{\mathbb{Q}}_\ell \simeq \oplus_{\pi \in \mathrm{Coh}^j(N)} W^j_\ell(\pi)^{m(\pi_f, N)},$$

where $W^j_\ell(\pi)$ is an irreducible $\overline{\mathbb{Q}}_\ell$-representation of $\mathrm{Gal}(\overline{\mathbb{Q}}/\mathbb{Q})$ of dimension 2 (resp. 1) when $j = 1$ (resp. $j = 0, 2$), and $m(\pi_f, N)$ is the dimension of the $K_1(N)$-invariants in the space of $\pi_f$. When a $\pi$ occurs in degree 1, then any of its Galois conjugates $\pi^\tau := \pi_\infty \otimes \pi_f^\tau$ also occurs. One knows the following at *any* place $v$ by the Eichler-Shimura theory and its refinement due to Igusa, Deligne, Langlands and Carayol([Ca]):

$$(4.5.3) \qquad L_v(s, W^0_\ell(\pi)) = L(s, \pi_v) = L_v(s+1, W^2_\ell(\pi)),$$

and

$$L_v(s, W^1_\ell(\pi)) = L(s - 1/2, \pi_v).$$

In fact this also holds for 1-dimensional twists. The analogous identities hold for the $\varepsilon$-factors as well, implying in particular the identification of the conductors.

Now let us turn our attention to the variety $V$ at hand. In the following $i$ will denote any 4-tuple $(i_1, i_2, i_3, i_4)$ of integers in $\{0, 1, 2\}$ with $\sum_{j=1}^{4} i_j = 4$, and for each such $i$, $\pi$ will denote $(\pi_1, \pi_2, \pi_3, \pi_4)$ with $\pi_j \in \text{Coh}^{i_j}(N_j)$ for each $j$. Let $m(i, \pi)$ signify $\prod_{j=1}^{4} m(\pi_{j,f}, N_j)$. Applying the Künneth formula in conjunction with the decomposition (4.5.2), we get

$$(4.5.4) \qquad H^4(V_{\overline{\mathbb{Q}}}, \overline{\mathbb{Q}}_\ell) \simeq \oplus_i \oplus_\pi H^i_\ell(\pi)^{m(i,\pi)},$$

where

$$H^i_\ell(\pi) \simeq W^{i_1}_\ell(\pi_1) \otimes W^{i_2}_\ell(\pi_2) \otimes W^{i_3}_\ell(\pi_3) \otimes W^{i_4}_\ell(\pi_4).$$

Applying Proposition 4.3.1 in conjunction with (4.5.3), we get moreover,

$$(4.5.5) \qquad L(s, H^i_\ell(\pi)) = L(s - 2, \pi_1 \times \pi_2 \times \pi_3 \times \pi_4),$$

and similarly for the $\varepsilon$-factors.

Since the cycle class maps are functorial for the action of correspondences, we can split the problem and prove the assertions of the theorem for the image of the codimension 2 cycles in $H^i_\ell(\pi)$, for every $i$ and $\pi$. Let $\text{Ta}^i_\ell(\pi)_k$ denote the space of Tate classes over $k$ in $H^i_\ell(\pi)$.

Suppose two of the indices, say $i_3, i_4$, are zero. Then $H^i_\ell(\pi)$ is

$$W^2_\ell(\pi_1) \otimes W^2_\ell(\pi_2) \otimes W^0_\ell(\pi_3) \otimes W^0_\ell(\pi_4),$$

which is 1-dimensional, in fact of the form $\mu \otimes \overline{\mathbb{Q}}_\ell(-2)$ with $\mu$ of finite order. Then for *any* number field $k$, we have $L(s, H^i_\ell(\pi)) = L(s - 2, \mu)$. The assertions (a), (b) follow by Hecke. The order of pole at $s = 3$ is 1 or 0 depending on whether or not $\mu$ is trivial when restricted to $\text{Gal}(\overline{Q}/k)$. When $\mu$ is trivial, the corresponding algebraic cycle is given by the intersection of the the $\pi$-components of the threefolds $\{P\} \times X_1(N_2) \times X_1(N_3) \times X_1(N_4)$ and $X_1(N_1) \times \{Q\} \times X_1(N_3) \times X_1(N_4)$, for points $P, Q$ on $X_1(N_1), X_1(N_2)$ respectively. So (d) also follows in this case.

Now, we may assume that $i = (1, 1, 1, 1)$ and consequently that each $\pi_j$ is cuspidal. In view of (4.5.5), parts (a) and (b) hold by Theorem 4.2.1. Note that the determinant of each $W^1_\ell(\pi_j)$ is necessarily of the form $\omega_j \otimes \overline{\mathbb{Q}}_\ell(-1)$, where $\omega_j$ is of finite order corresponding to the central character of $\pi_j$ by class field theory.

Suppose all the $\pi_j$ are dihedral, defined by a Hecke character $\chi_j$ of a (necessarily imaginary) quadratic field $K_j$. Then $H^i_\ell(\pi)$ is the $\overline{\mathbb{Q}}_\ell$-realization of a CM motive (in the category of pure motives for absolute Hodge cycles), and the assertion (c) is well-known (see [DMOS]). The conductor $M_j$, say, of each $\pi_j$ is the conductor of an abelian variety factor $A(\pi_j)$ over $\mathbb{Q}$, up to isogeny, of the Jacobian of the corresponding modular curve. One knows by a theorem of Fontaine ([Fon]) that there is no abelian variety over $\mathbb{Q}$ of good reduction everywhere. On the other hand, since $\pi_j$ is automorphically induced by $\chi_j$, it has no local component which is special. Then by the description of the conductor of $\pi_j$ in [Ge], for example, we see that $M_j$, and hence $N_j$ which it divides, cannot be square-free. So the hypothesis of part (d) precludes the case when *all* the $\pi_j$ are dihedral.

Thus we may, and will, assume that some $\pi_j$, say $\pi_1$, is *nondihedral*.

LEMMA 4.5.6.  *Let $i = (1,1,1,1)$ and $\pi = (\pi_1, \pi_2, \pi_3, \pi_4)$, with each $\pi_j$ cuspidal automorphic of weight 2 at infinity and let $\pi_1$ be nondihedral. Fix any number field $k$. Then the dimension of $\mathrm{Ta}_\ell^i(\pi)_k$ is at most 2. It is nonzero if and only if, possibly after $\{\pi_2, \pi_3, \pi_4\}$ is renumbered, there is a character $\mu$ of $\mathrm{Gal}(\overline{\mathbb{Q}}/k)$ such that*

(i) *as $\mathrm{Gal}(\overline{\mathbb{Q}}/k)$-modules,*

$$W_\ell^1(\pi_2) \simeq W_\ell^1(\pi_1) \otimes \mu$$

*and*

(ii) *$(\mu\omega_1)^{-1}$ occurs in $W_\ell^1(\pi_3) \otimes W_\ell^1(\pi_4)$ over $k$.*

*Moreover, when this happens, the dimension of $\mathrm{Ta}_\ell^i(\pi)$ is 1 unless one of the following situations holds*:

(a) *Over $k$, $W_\ell^1(\pi_3))$ is irreducible and nondihedral, with its symmetric square being isomorphic to $S^2(W_\ell^1(\pi_1)) \otimes \omega_3 \omega_1^{-1}$;*

(b) *$W_\ell^1(\pi_3)$ and $W_\ell^1(\pi_4)$ are both reducible over $k$, and $\dim \mathrm{Hom}(\mu\omega_1^{-1}, W_\ell^1(\pi_3) \otimes W_\ell^1(\pi_4)$ is 2.*

*Proof.* Since $\pi_1$ is nondihedral, we know by a result of Ribet ([Ri1]) that $W_\ell^1(\pi_1)$ is irreducible upon restriction to any open subgroup of $\mathrm{Gal}(\overline{\mathbb{Q}}/\mathbb{Q})$. It follows easily that there are no Tate classes over $k$ unless (i) holds after a possible renumbering. There will be a Tate cycle if and only if (ii) also holds. Now suppose there is a Tate class.

If $\pi_3$ is irreducible, (ii) will hold if and only if $\pi_4$ is also irreducible and

(ii') $$W_\ell^1(\pi_4) \simeq W_\ell^1(\pi_3) \otimes (\mu\omega_1\omega_3)^{-1}.$$

We claim that, since $\pi_1$ is nondihedral, $S^2(W_\ell^1(\pi_1))$ must be irreducible under restriction to any open subgroup. Indeed, if the claim is false, then over a number field the symmetric square would admit a 1-dimensional and would force $W_\ell^1(\pi_1)$ to be dihedral (over the extension field), contradicting the openness of the image of $\mathrm{Gal}(\overline{\mathbb{Q}}/\mathbb{Q})$ in $\mathrm{GL}(W_\ell^1(\pi_1))$. Hence the only way there could be more than one Tate class over $k$ is to have (a), and this is not possible if $\pi_3$ is dihedral. Suppose $\pi_3$ is nondihedral. Then each $\pi_j$ is nondihedral and $S^2(W_\ell^1(\pi_j))$ is irreducible (by the above argument) for *any* $j$. Then the number of (independent) such classes can evidently not be more than 2.

Finally, consider the case when (i), (ii), hold, *and* $\pi_3, \pi_4$ are both reducible. Looking at the Hodge-Tate types we see that the dimension of $\mathrm{Ta}_\ell^i(\pi)$ can never exceed 2. It is then also clear that to have dimension 2, we need (b).

*Remark* 4.5.7. The reason for requiring in Theorem 4.5.1 that some $X_j$ be $X_0(N_j)$ with $N_j$ square-free is just to rule out the case when all the $\pi_j$ are dihedral, in which case the number of Tate classes can at the extreme be six, and only four of those can be shown to be algebraic by the argument below. When one of the $\pi_j$ is nondihedral, the numbers of Tate and algebraic cycles in $H_\ell^i(\pi)$ coincide over every number field, not necessarily solvable. The reason for restricting to solvable extensions is to be able to appeal to base change in conjunction with Theorem 4.2.1 and deduce the standard properties of the $L$-function.

*Proof of Theorem* 4.5.1 (*contd.*). Let $k$ be a number field. We need to match the number of Tate classes with the order of pole. By using the standard results on Hecke $L$-functions, we may assume that one of the $\pi_j$, say $\pi_1$, is nondihedral. In view of (4.5.5), part (c) of the theorem will follow once we establish the following analog of Lemma 4.5.6:

LEMMA 4.5.8. *Let $i = (1,1,1,1)$ and $\pi = (\pi_1, \pi_2, \pi_3, \pi_4)$, with each $\pi_j$ cuspidal automorphic of weight $2$ at infinity and $\pi_1$ nondihedral. Fix any number field $k$, and denote (for all $j$) by $\pi_{j,k}$ the base change of $\pi$ to $k$. Then the order of pole, $m(\pi_k)$ say, of $L(s, \pi_k) := L(s, \pi_1 \times \pi_2 \times \pi_3 \times \pi_4)$ at $s = 1$ is at most $2$. It is nonzero if and only if, when $\{\pi_2, \pi_3, \pi_4\}$ is possibly renumbered, there exists an idele class character $\mu$ of $k$ such that*

(i) $\pi_{2,k} \simeq \pi_{1,k} \otimes \mu$, *and*

(ii) $(\mu_k \omega_{1,k})^{-1}$ *occurs in the isobaric sum decomposition of* $\pi_{3,k} \boxtimes \pi_{4,k}$.

*Moreover, when this happens, $m(\pi_k)$ is $1$ unless one of the following situations holds*:

(a) $\pi_{3,k}$ *is cuspidal and nondihedral, and its symmetric square is isomorphic to* $S^2((\pi_1)) \otimes \omega_{3,k} \omega_{1,k}^{-1}$;

(b) $\pi_{3,k}$ *and* $\pi_{4,k}$ *are both Eisensteinian over $k$, and the multiplicity of* $(\mu_k \omega_{1,k})^{-1}$ *in* $\pi_{3,k} \boxtimes \pi_{4,k}$ *is* $2$.

*Proof of Lemma* 4.5.8. The conditions (i) and (ii) are easy to verify if every $\pi_j$ is twist equivalent to $\pi_1$. Now suppose not. Then, after a renumbering, we may assume that $\pi_{1,k}\boxtimes\pi_{3,k}$ is cuspidal. Using Theorem M we can write $L(s,\pi_k)$ as $L(s,(\pi_1\boxtimes\pi_3)\times(\pi_2\boxtimes\pi_4))$ and by [JS2], it has a pole at $s=1$ if and only if we have

$$\pi_{2,k}\boxtimes\pi_{4,k}\simeq(\pi_{1,k}\boxtimes\pi_{3,k})\simeq\pi_{1,k}\boxtimes\pi_{3,k}\otimes(\omega_{1,k}\omega_{3,k})^{-1}.$$

Up to a renumbering, this can happen if and only if (i) and (ii) hold.

If $\pi_{3,k}$ is cuspidal, then (ii) will hold if and only if $\pi_{4,k}$ is also cuspidal and

(ii') $$\pi_{4,k}\simeq\pi_{3,k}\otimes(\mu_k\omega_{1,k}\omega_{3,k})^{-1}.$$

We claim that $S^2(\pi_{1,k})$ is cuspidal. Suppose not. Then by [GJ], $\pi_{1,k}$ will be dihedral. This would imply, by (4.5.3) and the functoriality of base change, that $W_\ell^1(\pi_1)$ is dihedral when restricted to $\mathrm{Gal}(\overline{\mathbb{Q}}/k)$, contradicting the openness of the image of Galois ([Ri1]). Hence the claim holds. Applying [JS2], we then see that $m(\pi_k)$ could be more than 1 if and only if we have (a), and this is evidently not possible if $\pi_3$ is dihedral. Then each $\pi_j$ is nondihedral and $m(\pi_k)$ cannot possibly be more than 2. when $\pi_3,\pi_4$ are both reducible, we see by looking at the infinity type that $m(\pi_k)$ can never exceed 2, and that it can be 2 if and only if (b) holds. □

It remains to show (d) of Theorem 4.5.1. Since we have assumed that some $X_j$ is $X_0(N_j)$ with $N_j$ square-free, say for $j=1$, $\pi_1$ will not be dihedral; so we are in the situation of Lemma 4.5.8. In view of (c), it suffices to show, over any number field $k$, that $\mathrm{Ta}_\ell^i(\pi)$ consists of algebraic classes. Suppose there is a nonzero Tate class. Then (i) holds and this is represented by a twisting correspondence $R_\mu$ on $X_1(N_1)_k\times X_1(N_2)_k$ (see [Ri2], [Mu]). If $\pi_{3,k}$ is also cuspidal, then (ii') holds and this is represented by another twisting correspondence $R_\nu$ on $X_1(N_3)\times X_1(N_4)$, with $\nu=(\mu\omega_1\omega_3)^{-1}$. We get a codimension 2 algebraic cycle on $V_k$, whose $\pi_k$-component has cycle class in $\mathrm{Ta}_\ell^i(\pi)$, by taking the intersection $\gamma$ of $R_\mu\times X_1(N_3)\times X_1(N_4)$ and $X_1(N_1)\times X_1(N_2)\times R_\nu$. Now suppose there are two independent Tate classes. If $\pi_{3,k}$ is cuspidal, then we would have the symmetric square twisting equivalence of (a), and this will give us a twisting correspondence $R_\lambda$ directly on the 4-fold $V$, with $\lambda=\omega_3\omega^{-1}$. This cycle is evidently not a multiple of $\gamma$. Thus, we are done in this case. If we are in the situation of (b), the abelian varieties $A(\pi_3)$ and $A(\pi_4)$ corresponding to $\pi_3$ and $\pi_4$ will necessarily be of CM type, and we will have two independent Tate classes in $\mathrm{Hom}(A(\pi_3),A(\pi_4))$ over $k$ which are algebraic by Faltings; they correspond to divisors $\delta_1,\delta_2$ on $X_1(N_3)\times X_1(N_4)$ rational over $k$. If we take $\xi_j$, for $j=1,2$ to be $R_\mu\times\delta_j$, they define the needed cycles on $V$. □

California Institute of Technology, Pasadena, CA
*E-mail address*: dinakar@its.caltech.edu